\documentclass[10pt]{article}

\usepackage{natbib}                     
\usepackage{amsmath, amssymb, amsthm} 

\usepackage{graphicx}                   
\usepackage{fullpage}                   
\usepackage{hyperref}                   
\usepackage{enumitem}                   
\usepackage{bm}                         
\usepackage{authblk}                    

\usepackage{multirow}
\usepackage{algorithm}
\usepackage{algorithmic}
\usepackage{booktabs}
\usepackage{xcolor}
\usepackage{caption}
\usepackage{tikz}
\usepackage{pgfplots}
\usepgfplotslibrary{groupplots}
\usepgfplotslibrary{fillbetween}
\pgfplotsset{compat=1.18}
\pgfplotsset{
    myplotstyle/.style={
    ylabel style={align=center, font=\bfseries\boldmath},
    xlabel style={align=center, font=\bfseries\boldmath},
    x tick label style={font=\bfseries\boldmath},
    y tick label style={font=\bfseries\boldmath},
    scaled ticks=false,
    every axis plot/.append style={thick},
    },
}
\ifpdf
  \DeclareGraphicsExtensions{.eps,.pdf,.png,.jpg}
\else
  \DeclareGraphicsExtensions{.eps}
\fi
\definecolor{darkblue}{rgb}{0.0, 0.0, 0.5}  

\hypersetup{
    colorlinks=true,
    linkcolor=darkblue,
    anchorcolor=darkblue,
    citecolor=darkblue,
    filecolor=darkblue,
    menucolor=darkblue,
    urlcolor=darkblue
}

\usepackage[nameinlink,noabbrev]{cleveref}

\usepackage{caption}

\newtheorem{theorem}{Theorem}
\newtheorem{lemma}{Lemma}
\newtheorem{proposition}{Proposition}
\newtheorem{corollary}{Corollary}
\theoremstyle{definition}

\newtheorem{remark}{Remark}

\title{A Unified Model for High-Resolution ODEs: \\ New Insights on Accelerated Methods\footnote{Portions of this work were previously published at NeurIPS 2023 as ``A Variational Perspective on High-Resolution ODEs." That work focused on deriving high-resolution ODEs for Nesterov's accelerated gradient method using the forced Euler-Lagrange equation. The present paper offers a substantial extension by introducing a broader and more general framework that unifies a wide range of high-resolution ODEs as special cases of a newly introduced model, providing new theoretical insights and improved convergence guarantees for various momentum-based methods. This work excludes discussions on stochastic extensions presented in the NeurIPS version and introduces novel results for the triple momentum HR-ODE and quasi-hyperbolic momentum.}}
\author[1]{Hoomaan Maskan}
\author[2]{Konstantinos C. Zygalakis}
\author[3]{Armin Eftekhari}
\author[4]{Alp Yurtsever}
\affil[1,3,4]{Department of Mathematics and Mathematical Statistics, Umeå University, Sweden}
\affil[2]{School of Mathematics, University of Edinburgh, Scotland} 
\affil[1]{\texttt{\{hoomaan.maskan\}@umu.se}}
\affil[2]{\texttt{\{k.zygalakis\}@ed.ac.uk}}
\affil[3]{\texttt{\{armin.eftekhari\}@gmail.com}}
\affil[4]{\texttt{\{alp.yurtsever\}@umu.se}}
\date{}

\begin{document}

\maketitle

\begin{abstract}
Recent work on high-resolution ordinary differential equations (HR-ODEs) captures fine nuances among different momentum-based optimization methods, leading to accurate theoretical insights. 
However, these HR-ODEs often appear disconnected, each targeting a specific algorithm and derived with different assumptions and techniques. 
We present a unifying framework by showing that these diverse HR-ODEs emerge as special cases of a general HR-ODE derived using the Forced Euler-Lagrange equation. 
Discretizing this model recovers a wide range of optimization algorithms through different parameter choices. 
Using integral quadratic constraints, we also introduce a general Lyapunov function to analyze the convergence of the proposed HR-ODE and its discretizations, achieving significant improvements across various cases, including new guarantees for the triple momentum method's HR-ODE and the quasi-hyperbolic momentum method, as well as faster gradient norm minimization rates for Nesterov's accelerated gradient algorithm, among other advances. 
\end{abstract}

\section{Introduction}
\label{sec:intro}
In this work, we focus on the unconstrained minimization problem
\begin{align} \label{minf}
\min_{x \in \mathbb{R}^{d}} f(x),
\end{align}
where $f:\mathbb R^d \rightarrow \mathbb{R}$ either belongs to the class $\mathcal{F}_{L}$, the set of $L$-smooth convex functions with $L$ being the Lipschitz constant of the gradient, or the class $\mathcal{F}_{\mu,L}$ of $L$-smooth $\mu$-strongly convex functions. 
We denote a solution to this problem by $x^*$. 

Many classical methods in convex optimization can be viewed through the lens of a sequential approximation strategy, where a challenging optimization problem is tackled by breaking it down into a series of simpler sub-problems that can be solved efficiently, often in closed form. 
These methods compute a step direction by minimizing a surrogate function that approximates the objective in the vicinity of the current estimate. 
The simplest such method is gradient descent 
\begin{align} \label{eq:GD}
    x_{k+1}=x_{k}-s\nabla f(x_{k}),
\end{align}
where $s>0$ is the step-size, that determines the next point by minimizing an isotropic quadratic upper bound. 
In a similar fashion, Newton's method minimizes the second-order Taylor approximation, and the Frank-Wolfe algorithm moves toward the minimizer of the first-order Taylor approximation over the feasible set. 

However, momentum-based accelerated methods do not strictly follow this sequential approximation framework. 
Instead of computing steps simply by minimizing local surrogates, they incorporate momentum terms that leverage past iterates to accelerate convergence. 
Two of the earliest algorithms in this category are the Heavy-Ball (HB) \citep{polyak1964some} and Nesterov's Accelerated Gradient (NAG) \citep{Nesterov1983AMF} methods. 
In particular, the Nesterov family of accelerated gradient methods are of the form
\begin{align} \label{eqn_Nest_alg}\tag{NAG}
\begin{aligned}
    x_{k+1}&=y_k -s\nabla f(y_k), \\
    y_{k+1}&= x_{k+1}+\beta_{k} (x_{k+1}-x_{k}).
\end{aligned}
\end{align}
where $s \leq 1/L$ and $\beta_{k}$ changes depending on the properties of the function $f(x)$. 
This algorithm is close in spirit to gradient descent, but the gradient is evaluated at a convex combination between $x_{k}$ and $x_{k+1}$. 
The HB method has the same structure as NAG in the strongly convex case with the only differences being the coefficient $\beta_{k}$, and the fact that the gradient is evaluated at $x_{k}$.  

In recent years, there has been a significant effort to better understand the essence of acceleration by analyzing the NAG algorithm from different perspectives. 
One particular line of research, initiated by \citet{ALVAREZ2002747} and later revived by \citet{JMLR:v17:15-084}, focuses on the continuous-time analysis of acceleration in the limit of an infinitesimal step-size. 
In the convex case, it can be shown that the discrete iterates $x_k$ converge to a continuous trajectory $X(t)$ as the step-size $s \to 0$, where $X(t)$ satisfies the following second-order ODE:
\begin{align}
    \ddot{X}_{t}+\frac{3}{t} \dot{X}_{t}+\nabla f(X_{t})=0,
\end{align}
while in the strongly convex case, one has the following second-order ODE:
\begin{align}
    \ddot{X}_{t}+2\sqrt{\mu} \dot{X}_{t}+\nabla f(X_{t})=0. \tag{Polyak ODE}
\end{align}
More recently, \citep{WibisonoE7351,wilson2021lyapunov} extended these results beyond the Euclidean setting through the \textit{variational perspective}. 
Instead of concentrating on the behavior of the algorithms in the limit of their step-sizes, the variational perspective discretizes a functional on continuous-time curves, called \textit{Bregman Lagrangian}, to achieve acceleration.

Even though HB and NAG can be viewed as discretizations of the same Polyak ODE, their practical behavior can differ significantly. 
This was noted by \citet{lessard2016analysis}, who showed that for standard choices of $\beta_{k}$, HB fails to guarantee convergence. 
Motivated by this observation,  \citet{Shi2021UnderstandingTA} introduced the High-Resolution ODE (HR-ODE) framework. 
Using a backward error analysis approach, they derived differential equations that closely approximate the dynamics of each method. 
More precisely, they derived the following ODEs:
\begin{align}
    &\Ddot{X}_t+2\sqrt{\mu}\dot{X}_t+(1+\sqrt{\mu s})\nabla f(X_t)=0, 
        \tag{HB-ODE} \label{HB_ODE} \\
    \Ddot{X}_t+2\sqrt{\mu}&\dot{X}_t+\sqrt{s}\nabla^2 f(X_t)\dot{X}_t+(1+\sqrt{\mu s})\nabla f(X_t)=0.
        \tag{NAG-ODE} \label{NAG_ODE}
\end{align}
The difference between these two ODEs is attributed to the \textit{gradient correction} term present in \eqref{NAG_ODE}. 
In a similar fashion, for the convex case, the following HR-ODE for NAG was derived:
\begin{align} \label{cont_rate_match_2new}
    \Ddot{X}_t+\left(\frac{3}{t}+\sqrt{s}\nabla^2 f(X_t)\right)\dot{X}_t+\left(1+\frac{3\sqrt{s}}{2t}\right)\nabla f(X_t)=0.
\end{align}
In follow-up work, \citet{shi2019acceleration} discretized the HR-ODEs associated with the NAG algorithm using various numerical integrators and demonstrated that the Semi-Implicit Euler (SIE) discretization leads to an accelerated convergence rate. 

Within the framework of HR-ODEs, two interesting questions emerge. 
First, can these equations be interpreted through the lens of a variational perspective, similar to \citep{WibisonoE7351,wilson2021lyapunov}? 
We partially addressed this for the convex case in \citep{maskan2024variational}. 
Second, what are their appropriate discretizations that guarantee acceleration, and how do their convergence rates compare to NAG?
For example, the SIE discretization of the  NAG-ODE achieves acceleration, but it does not match the convergence rate of NAG \citep{nesterov2003introductory}. 
Furthermore, in a similar vein, an intriguing question is whether NAG itself can be recovered through a suitable discretization of the NAG-ODE.

Our main contribution in this paper is a novel extension of the variational perspective for HR-ODEs (see \Cref{sec:exforce_ODE}). 
In particular, a direct combination of these two frameworks is challenging, as it remains unclear how the Lagrangian should be modified to recover HR-ODEs. 
To address this problem, we propose an alternative approach that preserves the Lagrangian but extends the variational perspective. 
More specifically, instead of relying on the conventional Euler-Lagrange equation, we leverage the forced Euler-Lagrange equation that incorporates external forces acting on the system. 
By representing the damped time derivative of the potential function gradients as an external force, our proposed variational perspective allows us to reconstruct various HR-ODEs through specific damping parameters. 

Another significant contribution is a novel general HR-ODE for smooth and strongly convex functions (see \Cref{section4}):
\begin{align} \label{proposed_ODE_conti}
\left\{
    ~~
    \begin{aligned}
        \dot{X}_t & = -m\nabla f(X_t) -n(X_t-V_t)  \\
        \dot V_t & = -p\nabla f(X_t) -q(V_t-X_t)
    \end{aligned}
\right.
\tag{GM2-ODE}
\end{align}
with $\dot{X}_t:=\dot X(t)$ and $\dot V_t:=\dot V(t)$ denoting the first time-derivative of $X_t:=X(t)$ and $V_t:=V(t)$, respectively, and $m,n,p,q$ are non-negative parameters. 
Using Lyapunov functions, we establish convergence guarantees for various settings of $m,n,p,q$.
By analyzing \eqref{GM2-ODE}, which unifies \eqref{HB_ODE} and \eqref{NAG_ODE}, we derive convergence rates for both methods, improving upon the rates previously reported in the literature \citep{Shi2021UnderstandingTA,zhang2021revisiting,shi2019acceleration}.

In addition, we apply the Semi-Implicit Euler (SIE) on \eqref{GM2-ODE}, generating the following family of methods:
\begin{align} \label{SIE_newalg} \tag{GM2}
\left\{
    ~~
    \begin{aligned}
        x_{k+1}-x_k & = -m\sqrt{s}\nabla f(x_k)-n\sqrt{s}(x_{k+1}-v_k), \\
        v_{k+1}-v_k & =-p\sqrt{s}\nabla f(x_{{k+1}})-q\sqrt{s}(v_k-x_{k+1}).
    \end{aligned}
\right.
\end{align}
For specific parameter choices, \eqref{SIE_newalg} recovers several well-known optimization algorithms, including HB, NAG, the Triple Momentum (TM) method \citep{7967721}, and the Quasi-Hyperbolic Momentum (QHM) method \citep{ma2018quasihyperbolic}.
We establish convergence rates using Lyapunov function analysises
recovering the known rates for HB and NAG, and improve upon the existing results for QHM. We further show that this connection extends to the TM method.

Finally, we derive an improved rate for the decay of the gradient norm for the NAG algorithm \eqref{eqn_Nest_alg} for $f \in \mathcal{F}_{L}$. 
This is achieved by reformulating NAG in a structure similar to \eqref{SIE_newalg}. 
Additionally, we propose an HR-ODE inspired by the rate-matching technique used in \citep{WibisonoE7351}, demonstrating that NAG can be viewed as an approximation of this technique when applied to a specific ODE  (see \Cref{section3}). 

\section{External Forces and High-Resolution ODEs}\label{sec:exforce_ODE}

Consider a Lagrangian function $\mathcal{L}(X_t,\dot{X}_t,t)$. 
Using variational calculus, we define the action for the curves $\{X_t:t\in \mathbb R\}$ as the functional 
\begin{align}
    \mathcal{A}(X_t)=\int_{\mathbb R}\mathcal{L}(X_t,\dot{X}_t,t)dt.
\end{align}
In the absence of external forces, a curve is a stationary point for the problem of minimizing the action $\mathcal{A}(X_t)$ \textit{if and only if} it satisfies the Euler Lagrange equation: 
\begin{align}
\label{eqn:standard-EL}
    {\frac{d}{dt} \left\{  \frac{\partial \mathcal{L}}{\partial \dot{X}_t}(X_t,\dot{X}_t,t)  \right\} =\frac{\partial \mathcal{L}}{\partial X_t}(X_t,\dot{X}_t,t)}. 
\end{align}
This approach is used to derive the low resolution ODEs for convex and strongly convex functions in \citep{WibisonoE7351,wilson2021lyapunov}. 
However, the relationship between the Euler-Lagrange equation and the HR-ODE framework remains unclear. 
To the best of our knowledge, no prior work has established a connection between the Euler-Lagrange equation and this framework. 

To address this challenge, we explore a more comprehensive model that includes an external non-conservative force $F$ in the system.
When such a dissipative force is present, the \textit{forced} Euler-Lagrange equation should be used: 
\begin{align} \label{fEL}
    \frac{d}{dt}\left\{  \frac{\partial \mathcal{L}}{\partial \dot{X}_t}(X_t,\dot{X}_t,t)  \right\}-\frac{\partial \mathcal{L}}{\partial X_t}(X_t,\dot{X}_t,t)=F(X_t,\dot{X}_t,t),
\end{align}
which itself is the result of integration by parts of Lagrange d'Alembert principle \citep{campos2021discrete}. 

In the following sections, we establish connections between the forced Euler-Lagrange equation and the HR-ODE framework, addressing convex and strongly convex potential functions separately. 

\subsection{HR-ODEs for Strongly Convex Functions}

For strongly convex functions, we utilize the following Lagrangian:
\begin{align} \label{strongly_cvx_lagrange}
    \mathcal{L}(X_t,\dot{X}_t,t)=e^{\alpha_t+\beta_t+\gamma_t}\left( \frac{Ce^{-\alpha_t}}{2}\|\dot{X}_t\|^2-f(X_t)\right),
\end{align}
where $\alpha_t,\beta_t,\gamma_t:\mathbb{T}\rightarrow \mathbb{R}$ are continuously differentiable functions of time, corresponding to the weights of velocity, the potential function $f$, and the overall damping, respectively. 
For this Lagrangian, the corresponding forced Euler-Lagrange equation in (\ref{fEL}) takes the following form:
\begin{align} \label{HR_FEL}
    \Ddot{X}_t + (\dot \gamma_t+ \dot{\beta}_t)\dot{X}_t+\frac{e^{\alpha_t}}{C}\nabla f(X_t)=\frac{F}{Ce^{ \gamma_t + \beta_t}}.
\end{align}
Inspired by \citet{wilson2021lyapunov}, who show that a constant $\alpha_t$ is sufficient to derive NAG in continuous-time, we consider a fixed $\alpha_t = \alpha$ in this setting. 

For strongly convex functions, we propose the following external force:
\begin{align}
    F = -C\sqrt{s} e^{\gamma_t}\frac{d}{dt}\left(e^{\beta_t}\nabla f(X_t)\right).
\end{align}
Substituting this force into \eqref{HR_FEL} gives
\begin{align} \label{sc_eqn1}
    \Ddot{X}_t+(\dot \gamma_t+ \dot{\beta}_t)\dot{X}_t +\sqrt{s}\nabla^2f(X_t)\dot{X}_t+\bigg(\frac{e^{\alpha}}{C}+\sqrt{s}\dot{\beta}_t\bigg)\nabla f(X_t)=0.
\end{align}
The next theorem establishes the convergence guarantees for the trajectory $X_t$ to the unique minimizer $x^*$.
\begin{theorem} \label{Theorem3_1}
    Let $f \in \mathcal{F}_{\mu,L}$. 
    We consider two different parameter settings:\\
    (a) If the following scaling conditions hold:
    \begin{align} \label{mod-idela-scaling-cond-strcvx-a}
        ~\dot \gamma_t = \dot \gamma, ~~~\text{such that}~~~ 0 \leq \dot{\beta}_t \leq \dot \gamma = e^{\alpha} \leq \tfrac{\mu}{C},
    \end{align}
    then the trajectory $X_t$ from \eqref{sc_eqn1} satisfies
    \begin{align} \label{Theorem32_eqn1}
        f(X_t)-f(x^*)\leq \mathcal{O}(e^{-\beta_t}). 
    \end{align}
    (b) If the following modified scaling conditions hold:
    \begin{align} \label{mod-idela-scaling-cond-strcvx-b}
        ~
        \dot \gamma_t = \dot \gamma,
        ~~~\text{such that}~~~
        \dot \gamma = e^\alpha \leq 2 \dot{\beta} \leq \tfrac{2\mu}{C},~~~\text{and}~~~ C\sqrt{s}\leq 2 
    \end{align}
    then the trajectory $X_t$ from \eqref{sc_eqn1} satisfies
    \begin{align} \label{Theorem32_eqn2}
        f(X_t)-f(x^*)\leq \mathcal{O}(e^{-\gamma_t}). 
    \end{align}
\end{theorem}

\begin{proof}[Proof of \Cref{Theorem3_1}]
    At the center of our analysis is the construction of a suitable Lyapunov function for the corresponding ODE, serving as a non-increasing energy measure along the trajectory to ensure convergence. 
    See \citep{siegel2019accelerated,shi2019acceleration,attouch2020first,attouch2021convergence} for similar uses of Lyapunov functions in optimization. 
    Here, we present the proposed Lyapunov functions and defer the technical details of proving their decrease to \Cref{prf:thm3_1} in the supplementary material. 
    
    For the first setting, in (a), we utilize the following Lyapunov function: 
    \begin{align} \label{eqn:hrode-strcnvx-lyap1}
        \varepsilon(t) = e^{\beta_t}\left(\frac{Ce^{\alpha}}{2}\|X_t-x^*+e^{-\alpha}\dot{X}_t+\sqrt{s}e^{-\alpha}\nabla f(X_t)\|^2+f(X_t)-f(x^*)\right).
    \end{align}
    This Lyapunov function stems from an Integral Quadratic Constrained (IQC) based analysis presented in detail in \Cref{app_subsec_QCI}. 
    Under the conditions listed in \eqref{mod-idela-scaling-cond-strcvx-a}, we can show that $d\varepsilon/dt \leq 0$ along the trajectory $X_t$ from \eqref{sc_eqn1}. 
    As a result, we have $e^{\beta_t}(f(X_t)-f(x^*))\leq \varepsilon(t)\leq \varepsilon(t_0)$, which yields the desired result.
    
    For setting (b), we propose the following Lyapunov function:
    \begin{align} \label{eqn:hrode-strcnvx-lyap2}
        \begin{aligned}
            \varepsilon(t) = e^{\gamma_t}\Bigg(\frac{Ce^{\alpha}}{2}\|X_t-x^* & +e^{-\alpha}\dot{X}_t+\sqrt{s}e^{-\alpha}\nabla f(X_t)\|^2-C\Big(\frac{e^{\alpha}-\dot\beta}{2}\Big)\|X_t-x^*\|^2+f(X_t)-f(x^*)\Bigg).
        \end{aligned}
    \end{align}
    We discovered this Lyapunov function through a trial-and-error refinement process, leveraging the strong convexity of $f$ to accommodate a negative quadratic term while ensuring the function remains positive, and we can show that $d\varepsilon(t)/dt \leq 0$ along $X_t$ under the conditions listed in \eqref{mod-idela-scaling-cond-strcvx-b}.  
    
    Positivity is ensured by strong convexity, as we have
    \begin{align*}
        f(X_t)-f(x^*) \geq \tfrac{\mu}{2}\|X_t-x^*\|^2\geq C(\frac{e^{\alpha}-\dot\beta}{2})\|X_t-x^*\|^2,
    \end{align*}
    and, combined with \eqref{eqn:hrode-strcnvx-lyap2}, this leads to
    \begin{align} \label{eqn:theorem3.1_b_conv} 
        e^{\gamma_t}\Bigg(\frac{Ce^{\alpha}}{2}\|X_t-x^*+e^{-\alpha}\dot{X}_t+\sqrt{s}e^{-\alpha}\nabla f(X_t)\|^2\Bigg)\leq \varepsilon(t)\leq \varepsilon(t_0),
    \end{align}
establishing the result in \eqref{Theorem32_eqn2}. 

Next, we define 
    \begin{align*}
        V_t &:= X_t+e^{-\alpha}\dot{X}_t+\sqrt{s}e^{-\alpha}\nabla f(X_t),\qquad \dot \nonumber\\
        V_t &= \dot{X}_t + e^{-\alpha}\ddot X_t + \sqrt{s}e^{-\alpha}\nabla^2 f(X_t)\dot{X}_t. 
    \end{align*}
    Using this definition in \eqref{sc_eqn1}, we get
    \begin{align*}
       e^{\alpha}\dot V_t +  \dot\beta \dot{X}_t + (\frac{e^{\alpha}}{C}+\sqrt{s} \dot\beta)\nabla f(X_t)= 0.
    \end{align*}
    In addition, if we have $X_t-x^*+e^{-\alpha}\dot{X}_t+\sqrt{s}e^{-\alpha}\nabla f(X_t) = 0$, then $\dot V_t=0$ and
    \begin{align*}
        X_t-x^*-\frac{1}{C\dot{\beta}}\nabla f(X_t)= 0.
    \end{align*}
    Through the definition of strong convexity, we get
    \begin{align*}
        f(x^*)-f(X_t)\geq  -C\dot{\beta}\|X_t-x^*\|^2 +  \frac{\mu}{2}\|X_t-x^*\|^2\geq -  \frac{\mu}{2}\|X_t-x^*\|^2 ,
    \end{align*}
    where the last inequality holds due to $C\dot{\beta} \leq \mu$. Also, from the definition of strong convexity, we know $f(X_t)-f(x^*)\geq \frac{\mu}{2}\|X_t-x^*\|^2$.
    This implies $f(X_t) - f(x^*) = \frac{\mu}{2}\|X_t-x^*\|^2$, indicating that $X_t = x^*$.
    Thus, if $V_t = x^*$, it follows that $X_t =x^*$.
\end{proof}

\begin{remark} \label{rem:nag-ode}
    By setting $\alpha=\log(\sqrt{\mu})$, $C=\sqrt{\mu}$ and $\beta_t=\gamma_t=\sqrt{\mu}t$ in \eqref{sc_eqn1}, we recover the HR-ODE that corresponds to the NAG algorithm for strongly convex functions \citep{Shi2021UnderstandingTA}, see also \eqref{NAG_ODE}. 
    Moreover, the corresponding Lagrangian in \eqref{strongly_cvx_lagrange} matches the one proposed in \citep{wilson2021lyapunov}. 
    Finally, when $e^{\alpha}=\dot{\beta}$ the Lyapunov functions in \eqref{eqn:hrode-strcnvx-lyap1} and \eqref{eqn:hrode-strcnvx-lyap2} coincide. 
\end{remark}

In \Cref{section4}, we will demonstrate that our HR-ODE in \eqref{sc_eqn1} captures many existing accelerated methods for strongly convex functions as special cases. 
Prior to that, in the next subsection, we present our HR-ODE for general convex objectives. 

\subsection{HR-ODEs for Convex Functions}

For convex functions, we consider the following Lagrangian:
\begin{align} \label{Lagrangian2}
    \mathcal{L}(X_t,\dot{X}_t,t) =e^{\alpha_t+\gamma_t}\left(\frac{1}{2}\|e^{-\alpha_t}\dot{X}_t\|^2-e^{\beta_t}f(X_t)\right).
\end{align}
From the definition, it is easy to compute 
\begin{align} \label{Lagrange_par_der}
    \frac{\partial \mathcal{L}}{\partial \dot{X}_t}(X_t,\dot{X}_t,t)  = e^{\gamma_t-\alpha_t}\dot{X}_t,
    \qquad 
    \frac{\partial \mathcal{L}}{\partial X_t}(X_t,\dot{X}_t,t)= -e^{\gamma_t+\alpha_t+\beta_t}\nabla f(X_t).
\end{align}
By substituting these into (\ref{fEL}), we obtain the forced Euler-Lagrange equation for convex functions as follows:
\begin{align} \label{HR_general_ODE_F}
    \Ddot{X}_t + (\dot{\gamma}_t-\dot{\alpha}_t)\dot{X}_t+e^{2\alpha_t+\beta_t}\nabla f(X_t) =e^{\alpha_t-\gamma_t} F.
\end{align}
In the sequel, we propose and analyze two different choices for the external force $F$. 
The first choice recovers the structure of the HR-ODE studied in \citep{pmlr-v108-laborde20a}, while the second corresponds to the HR-ODE from \citep{shi2019acceleration}. 

\subsubsection{External Force Inspired by Laborde’s Approach}
Let us first consider the following external force:
\begin{align}
    F = -\sqrt{s}e^{\gamma_t}\frac{d}{dt}[e^{-\alpha_t}\nabla f(X_t)],
\end{align}
where $s \geq 0$ is a constant parameter. 
In this case, the forced Euler-Lagrange equation in (\ref{HR_general_ODE_F}) becomes
\begin{align} \label{HR_general_ODE_F_laborde}
    \Ddot{X}_t + (\dot{\gamma}_t-\dot{\alpha}_t)\dot{X}_t+e^{2\alpha_t+\beta_t}\nabla f(X_t) = -\sqrt{s}e^{\alpha_t}\frac{d}{dt}\left[e^{-\alpha_t}\nabla f(X_t)\right].
\end{align}

It is possible to demonstrate the convergence of $X_t$ to $x^*$ and establish the convergence rate, as stated in the following theorem. 

\begin{theorem}\label{Theorem_ODE_laborde}
    Let $f \in \mathcal{F}_L$. 
    If the ideal scaling conditions $\dot\beta_t\leq \dot \gamma_t = e^{\alpha_t}$ hold, then the trajectory $X_t$ from (\ref{HR_general_ODE_F_laborde}) satisfies 
    $$f(X_t)-f(x^*)\leq \mathcal{O}(e^{-\beta_t}).$$
\end{theorem}

\begin{proof}[Proof of \Cref{Theorem_ODE_laborde}]
We define our Lyapunov function as follows: 
\begin{align}
    \varepsilon(t)=\frac{1}{2}\|X_t+e^{-\alpha_t}\dot{X}_t-x^*+\sqrt{s}e^{-\alpha_t}\nabla f(X_t)\|^2+e^{\beta_t}(f(X_t)-f(x^*)),
\end{align}
and show that it satisfies $d\varepsilon/dt \leq 0$ along the trajectory $X_t$ from \eqref{HR_general_ODE_F_laborde}. 
Then, denoting the initial point by $t_0$, we have 
$e^{\beta_t}(f(X_t)-f(x^*))\leq \varepsilon(t)\leq \varepsilon(t_0),$
leading to the desired result. 
Technical details are deferred to \Cref{thm1_proof}. 
\end{proof}

\begin{remark}
We can derive existing HR-ODEs from \eqref{HR_general_ODE_F_laborde}. 
For example, by following the convention used in \citep{WibisonoE7351}, we reparameterize the equation in terms of the functions $n(t)$ and $q(t)$, by choosing 
\begin{align} \label{prams_general}
        \alpha_t=\log (n(t)),\quad
    \beta_t=\log (q(t)/n(t)),\quad 
    \dot\gamma_t=e^{\alpha_t}=n(t).
\end{align}
In this setting, using $n(t)=p/t$ and $q(t)=Cpt^{p-1}$, equation (\ref{HR_general_ODE_F_laborde}) becomes
\begin{align} \label{HR_cnts_general_npq2}
    \Ddot{X}_t + \left(\frac{p+1}{t}+\sqrt{s}\nabla^2 f(X_t)\right)\dot{X}_t+\left(Cp^2t^{p-2}+\frac{\sqrt{s}}{t}\right)\nabla f(X_t)=0.
\end{align}
With $p=2$ and $C=1/4$, this equation captures the HR-ODE from \citep{pmlr-v108-laborde20a} as a special case.\end{remark}

\subsubsection{External Force Inspired by Shi’s Approach}
Consider the external force given by
\begin{align}
    F = -\sqrt{s}e^{\gamma_t-\beta_t}\frac{d}{dt}\left[e^{-(\alpha_t-\beta_t)}\nabla f(X_t)\right]. 
\end{align}
In this case, replacing $F$ in the forced Euler-Lagrange equation in (\ref{HR_general_ODE_F}) gives
\begin{align} \label{HR_general_ODE_F_Shi}
    \Ddot{X}_t + (\dot{\gamma}_t-\dot{\alpha}_t)\dot{X}_t+e^{2\alpha_t+\beta_t}\nabla f (X_t)= -\sqrt{s}e^{\alpha_t-\beta_t}\frac{d}{dt}[e^{-(\alpha_t-\beta_t)}\nabla f(X_t)].
\end{align}

\begin{theorem}\label{Theorem_ODE_Shi}
    Let $f \in \mathcal{F}_L$. 
    If the modified ideal scaling conditions 
    \begin{align*}
        \dot{\beta}_t \leq \dot{\gamma}_t = e^{\alpha_t} ~~~\text{and}~~~ \Ddot{\beta}_t \leq e^{\alpha_t} \dot{\beta}_t + 2 \dot{\alpha}_t \dot{\beta}_t
    \end{align*}
    hold, then the trajectory $X_t$ from (\ref{HR_general_ODE_F_Shi}) satisfies
    \begin{align*}
        f(X_t)-f(x^*)\leq \mathcal{O}\left((e^{\beta_t}+\sqrt{s}e^{-2\alpha_t}\dot{\beta}_t)^{-1}\right).
    \end{align*}
\end{theorem}

\begin{proof}[Proof of \Cref{Theorem_ODE_Shi}]
Our main tool in the proof is the following Lyapunov function:
\begin{align*}
    \varepsilon(t)=\frac{1}{2}\|X_t+e^{-\alpha_t}\dot{X}_t-x^*+\sqrt{s}e^{-\alpha_t}\nabla f(X_t)\|^2+(e^{\beta_t}+\sqrt{s}e^{-2\alpha_t}\dot{\beta}_t)(f(X_t)-f(x^*)).
\end{align*}
Under the conditions listed, we can show that $d\varepsilon/dt \leq 0$ along $X_t$ from \eqref{HR_general_ODE_F_Shi}. 
As a result, we have $(e^{\beta_t}+\sqrt{s}e^{-2\alpha_t}\dot{\beta}_t)(f(X_t)-f(x^*))\leq \varepsilon(t)\leq \varepsilon(t_0),$ where $t_0$ is the initial time. 
This yields the desired result. 
We defer the technical details to \Cref{thm3_proof}. 
\end{proof}

\begin{remark}
    \Cref{Theorem_ODE_Shi} establishes a faster convergence rate than \Cref{Theorem_ODE_laborde}.
\end{remark}

\begin{remark}
    We can connect the forced Euler-Lagrange equation in \eqref{HR_general_ODE_F_Shi} to the known HR-ODEs. 
    Using the same parametrization in (\ref{prams_general}), we obtain  
    \begin{align*}
        \Ddot{X}_t + \left(\frac{p+1}{t}+\sqrt{s}\nabla^2 f(X_t)\right)\dot{X}_t+\left(Cp^2t^{p-2}+\frac{\sqrt{s}(p+1)}{t}\right)\nabla f(X_t)=0.
    \end{align*}
    With $p=2$ and $C=1/4$, this equation gives the following HR-ODE:
    \begin{align} \label{HR_cnts_general_stable_ODE}
        \Ddot{X}_t + \left(\frac{3}{t}+\sqrt{s}\nabla^2 f(X_t)\right)\dot{X}_t+\left(1+\frac{3\sqrt{s}}{t}\right)\nabla f(X_t)=0.
    \end{align}
    This HR-ODE was studied in \citep{shi2019acceleration}, where it was discretized using the Semi-Implicit Euler (SIE) and Implicit Euler (IE) schemes, resulting in optimization algorithms with accelerated convergence rates. 
\end{remark}

\section{Implications of the HR-ODEs for Strongly Convex Functions}
\label{section4}

In this section, we examine the HR-ODE in \eqref{sc_eqn1}. 
For simplicity and to maintain notational consistency with prior work, we reparametrize \eqref{sc_eqn1} using the non-negative constants $m,n,p$ and~$q$:
\begin{align} \label{SC_ODE_represented}
    \Ddot{X}_t+(n+q)\dot{X}_t + m\nabla^2f(X_t)\dot{X}_t + (np+mq)\nabla f(X_t)=0,
\end{align}
which follows by choosing $\alpha = \log n,$ $ \beta = qt,$ $\gamma = nt,$ $m=\sqrt{s},$ and $C=1/p$. 
Equation \eqref{SC_ODE_represented} can be restated in the compact form as follows: 
\begin{align} \label{GM2-ODE} \tag{GM2-ODE}
    \left\{
    ~~
    \begin{aligned}
        \dot{X}_t & = -m\nabla f(X_t) -n(X_t-V_t)  \\
        \dot V_t & = -p\nabla f(X_t) -q(V_t-X_t)
    \end{aligned}
    \right.
\end{align}
In what follows, we explore the implications of \eqref{proposed_ODE_conti} in both continuous and discrete-time settings. 

\subsection{Continuous-Time Analysis}
\label{subsec_continuoustime}

The following theorem presents a reformulation of the convergence result in \Cref{Theorem3_1} for \eqref{GM2-ODE}, offering an alternative proof that provided key insights into the choice of Lyapunov functions (see \Cref{rmk:insight}).

\begin{theorem} \label{theorem2}
    Let $f$ be a $\mu$-strongly convex function. 
    Assume that $ m, q \geq 0$ and $n,p > 0$. 
    (a) Suppose that $n = q \leq p \mu$. 
    Then, trajectories $X_t$ and $V_t$ are globally asymptotically stable for \eqref{proposed_ODE_conti}, satisfying 
    \begin{align} \label{Lyap1}
        \varepsilon(t) \leq e^{-qt}\varepsilon(0) \quad \text{with} \quad \varepsilon (t)= f(X_t)-f(x^*)+\tfrac{q}{2p}\|V_t-x^*\|^2.
    \end{align}
    (b) Suppose that $n\leq 2q \leq 2 p \mu$ and $m\leq 2p$. 
    Then, trajectories $X_t$ and $V_t$ are globally asymptotically stable for \eqref{proposed_ODE_conti}, satisfying
    \begin{align} \label{eqn:new-lyapunov}
        \varepsilon(t) \leq e^{-nt}\varepsilon(0) ~~ \text{with} ~~ \varepsilon(t)=f(X_t)-f(x^*)-\tfrac{n-q}{2p}\|X_t-x^*\|^2+\tfrac{n}{2p}\|V_t-x^*\|^2.
    \end{align}
    Consequently, both $X_t$ and $V_t$ converge to $x^*$ as $t \to \infty$. 
\end{theorem}

\begin{proof}[Proof of \Cref{theorem2}]
    (a) This result can be obtained by choosing $\alpha = \log n,$ $ \beta = qt,$ $\gamma = nt,$ $m=\sqrt{s},$ and $C=1/p$ in \Cref{Theorem3_1} (a). 
    (b) Similarly, the result follows directly by applying the same parameter settings $\alpha = \log n,$ $ \beta = qt,$ $\gamma = nt,$ $m=\sqrt{s},$ and $C=1/p$ in \Cref{Theorem3_1} (b). 
\end{proof}

\begin{remark} \label{rmk:insight}
    Although \Cref{theorem2} follows directly from \Cref{Theorem3_1} by selecting appropriate parameters, part~(a) of \Cref{theorem2} was, in fact, derived first in our development and served as a key insight for \Cref{Theorem3_1}. 
    For clarity and instructional value, we outline this initial proof:     
    Our construction of the Lyapunov function in this setting is inspired by the connections between integral quadratic constraints (IQC) from stability analysis in control theory and Polyak's ODE \citep{doi:10.1137/20M1364138, doi:10.1137/17M1136845}. 
    Adopting this perspective, we reformulate \eqref{proposed_ODE_conti} in the form of a non-linear feedback control system, allowing us to utilize the analytical framework established in \citep{doi:10.1137/17M1136845}. 
    Specifically, we apply \citep[Theorem~6.4]{doi:10.1137/17M1136845}, with the constraints on the parameters $m,n,p$ and $q$ in our theorem designed to satisfy the requirements of the referenced analysis. 
    See \Cref{prf:conti_conv} in the supplementary material for the complete proof.
    Additionally, while the Lyapunov function in \eqref{eqn:new-lyapunov} covers the one in \eqref{Lyap1} as a special case, we present both for clarity. 
    The function in \eqref{Lyap1} provides a clear dynamical systems interpretation, as outlined above, whereas \eqref{eqn:new-lyapunov} was obtained by extending \eqref{Lyap1} under the given parameter conditions and the strong convexity of $f$, which ensure the non-negativity of the first two terms. 
\end{remark}

In what follows, we discuss the implications of our continuous-time analysis and compare our results with related works.

\begin{remark} \label{rem1}
    Let us rewrite (\ref{Lyap1}) by using the definition of $V_t$ from \eqref{proposed_ODE_conti}:
    \begin{align} \label{Lyap2}
        \varepsilon (t)= f(X_t)-f(x^*)+\tfrac{q}{2p}\|\tfrac{\dot{X}_t}{n}+X_t-x^* + \tfrac{m}{n}\nabla f(X_t)\|^2. 
    \end{align}
    If we eliminate the gradient term by setting $ m=0 $, we recover the ODE presented in \citep[Theorem~4.3]{doi:10.1137/20M1364138}. 
    Therefore, our Lyapunov function in \eqref{Lyap1} represents a non-trivial extension of its counterpart in \citep[Theorem~4.3]{doi:10.1137/20M1364138} incorporating a high-resolution correction term. 
\end{remark}

\begin{remark}
    For the specific case of \eqref{NAG_ODE}, \Cref{theorem2} proves a rate of $\mathcal{O}\left(e^{-\sqrt{\mu}t}\right)$. 
    This rate is faster than previous convergence rates for this ODE. 
    Specifically, \citet{Shi2021UnderstandingTA} and \citet{zhang2021revisiting} showed convergence rates of $\mathcal{O}\left(e^{-\frac{\sqrt\mu}{4}t}\right)$ and $\mathcal{O}\left(e^{-\frac{\sqrt\mu}{2}t}\right)$, respectively. 
    With different parameter choices, various ODEs can be constructed that achieve the same rate of convergence as (\ref{NAG_ODE}). 
    For example, for fixed values of $n$ and $q$, increasing $p$ would result in different ODEs, but all with the same convergence rate due to the fixed $n$. 
\end{remark}

\begin{remark}
    The continuous-time behavior of many methods emerges as a special case of \eqref{GM2-ODE} with an appropriate selection of the parameters $m$, $n$, $p$, and $q$. 
    Notably, \eqref{GM2-ODE} recovers both \eqref{NAG_ODE} and \eqref{HB_ODE} for specific parameter choices, as reported in \Cref{table_compare1}. 
    Similarly, the continuous-time dynamics of other methods, such as gradient descent and H-NAG, can also be recovered. 
    Moreover, \Cref{theorem2} provides a means to determine the continuous-time convergence rates of many of these methods (see \Cref{table_compare1}).
\end{remark}

\begin{table}
  \caption{Various ODEs and their convergence rates from the literature are recovered from \eqref{proposed_ODE_conti} and \Cref{theorem2}.}
  \label{table_compare1}
  \centering
  \resizebox{\textwidth}{!}{
  \begin{tabular}{ccc}
    \toprule
    \multicolumn{1}{c}{ODE}&\multicolumn{1}{c}{Parameters \eqref{proposed_ODE_conti}}   &     \multicolumn{1}{c}{Convergence Rate (\Cref{theorem2})}               \\
    \midrule
    
    \multirow{2}{*}{Gradient Flow}  & 
 $n=0,m =1$, & \multirow{2}{*}{$\mathcal{O}(e^{-2\mu t})$}  \\
& any $p,q$&     \\

\multirow{2}{*}{Polyak's ODE}  & 
 $q=n=\sqrt{\mu}$& \multirow{2}{*}{$\mathcal{O}(e^{- \sqrt{\mu}t})$}\\
  &  $m=0,p=\frac{1}{\sqrt{\mu}}$ & \\
   \multirow{2}{*}{\eqref{HB_ODE} } & 
 $n=q=\sqrt{\mu}$, & \multirow{2}{*}{$\mathcal{O}(e^{-\sqrt
 {\mu}t})$ } \\
&$m=0,p=\tfrac{1}{\sqrt{\mu}}+\sqrt{s}$&    \\

     \multirow{2}{*}{\eqref{NAG_ODE}}  & 
 $n=q=\sqrt{\mu}$, & \multirow{2}{*}{$\mathcal{O}(e^{-\sqrt
 {\mu}t})$} \\
&$m=\sqrt{s},p=\tfrac{1}{\sqrt{\mu}}$&     \\

Continuous H-NAG\citep{chen2019first} & $n=1, q=\tfrac{\mu}{\gamma}$ & \multirow{2}{*}{$\mathcal{O}(e^{- t})$} \\
($\gamma=\mu$)& $p=\frac{1}{\gamma},m=\beta$ & \\

    \bottomrule
  \end{tabular}
  }
\end{table}

\begin{remark}
    \eqref{GM2-ODE} is not the first attempt toward unification of accelerated methods. 
    An earlier general ODE, denoted as \eqref{GM-ODE}, was proposed by \citet{zhang2021revisiting}:
    \begin{align} \label{GM-ODE}
        \left\{
        ~~
        \begin{aligned}
            \dot{U}_t &= -m'\nabla f(U_t)-n'W_t, \\
            \dot{W}_t &= \nabla f(U_t)-q'W_t,
        \end{aligned}
        \right.
        \tag{GM-ODE}
    \end{align}
    
    where $m',n',q'\geq 0$, $U_t=U(t)$, and $W_t=W(t)$. 
    A comparison between \eqref{GM-ODE} and \eqref{GM2-ODE} reveals that our general model achieves a better convergence rate. 
    More details can be found in \Cref{sec:gmvsgm2}. 
\end{remark}

\begin{figure}
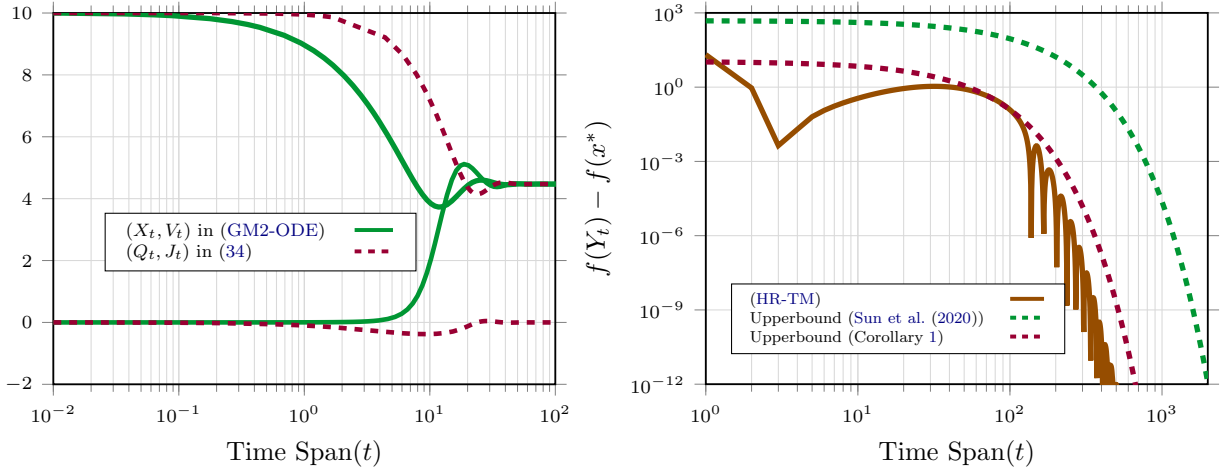

\centering
\definecolor{mycolor1}{rgb}{0.59,0.3,0.0}
\definecolor{mycolor2}{rgb}{0,0.59,0.2}
\definecolor{mycolor3}{rgb}{0.6,0,0.2}
 
\caption{Continuous-time simulation using $f(x)=4(L-\mu)\log(1+e^{-x})+\frac{\mu}{2}x^2$, (left) comparison of trajectories $(X_t,V_t)$ from \eqref{proposed_ODE_conti} and $(Q_t,J_t)$ from (\ref{proposed_ODE_conti2}). The simulation is done for ${L=1,\mu=10^{-2}}$ and $X(0)=Q(0)=10, V(0)=J(0)=0$, (right) trajectory of (\ref{HR_TM}) ODE for $L=10, \mu=10^{-3}, \text{ random }Y(0),\dot Y(0)=0$  and its corresponding upper bounds.}
        \label{fig1}
\end{figure}

\begin{remark}
    The phase space representation of \eqref{NAG_ODE} was studied in \citep{Shi2021UnderstandingTA} and is given by:
    \begin{align} \label{proposed_ODE_conti2}
        \left\{
        ~~
        \begin{aligned}
            \dot Q_t & = J_t,   \\ 
            \dot J_t & = -(2\sqrt{\mu}+\sqrt{s} \nabla^2 f(Q_t))J_t-(1+\sqrt{\mu s})\nabla f(Q_t) .
        \end{aligned}
        \right.
    \end{align}
    There are two major differences between \eqref{proposed_ODE_conti2} and the representation \eqref{NAG_ODE} through \eqref{GM2-ODE}. 
    First, when $Q_t$ in \eqref{proposed_ODE_conti2} reaches the stationary point of $f$, $J_t$ converges to zero. 
    In contrast, both $X_t$ and $V_t$ in \eqref{GM2-ODE} converge to the stationary point $x^*$ (see \Cref{fig1}). 
    Second, the absence of the Hessian term in \eqref{GM2-ODE} is advantageous, particularly when discretizing HR-ODEs to recover discrete-time algorithms. 
    A more detailed discussion is deferred to \Cref{sec:SIE_rec_NAG}.
\end{remark}

\subsubsection{New Results for the TM Method}

TM is recognized as the first-order method with the fastest convergence guarantees for minimizing smooth and strongly convex functions \citep{7967721}. 
The state space presentation of the TM method is given by:
\begin{align} \label{TM_method}
    \left\{\begin{array}{rcl}
        \epsilon_{k+1}&=&(1+\beta)\epsilon_k-\beta\epsilon_{k-1}-\alpha\nabla f(y_k),  \\
        y_k &=&(1+\gamma)\epsilon_k-\gamma \epsilon_{k-1},\\
        x_k&=&(1+\delta)\epsilon_k-\delta\epsilon_{k-1}, 
    \end{array}\right.\tag{TM-Method}
\end{align}
with initialization $\epsilon_0,\epsilon_{-1}\in\mathbb{R}^d$ and the parameters $\alpha = \frac{1+\rho}{L}$, $\beta = \frac{\rho^2}{2-\rho}$, $\gamma = \frac{\rho^2}{(1+\rho)(2-\rho)}$, and $\delta = \frac{\rho^2}{1-\rho^2}$, where $\rho=1-\sqrt{1/\kappa}$ and $\kappa=L/\mu$. 
At first glance, the connection between this method and \eqref{SIE_newalg} is not immediately apparent. 
However, if we eliminate $\epsilon_k$, we obtain  
    \begin{align}
        x_{k+1} &=\frac{\sqrt{L}-\sqrt{\mu}}{\sqrt{L}+\sqrt{\mu}}\left(x_{k}-\frac{1}{L}\nabla f(x_{k})\right)+\left(\frac{2\sqrt{\mu }}{\sqrt{L}+\sqrt{\mu}}\right)v_k,\nonumber\\
        v_{k+1}&=\sqrt{\frac{\mu}{L}}\left(x_{k+1}-\frac{1}{\mu}\nabla f(x_{k+1})\right)+\left(1-\sqrt{\frac{\mu}{L}}\right)v_{k}, \nonumber
    \end{align}
which is in the form of \eqref{SIE_newalg} with the choice of parameters:
\begin{align} \label{coeff_TM}
    m=\frac{1}{\sqrt{L}}, ~n=\frac{2\sqrt{\mu L}}{\sqrt{L}-\sqrt{\mu}},~ q=\sqrt{\mu},~  p=\frac{1}{\sqrt{\mu}},~\text{and } v_k:=y_k.
\end{align}
Inserting these parameters in \eqref{GM2-ODE} gives:
\begin{align} \label{TM_new_high_res}
    \Ddot{X}_t+\sqrt{\mu}\left(\frac{3-\sqrt{\tfrac{\mu}{L}}}{1-\sqrt{\tfrac{\mu}{L}}}\right)\dot{X}_t +\sqrt{s}\nabla^2 f(X_t)\dot{X}_t+\left(\frac{2}{1-\sqrt{\tfrac{\mu}{L}}}+\sqrt{s\mu}\right)\nabla f(X_t)=0
\end{align}
with $s=1/L$. 

To satisfy the parameter conditions in \Cref{theorem2}, we approximate $n$ in \eqref{coeff_TM} as $n\approx 2\sqrt{\mu}$, which leads to the following HR-ODE for the TM method:
\begin{align} \label{TM_new_high_res2}
    \Ddot{X}_t+3\sqrt{\mu}\dot{X}_t +\sqrt{s}\nabla^2 f(X_t)\dot{X}_t+\left(2+\sqrt{s\mu}\right)\nabla f(X_t)=0. \tag{HR-TM2}
\end{align}
The following corollary of \Cref{theorem2} presents the convergence rate of \eqref{TM_new_high_res2}. 

\begin{corollary} \label{corollary_hrtm2}
    With the choice of $m = \sqrt{s}$, $q = \sqrt{\mu}$, $n = 2\sqrt{\mu}$, and $p = 1/\sqrt{\mu}$, \eqref{proposed_ODE_conti} reduces to \eqref{TM_new_high_res2}. 
    Furthermore, invoking \Cref{theorem2} with these parameters yields the following convergence guarantee:
    \begin{align}
        f(X_t) - f(x^*) \leq C'_{TM} e^{-2 \sqrt{\mu} t},
    \end{align}
    where $C'_{TM}$ is a positive constant.
\end{corollary}

\begin{figure}
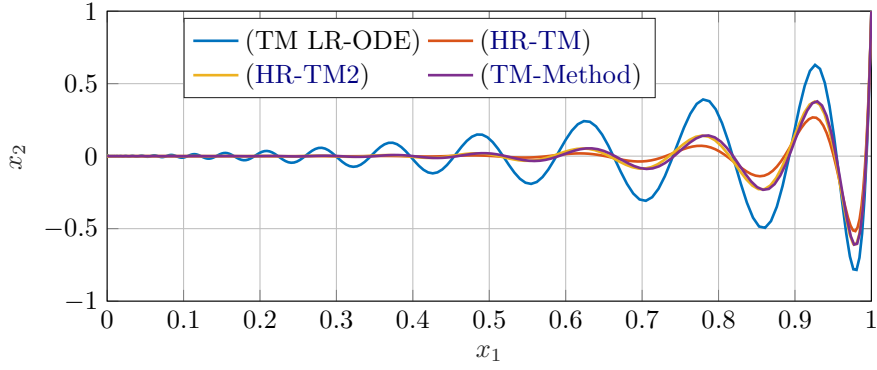

\centering
\definecolor{mycolor1}{rgb}{0.00000,0.44700,0.74100}
\definecolor{mycolor2}{rgb}{0.85000,0.32500,0.09800}
\definecolor{mycolor3}{rgb}{0.92900,0.69400,0.12500}
\definecolor{mycolor4}{rgb}{0.49400,0.18400,0.55600}
 
\caption{Trajectory of TM method and various proposed ODEs for $f(x_1,x_2)= 5\times 10^{-3}x_1^2+x_2^2$ with starting point $(x_1(0),x_2(0))=(1,1)$. In this figure, (TM LR-ODE) corresponds to the low-resolution ODE of the TM method ((\ref{TM_new_high_res2}) with $\sqrt{s}=0$). The step-size was $s=0.16$; (\ref{HR_TM}) and (\ref{TM_new_high_res2}) were according to Corollaries \ref{corl_TM} and \ref{corollary_hrtm2} respectively.}
    \label{TM_method_fig}
\end{figure}

Recently, \citet{DBLP:conf/cdc/SunGK20} proposed another HR-ODE for the TM method, using $y_k=Y(t_k)$ and $t_k=k\sqrt{\alpha}$:
\begin{align} \label{HR_TM}
    \Ddot{Y}_t&+2\sqrt{M}\dot Y_t+\gamma(1+\sqrt{M \alpha})\sqrt{\alpha}\nabla^2 f(Y_t)\dot Y_t
    +(1+\sqrt{M \alpha})\nabla f(Y_t)=0 ,\tag{HR-TM}
\end{align}
where $M=\big(\frac{1-\beta}{\sqrt{\alpha}(1+\beta)}\big)^2$ is a $\kappa$-dependent constant. 
It is easy to verify $\mu \leq M \leq 1.3661\mu$.
The convergence rate for \eqref{HR_TM} was reported in \citep[Theorem 4.1]{DBLP:conf/cdc/SunGK20} as 
\begin{align} \label{HR_TM_old_rate}
    {f(Y_t)-f(x^*)}\leq \tilde{C}_{TM} e^{-p_{TM}^* t}
\end{align}
where $\tilde{C}_{TM}$ and $p^*_{TM}$ are positive constants, with $p^*_{TM}\leq \sqrt{M}/2$. 
Our next corollary, which follows from \Cref{theorem2}, improves upon this convergence rate. 
The proof is provided in \Cref{proof_TM}.
\begin{corollary} \label{corl_TM}
    \eqref{proposed_ODE_conti} reduces to (\ref{HR_TM}) with the choice of $$q=\xi\sqrt{M},~m=\gamma\sqrt{\alpha}(1+\sqrt{M \alpha}),~n=(2-\xi)\sqrt{M}, \text{ and }~~ {p=\frac{(1-\xi\gamma \sqrt{M \alpha})(1+\sqrt{M \alpha})}{(2-\xi)\sqrt{M}}},$$
    where $\xi$ is a universal constant equal to $2/3$. 
    Furthermore, invoking \Cref{theorem2} (b) with these parameters leads to 
    \begin{align} \label{HR_TM_our_rate}
        f(Y_t)-f(x^*)\leq C_{TM}e^{-(2-\xi) \sqrt{M} t}, 
    \end{align}
    where $C_{TM}$ is a positive constant. 
\end{corollary}

\begin{remark}
    The rate in \eqref{HR_TM_our_rate} is faster than the rate in \eqref{HR_TM_old_rate}, since
    \begin{align*}
        e^{-(2-\xi) \sqrt{M} t} = e^{-\tfrac{4}{3} \sqrt{M} t} < e^{-\tfrac{1}{2} \sqrt{M} t} \leq e^{-p_{TM}^* t}.
    \end{align*}
    A numerical example illustrating this result is provided in \Cref{fig1}.
\end{remark}

\begin{remark}
    In contrast to \eqref{HR_TM}, the proposed \eqref{TM_new_high_res2} recovers the low-resolution ODE of the TM method (as introduced by \citet{kim2024convergence}) as the resolution parameter $\sqrt{s}$ approaches zero. 
    Also note that \eqref{TM_new_high_res2} exhibits a better convergence rate than \eqref{HR_TM}, compare the rates in \Cref{corollary_hrtm2} and \Cref{corl_TM}. 
    A numerical illustration of these results is presented in \Cref{TM_method_fig}, showing that \eqref{TM_new_high_res2} follows the trajectory of the TM method more closely than both \eqref{HR_TM} and the low-resolution ODE.
\end{remark}

\subsection{Discrete-Time Analysis}\label{subsec_discrete_SC}

In this section, we derive a new general accelerated first-order method from our continuous-time characterization. 
By applying SIE discretization to \eqref{proposed_ODE_conti}, we introduce the following iterative method: 
\begin{align} \tag{GM2}
\left\{
    ~~
    \begin{aligned}
        x_{k+1}-x_k & = -m\sqrt{s}\nabla f(x_k)-n\sqrt{s}(x_{k+1}-v_k), \\
        v_{k+1}-v_k & =-p\sqrt{s}\nabla f(x_{{k+1}})-q\sqrt{s}(v_k-x_{k+1}).
    \end{aligned}
\right.
\end{align}
The following theorem provides convergence guarantees for \eqref{SIE_newalg}. 

\begin{theorem}\label{theorem4}
    Let $f \in \mathcal{F}_{\mu,L}$, and consider the iterative method described in (\ref{SIE_newalg}), initialized with an arbitrary point $x_0$ and $v_0=x_0-\frac{m}{n}\nabla f(x_0)$. 
    If the following conditions hold: $q/p\leq \mu$, $0\leq nps \leq m\sqrt{s}\leq 1/L$, $0\leq q\sqrt{s} < 1$, $p>0$, and $n=q$, then the Lyapunov function 
    \begin{align} \label{lyap_discrete}
        \varepsilon(k)= f(x_k)-f(x^*)+\frac{n}{2p}\|v_k-x^*\|_2^2-\frac{nps}{2}\|\nabla f(x_{k})\|^2,
    \end{align}
    decreases as $\varepsilon(k)\leq (1-q\sqrt{s})^k\varepsilon(0)$. 
    Furthermore, we have
    \begin{align}
        f(x_{k})-f(x^*)\leq C''_{GM}(1-q\sqrt{s})^k,
    \end{align}
    where $C''_{GM}$ is a positive constant.
\end{theorem}

\begin{remark}
    The Lyapunov function in \eqref{lyap_discrete} converges to its continuous-time counterpart in \eqref{Lyap2} as the step-size $s$ approaches zero. 
    The terms in the Lyapunov function are the same as in \eqref{Lyap2}, except for the negative quadratic term. 
    Note that the smoothness assumption ensures that the Lyapunov function remains positive despite this term. 
\end{remark}

\begin{proof}[Proof of \Cref{theorem4}]
    We begin by bounding the update differences of the Lyapunov function, leading to $\varepsilon(k)\leq (1-q\sqrt{s})^k\varepsilon(0).$ 
    Although the convergence of the Lyapunov function does not directly imply the convergence of $f(x_k)-f(x^*)$ due to the negative quadratic term, it does establish the convergence rate for $\|v_k-x^*\|^2$. 
    A limit analysis on the Lyapunov function and the algorithm \eqref{SIE_newalg} then completes the proof by showing that convergence in $v_k$ implies convergence in $f(x_k)$. 
    See \Cref{prf:disc_conv} for details. 
\end{proof}

\begin{table}
  \caption{Various algorithms and their convergence rates from the literature are recovered from \eqref{SIE_newalg} and \Cref{theorem4}. Specifically, \eqref{HB_ODE} is used with $\alpha = \tfrac{1-\sqrt{\mu s}}{1+\sqrt{\mu s}}$. For H-NAG, we consider $\gamma$ as a constant.}
  \label{table_compare2}
  \centering
 
    \begin{tabular}{ccc}
    \toprule
    \multicolumn{1}{c}{Algorithm}&\multicolumn{1}{c}{Parameters \eqref{SIE_newalg}}   &     \multicolumn{1}{c}{Convergence Rate (\Cref{theorem4})}               \\
    \midrule
   \multirow{2}{*}{Heavy Ball } & 
 $n=q=\tfrac{1-\alpha}{\sqrt{s}(1+\alpha)},$, & \multirow{2}{*}{No Convergence} \\

&$m =0$,$p=\tfrac{1}{n}+\sqrt{s},$&    \\
     \multirow{2}{*}{NAG}  & 
 $n=q=\sqrt{\mu}$, & \multirow{2}{*}{$\mathcal{O}\left( (1-\sqrt{\tfrac{1}{\kappa}})^k\right)$} \\

&$m=\sqrt{s},p=\tfrac{1}{\sqrt{\mu}}$&     \\

 \multirow{2}{*}{H-NAG ($\gamma = \mu(1-\alpha)$)}  & 
 $n= 1 ,q=\tfrac{\mu}{\gamma+\mu\alpha}$, & \multirow{2}{*}{$\mathcal{O}\left( (1-\sqrt{\tfrac{1}{\kappa}})^k\right)$}  \\
& $p=\tfrac{1}{\gamma+\mu\alpha},m=\beta$&     \\

    \bottomrule
  \end{tabular}
  
\end{table}

Several well-known optimization algorithms, such as NAG and HB, can be recovered through \eqref{SIE_newalg} by selecting appropriate parameters. 
Using these parameters together with \Cref{theorem4}, we can compute the convergence rates of these methods. 
\Cref{table_compare2} summarizes some of these algorithms, their corresponding parameter selections, and their convergence rates as derived from \Cref{theorem4}. 
Notably, these parameters match those used to characterize the continuous-time behavior of these methods, a property we refer to as \emph{coefficient consistency}.

Next, we discuss the implications of our discrete-time analysis and compare our results with recent works.
\begin{remark}
    Although the algorithm \eqref{SIE_newalg} includes the HB method as a special case, the corresponding parameter setting (see \Cref{table_compare1}) results in $s=0$, which consequently nullifies the convergence rate established in \Cref{theorem4}. 
    This observation aligns with recent findings in \citep{goujaud2023provable}, which demonstrate the existence of an $L$-smooth, $\mu$-strongly convex objective function and an initialization where the HB method fails to converge, instead cycling over a finite number of iterates. 
\end{remark}

 \begin{remark}
     \eqref{GM-ODE} was also discretized using SIE in \citep{zhang2021revisiting}. 
     After recovering accelerated methods such as NAG, HB, and QHM, they analyzed the convergence of these methods. 
     However, the rates found did not match the previously known rates for NAG and HB. 
     Moreover, they used slightly different coefficients than their continuous-time analysis to recover these methods. 
     In other words, the SIE discretization of \eqref{GM-ODE} was not coefficient consistent. 
     Further details and a comparison of convergence rates are provided in \Cref{sec:gmvsgm2}.
 \end{remark}
 
\begin{remark}
    Another key insight is that the SIE discretization of \eqref{NAG_ODE} exactly recovers the NAG method.  
    Notably, previous attempts to derive accelerated methods using SIE discretization \citep{shi2019acceleration} did not recover the precise NAG algorithm. 
    See \Cref{sec:SIE_rec_NAG} for a more detailed discussion.
\end{remark}

\begin{remark}
    The Hessian-driven Nesterov Accelerated Gradient algorithm (H-NAG, see equation~(42) in \citep{chen2019first}) achieves a convergence rate of $\mathcal{O}\big((1+\sqrt{\frac{\min\{\gamma,\mu\}}{L}})^{-k}\big)$. 
    H-NAG itself is not coefficient consistent. 
    This is why H-NAG has different parameters in \Cref{table_compare1} and \Cref{table_compare2}. 
    Here, we focus on the specific case where H-NAG is coefficient consistent, with ${\gamma=\mu(1-\alpha)}$ for some $\alpha\geq \sqrt{1/\kappa}$, in which case the convergence rate becomes $\mathcal{O}\big((1+\sqrt{(1-\alpha)/\kappa})^{-k}\big)$. 
    Our method in \eqref{SIE_newalg} can recover H-NAG (see \Cref{table_compare2}). 
    In particular, for this setting we have $q=n=1$, and for a step-size $\alpha = \sqrt{1/\kappa}$, the convergence rate becomes $\mathcal{O}\big((1-\sqrt{1/\kappa})^k\big)$. 
    This is clearly an improvement since 
    \begin{align}
        1-\sqrt{1/\kappa}\leq 1-\sqrt{(1-\alpha)/\kappa} \leq \frac{1-\sqrt{(1-\alpha)/\kappa}}{1-(1-\alpha)/\kappa} = \frac{1}{1+\sqrt{(1-\alpha)/\kappa}}.
    \end{align}
    In addition, our method can also directly recover NAG, eliminating the need for the specific parameter choices required for recovery through H-NAG. 
    In this case, we further improve H-NAG's rate to $\mathcal{O}\big((1-\sqrt{1/\kappa})^k\big)$, which matches the original rates of the NAG algorithm shown in \citep{nesterov2003introductory}.
\end{remark} 

\begin{figure}
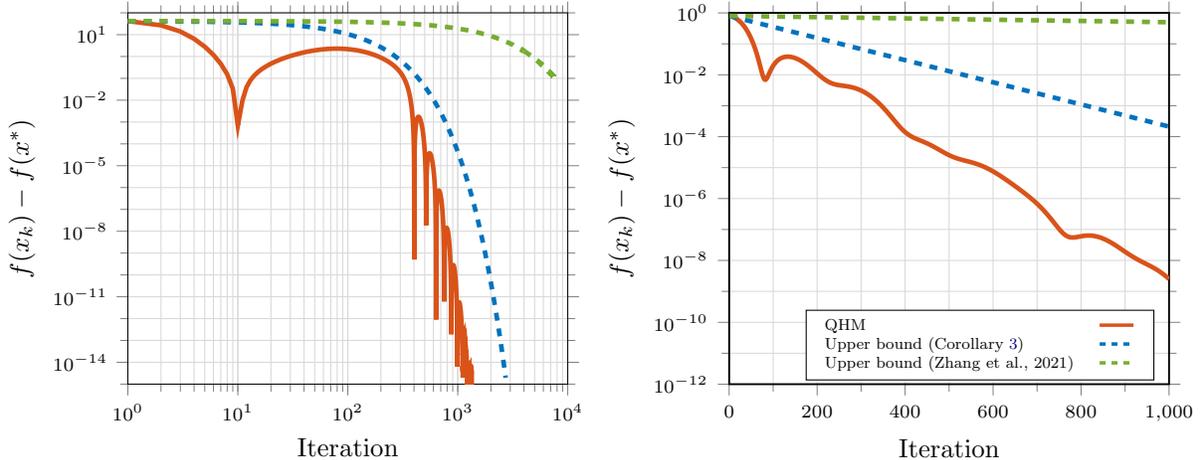

    \centering
\definecolor{mycolor4}{rgb}{0.00000,0.44700,0.74100}
\definecolor{mycolor3}{rgb}{0.85000,0.32500,0.09800}
\definecolor{mycolor2}{rgb}{0.92900,0.69400,0.12500}
\definecolor{mycolor1}{rgb}{0.49400,0.18400,0.55600}
\definecolor{mycolor5}{rgb}{0.46600,0.67400,0.18800}

\caption{Discrete-time simulation; comparison of QHM method performance under various settings with existing upper bounds for (a) $f(x)=4(L-\mu)\log(1+e^{-x})+\frac{\mu}{2}x^2$, $L=10,\mu=10^{-3}$ and (b) 10-dimensional regularized binary classification problem with logistic loss for random data and labels of length 1000. The regularization parameter is chosen $\mu = 10^{-3}$. QHM parameters are chosen so that the best possible theoretical convergence rate is achieved.}
    \label{fig2_new}
\end{figure}

\begin{remark}
    \Cref{theorem4} provides various convergence rates depending on the choice of parameters. 
    A natural question is whether a faster convergence rate than the rate of NAG can be achieved. 
    Setting $s=\left({(1-\rho)}/{q}\right)^2$ so that $\rho=1-q\sqrt{s}$, and substituting into the conditions of \Cref{theorem4} yields $\rho\geq \max\left\{1-{m}/{p},1-{q}/{(mL)}\right\}$. 
    These two expressions are inversely dependent on $m$. 
    Thus, their equality determines the optimal choice, $m=\sqrt{{(qp)}/{L}}$. 
    Substituting this value of $m$ in the bounds gives $\rho\geq 1-\sqrt{{q}/{(pL)}}$. 
    Since $q/p\leq\mu$, the right-hand side cannot be smaller than $1-\sqrt{{\mu}/{L}}$.
\end{remark}

\subsubsection{New Results for the QHM Method}

QHM \citep{ma2018quasihyperbolic} is an accelerated method commonly used in machine learning, particularly for training neural networks. 
A rescaled version of this numerical scheme is as follows:
\begin{align} \label{QHM} \tag{QHM}
    \left\{
    ~~
    \begin{aligned}
        x_{k+1}-x_k & =   -s(1-a)\nabla f(x_k)- s a  g_{k+1}, \\
        g_{k+1} & = bg_k+ \nabla f(x_{k}),
    \end{aligned}
    \right.
\end{align}
where $0 \leq a,b\leq 1$. The following corollary, derived from \Cref{theorem4}, provides new convergence guarantees for the QHM method:

\begin{corollary} \label{coro_QHM}
    With the choice of $m=(1-a)\sqrt{s}$, $n=q=\sqrt{a\mu}$, $p=\frac{a}{q}+\sqrt{s}$, $q\sqrt{s}\leq \frac{1}{2}$, and $s\leq \frac{4}{3L}$, for some $0<a\leq \frac{1}{4}$, algorithm \eqref{SIE_newalg} reduces to \eqref{QHM} with $b=\frac{1-q\sqrt{s}}{1+q\sqrt{s}}$. 
    Furthermore, invoking \Cref{theorem4} with these parameters yields the following convergence guarantee:
    \begin{align}
        f(x_{k})-f(x^*)\leq C_{QHM}(1-\sqrt{a \mu s})^k,
    \end{align}
where $C_{QHM}$ is a positive constant.
\end{corollary}

The above corollary is proved in \Cref{proof_QHM}. 
To our knowledge, the previous best convergence rate guarantee for QHM was $\mathcal{O}((1+\sqrt{1/\kappa}/40)^{-k})$), but this result was established only for $\kappa\geq 9$, where $\kappa$ is the condition number \citep{zhang2021revisiting}. 
In \Cref{coro_QHM}, we not only improve this convergence rate but also remove the restrictive assumption on $\kappa$. 
We demonstrate this improvement with numerical examples in \Cref{fig2_new}. 

\begin{remark}
    The general method, \eqref{SIE_newalg}, recovers other accelerated methods beyond those listed in \Cref{table_compare2}. 
    One such example is \eqref{TM_method}. 
    However, \Cref{theorem4} does not recover the convergence rate $\mathcal{O}\left((1-\sqrt{\tfrac{1}{\kappa}})^{2k}\right)$ associated with \eqref{TM_method}. 
    We believe that a more refined Lyapunov analysis could address this gap. 
\end{remark}

Finally, we conclude this section with \Cref{table_compare3}, which summarizes previously known convergence rates related to this work alongside the rates reported here.

\begin{table}
  \caption{Comparison of the convergence rates proposed in this work with existing results.}
  \label{table_compare3}
  \centering
  \resizebox{\textwidth}{!}{
  \begin{tabular}{c c c c c c}
    \toprule
    \multicolumn{1}{c}{ }&\multicolumn{1}{c}{\citep{Shi2021UnderstandingTA}}   &     \multicolumn{1}{c}{\citep{zhang2021revisiting}} &      \multicolumn{1}{c}{\citep{DBLP:conf/cdc/SunGK20}} &\multicolumn{1}{c}{H-NAG\citep{chen2019first}} &\multicolumn{1}{c}{This work} \\     
    \midrule
   \multirow{1}{*}{\ref{NAG_ODE}} & 
 $\mathcal{O}\left(e^{-\tfrac{\sqrt{\mu}t}{4}}\right)$ & $\mathcal{O}\left(e^{-\tfrac{\sqrt{\mu}t}{2}}\right)$
&  -& - & $\mathcal{O}\left(e^{-\sqrt{\mu}t}\right)$  \\

 \multirow{1}{*}{\ref{HR_TM}}  & 
 - & -&$\mathcal{O}\left(e^{-\tfrac{\sqrt{M}t}{2}}\right)$
& -
&  $\mathcal{O}\left(e^{-\tfrac{4\sqrt{M}t}{3}}\right)$   \\

\multirow{1}{*}{NAG Algorithm}  & 
 $\mathcal{O}\left((1+\tfrac{1}{9}\sqrt{\tfrac{1}{\kappa}})^{-k}\right)$ & $\mathcal{O}\left((1+\tfrac{1}{30}\sqrt{ \tfrac{1}{\kappa}})^{-k}\right)$ 
&  -& 
$\mathcal{O}\left((1+\sqrt{\tfrac{\sqrt{\kappa}-1}{\kappa\sqrt{\kappa}}})^{-k}\right)$
& $\mathcal{O}\left((1-\sqrt{ \tfrac{1}{\kappa}})^{k}\right)$  \\
 \multirow{1}{*}{\ref{QHM}}  & 
  - & $\mathcal{O}\left((1+\tfrac{1}{40}\sqrt{\tfrac{1}{\kappa}})^{-k}\right),\kappa\geq 9$ 
& -
&  -&$\mathcal{O}\left((1-\sqrt{\tfrac{1}{3\kappa}})^{k}\right)$    \\
    \bottomrule
  \end{tabular}
  }
\end{table}

\section{Implications of the HR-ODEs for Convex Functions}\label{section3}

In this section, we examine the HR-ODEs for convex functions defined in equations \eqref{HR_general_ODE_F_laborde} and \eqref{HR_general_ODE_F_Shi}. 

\subsection{Gradient Norm Minimization of NAG}\label{subsec_gradientnorm}

Reformulating \eqref{HR_general_ODE_F_Shi} with the parameters $\alpha_t$, $\beta_t$ and $\gamma_t$ defined in (\ref{prams_general}) gives
\begin{align} \label{Two_oneline_ODE_perturbed}
    \left\{
    ~~
    \begin{aligned}
        \dot{X_t} &= n(t)(V_t-X_t)-\sqrt{s}\nabla f(X_t)\\
        \dot{V_t} &= -q(t)\nabla f(X_t) - \sqrt{s}\frac{\dot q(t)n(t)-\dot n(t)q(t)}{n^2(t)q(t)} \nabla f(X_t).
    \end{aligned}
    \right.
\end{align}
Applying the SIE discretization to this ODE for $X_t \approx X(t_k)$, $V_t\approx V(t_k)$, and setting $n(t_k)=p/t_k$, $q(t_k)=Cpt_k^{p-1}$ with $p=2$, $t_k=k\sqrt{s}$, and $C=1/4$, we obtain 
\begin{align} \label{new_algorithm}
    \left\{
    ~~
    \begin{aligned}
        x_{k+1} &= x_{k} + \tfrac{2}{k}(v_k-x_{k+1})-{s}\nabla f(x_k),\\
        v_{k+1} &= v_k -\tfrac{1}{2} sk \nabla f(x_{k+1})-s\nabla f(x_{k+1}), 
    \end{aligned}
    \right.
\end{align}
which exactly recovers the NAG algorithm. 

To our knowledge, the interpretation of NAG as the SIE discretization of \eqref{Two_oneline_ODE_perturbed} is new (see \citep{ahn2022understanding} for other four common representations of NAG), and this precise connection motivated our choice of the Lyapunov function in the next result, leading to an improved convergence rate.

\begin{theorem}\label{theorem_GDnorm}
    Let $f\in\mathcal{F}_L$, and  consider the algorithm defined in \eqref{new_algorithm}, initialized with $v_0=x_1 \in \mathbb{R}^d$ and a step-size $0 < s\leq 1/L$. 
    Then, the sequence $x_k$ satisfies 
    \begin{gather}
        \frac{k^3s^2}{12}\min_{0\leq i\leq k-1}\|\nabla f(x_i)\|^2  \leq \|x_1-x^*\|^2 \label{eqn:Lyapunocv-convex-thm-1}, \\
         \frac{sk(k+2)}{2}f(x_k)-f(x^*) \leq\|x_1-x^*\|^2. \label{eqn:Lyapunocv-convex-thm-2}
    \end{gather}
    
\end{theorem}

\begin{proof}[Proof of \Cref{theorem_GDnorm}]
The proof is based on the discrete Lyapunov analysis with the function
\begin{align} \label{Lyapunov_function}
    \varepsilon(k) = \frac{s(k+2)k}{4}(f(x_k)-f(x^*)) + \frac{1}{2}\big\|x_{k+1}-x^*+\tfrac{k}{2}(x_{k+1}-x_k)+\tfrac{ks}{2}\nabla f(x_k)\big\|^2.
\end{align}
$\varepsilon(k)$ is clearly non-negative. 
We show that, under the conditions listed in \Cref{theorem_GDnorm}, it decreases as
\begin{align} \label{eqn:Lyapunov-bound}
    \varepsilon(k)-\varepsilon(0)\leq -\frac{k^3s^2}{24}\min_{0\leq i\leq k}\|\nabla f(x_i)\|^2
\end{align}
See \Cref{thm5_proof} for the complete proof. 

By combining \eqref{Lyapunov_function} and \eqref{eqn:Lyapunov-bound}, we obtain 
\begin{align*}
    \epsilon(0)=\frac{1}{2}\|x_1 - x^*\|^2 
    \geq \frac{k^3s^2}{24}\min_{0\leq i\leq k}\|\nabla f(x_i)\|^2,
\end{align*}
which yields the result in \eqref{eqn:Lyapunocv-convex-thm-1}. 
Similarly, we have 
\begin{align*}
    \frac{sk(k+2)}{4}(f(x_k)-f(x^*)) 
    \leq \epsilon(k) 
    \leq \epsilon(0) 
    = \frac{1}{2}\|x_1 - x^*\|^2 ,
\end{align*}
which gives the bound in \eqref{eqn:Lyapunocv-convex-thm-2}. 
\end{proof}

\begin{remark}
    The rate in \Cref{theorem_GDnorm} improves upon the previously established rate in  \citep{Shi2021UnderstandingTA}, which is given by 
    \begin{align}
        \min_{0\leq i\leq k}\|\nabla f(x_i)\|^2 \leq \frac{8568}{(k+1)^3s^2}\|x_0-x^*\|^2,   
    \end{align}
    for $0< s\leq 1/(3L)$ and $k\geq0$. 
    This improvement not only reduces the constant but also relaxes the conditions on the step-size. 
\end{remark}

\begin{remark}
    Concurrently with our prior work \citep{maskan2024variational}, \citet{chen2022gradient} demonstrated a similar convergence rate using an \textit{implicit velocity} perspective on HR-ODE.
\end{remark}

\subsection{Rate-Matching Approximates the NAG  Algorithm}\label{subsecsec_rate}

In this section, we establish a connection between the rate-matching technique proposed in \citep{WibisonoE7351} and the NAG algorithm.
Specifically, applying the rate-matching technique to the LR-ODE formulation of NAG from \citep{JMLR:v17:15-084} results in an algorithm that exhibits behavior similar to \eqref{HR_cnts_general_npq2}. 
Using this observation, we then approximately recover the NAG algorithm through rate-matching discretization.

\citet{WibisonoE7351} showed that discretizing the LR-ODE of the NAG method using the rate-matching technique gives an accelerated method:
\begin{align} \label{rate_match_2}
    \left\{
    ~~
    \begin{aligned}
        x_{k+1} & =\tfrac{2}{k+2}z_k+\tfrac{k}{k+2}y_k,\\
        y_{k} & =x_k-s\nabla f(x_k),   \\
        z_{k} & =z_{k-1} -\tfrac{1}{2} sk \nabla f(y_k),
    \end{aligned}
    \right.
\end{align}
which has a convergence rate of $\mathcal{O}(1/(s k^2))$. 
In the following proposition, we analyze the behavior of (\ref{rate_match_2}) in the limit as $s \rightarrow 0$.
  
\begin{proposition}\label{prop1}
    The continuous-time behavior of (\ref{rate_match_2}) is given by
    \begin{align} \label{cont_rate_match_2}
        \Ddot{X}_t+\left(\frac{3}{t}+\sqrt{s}\nabla^2 f(X_t)\right)\dot{X}_t+\left(1+\frac{\sqrt{s}}{2t}\right)\nabla f(X_t)=0.
    \end{align}
    We also consider the following modified HR-ODE:
    \begin{align} \label{cont_rate_match_2_mod}
        \Ddot{X}_t+\left(\frac{3}{t}+\sqrt{s}\nabla^2 f(X_t)\right)\dot{X}_t+\left(1+\frac{\sqrt{s}}{t}\right)\nabla f(X_t)=0,
    \end{align}
    which is obtained by replacing the last term in \eqref{cont_rate_match_2} using the approximation $(k+1)/(k+1/2)\approx (k+1)/(k)$. 
    A similar modification was made in \citep[Section 4]{shi2019acceleration}. 
\end{proposition}

\begin{proof}[Proof of \Cref{prop1}]
    The proof follows from a one-line representation of \eqref{rate_match_2}. 
    Taking $X(t)\approx X(t_k)=x_k$ and $Y(t)\approx Y(t_k)=y_k$ at $t_k=(k+1/2)\sqrt{s}$, and using Taylor approximations, we obtain \eqref{cont_rate_match_2}. 
    See \Cref{proof_prop1} for the detailed proof. 
\end{proof}

The ODE (\ref{cont_rate_match_2_mod}) coincides with \eqref{HR_cnts_general_npq2} for $p=2,C=1/4$. 
This ODE corresponds to a perturbed LR-ODE for the NAG method. 
Indeed, it can be rewritten as
\begin{align} \label{cont_rate_match_2_perturb}
    \left\{
    ~~
    \begin{aligned}
         \dot{X}_t & = \tfrac{2}{t}(Z_t-X_t)-\sqrt{s}\nabla f(X_t),  \\
          \dot Z_t & = -\tfrac{t}{2}\nabla f(X_t).
    \end{aligned}
    \right.
\end{align}
which can be viewed as a perturbation of the LR-ODE
\begin{align} \label{cont_rate_match_2_perturb2}
    \left\{
    ~~
    \begin{aligned}
        \dot{X}_t & = \tfrac{2}{t} (Z_t-X_t),  \\
        \dot Z_t & = -\tfrac{t}{2} \nabla f(X_t).
    \end{aligned}
    \right.
\end{align}
In this sense, rate-matching implicitly perturbs the LR-ODE. 

A natural question arises: can we recover the HR-ODE \eqref{HR_cnts_general_stable_ODE} (which corresponds to the NAG algortihm through the SIE discretization) from the rate-matching technique? 
To address this, we first introduce a perturbation in the second line of \eqref{cont_rate_match_2_perturb2}, then apply the rate-matching discretization. 
Perturbing \eqref{cont_rate_match_2_perturb2} yields
\begin{align} \label{cont_rate_match_2_perturb3}
    \left\{
    ~~
    \begin{aligned}
        \dot{X}_t & = \tfrac{2}{t}(Z_t-X_t),  \\
        \dot Z_t & = -\tfrac{t}{2}\nabla f(X_t)-\sqrt{s}\nabla f(X_t).
    \end{aligned}
    \right.
\end{align}
Discretizing \eqref{cont_rate_match_2_perturb3} using the rate-matching technique with $t_k=k\sqrt{s}$ gives
\begin{align} \label{rate_match_3}
    \left\{
    ~~
    \begin{aligned}
        x_{k+1} &=\tfrac{2}{k+2}z_k+\tfrac{k}{k+2}y_k,\\
        y_{k} &=x_k-s\nabla f(x_k),   \\
        z_{k} &=z_{k-1} -\tfrac{s (k+2)}{2} \nabla f(y_{k}),
    \end{aligned}
    \right.
\end{align}
which closely resembles the NAG algorithm. 
Indeed, replacing $\nabla f(y_k)$ with $\nabla f(x_k)$ in the third line of \eqref{rate_match_3} gives exactly the NAG algorithm. 

Typically, $x_k$ and $y_k$ remain close in practice, since their difference is proportional to $s$. 
In continuous-time, as $s\rightarrow 0$, we have $X(t_k)\approx Y(t_k)$, implying that the HR-ODE of (\ref{rate_match_3}) is (\ref{HR_cnts_general_stable_ODE}). 
Thus, the corresponding HR-ODE for \eqref{rate_match_3} is given by
\begin{align} \label{cont_rate_match_2_perturb4}
    \left\{
    ~~
    \begin{aligned}
        \dot{X}_t & = \tfrac{2}{t}(Z_t-X_t)-\sqrt{s}\nabla f(X_t),  \\
        \dot Z_t  & = -\tfrac{t}{2}\nabla f(X_t)-\sqrt{s}\nabla f(X_t).
    \end{aligned}
    \right.
\end{align}
which is a perturbed version of (\ref{cont_rate_match_2_perturb3}) and the HR-ODE associated with the NAG algorithm.

\section{Conclusion and Future Directions}

In this work, we developed a unified HR-ODE framework for analyzing accelerated gradient methods in smooth convex and smooth strongly convex optimization. 
Using the Forced Euler-Lagrange equation, we derived a general HR-ODE that captures the continuous-time behavior of many important momentum-based methods. 
Through this framework, we recovered existing HR-ODEs for these methods or introduced new HR-ODE formulations that provide tighter convergence guarantees. 
We analyzed the convergence rate of the proposed HR-ODEs and their discretizations through generic Lyapunov functions. 

In the strongly convex regime, our framework revealed a consistent connection between the proposed general HR-ODE \eqref{GM2-ODE} and its SIE discretization. 
The proposed general method precisely recovers several important accelerated methods such as HB, NAG, TM, and QHM. 
Among existing general formulations, we showed that \eqref{GM2-ODE} achieves the fastest rate of convergence for each recovered algorithm. 
Furthermore, in the convex regime, our approach led to improved gradient norm minimization rates for NAG. 
Additionally, our unified framework captures the HR-ODE of the accelerated method proposed by \cite{WibisonoE7351}. 
We also demonstrated that NAG can be viewed as an approximation of the rate-matching discretization applied to a specific ODE. 

We highlight two key directions for future work. 
First, in our application of the Forced Euler-Lagrange equation, the proposed forces are defined only in Euclidean space. 
Extending this framework to non-Euclidean spaces would be a natural next step. 
Second, in the strongly convex setting, although the proposed general model (\ref{SIE_newalg}) recovers the TM method, our discrete-time analysis does not achieve the state-of-the-art convergence rate of TM. 
While we can achieve faster rates for TM by using \eqref{TM_new_high_res2}, this formulation relies on an approximation step. 
It remains open whether similar convergence rates can be shown for the precise formulation in \eqref{TM_new_high_res}. 

\section*{Acknowledgments}
This work received support from the Wallenberg AI, Autonomous Systems and Software Program (WASP) funded by the Knut and Alice Wallenberg Foundation.

\bibliographystyle{icml2021}
\bibliography{bibliography}

\newpage
\appendix

\section{Appendix}

For better accessibility, we divided the supplementary material into 8 sections.
In Section \ref{app1}, we present only the proofs of the results in the main paper.
In Section \ref{sec:gmvsgm2} \eqref{GM-ODE} is compared with \eqref{GM2-ODE}. Exact recovery of NAG through SIE discretization is discussed in Section \ref{sec:SIE_rec_NAG}. In Sections \ref{app2} and \ref{app3}, we conducted continuous-time and discrete-time convergence analysis of the proposed method \eqref{SIE_newalg} on quadratics.
We believe these results are insightful.
Section \ref{app_subsec_QCI} provides a short background on the convergence analysis of the optimization methods through dynamical systems.
Following the line of \citep{zhang2021revisiting}, Section \ref{appendix_EE} shows how EE discretization can lead to acceleration.
Finally, Section \ref{app6} provides additional numerical experiments.
\subsection{Proofs}\label{app1}
In the following section, we present the detailed proofs of the results proposed in the main body of our work. 

\subsubsection{Proof of \Cref{Theorem3_1}}\label{prf:thm3_1}

\noindent (a) Consider the Lyapunov function
\begin{align*} 
    \varepsilon(t) = e^{\beta_t}\left(\frac{Ce^{\alpha}}{2}\|X_t-x^*+e^{-\alpha}\dot X_t+e^{-\alpha}\sqrt{s}\nabla f(X_t)\|^2+f(X_t)-f(x^*)\right).
\end{align*}
Taking derivative with respect to time gives
\begin{align} \label{thm32_eqn_2}
    & \frac{d\varepsilon(t)}{dt} 
    = \dot \beta e^{\beta_t}\left(\frac{Ce^{\alpha}}{2}\|X_t-x^*+e^{-\alpha}\dot X_t+e^{-\alpha}\sqrt{s}\nabla f(X_t)\|^2+f(X_t)-f(x^*)\right) \\
    & \quad + C e^{\alpha+\beta_t}\left\langle \dot X_t  + e^{-\alpha}\Ddot{X}_t +e^{-\alpha}\sqrt{s}\nabla ^2f(X_t)\dot X_t  ,X_t-x^*+e^{-\alpha}\dot X_t+e^{-\alpha}\sqrt{s}\nabla f(X_t)\right\rangle+e^{\beta_t}\langle \nabla f(X_t),\dot X_t \rangle. \nonumber
\end{align}
Next, we will use (\ref{sc_eqn1}) in (\ref{thm32_eqn_2})
\begin{align*} 
    \frac{d\varepsilon(t)}{dt} 
    & = \dot \beta_t e^{\beta_t}\left(\frac{Ce^{\alpha}}{2}\Bigg[\|X_t-x^*\|^2+e^{-2\alpha}\|\dot X_t\|^2+\|\sqrt{s}e^{-\alpha}\nabla f(X_t)\|^2\right.\nonumber\\
    & \qquad +2\sqrt{s}e^{-\alpha}\langle X_t-x^*,\nabla f(X_t) \rangle+ 2\sqrt{s}e^{-2\alpha}\langle \nabla f(X_t),\dot X_t \rangle+2e^{-\alpha}\langle X_t-x^*,\dot X_t\rangle\Bigg]\Bigg)\nonumber\\
    & \qquad +C e^{\alpha+\beta_t}\Bigg[ (1-e^{-\alpha}(\dot \gamma_t+\dot \beta_t))\langle \dot X_t,X_t-x^*\rangle+ e^{-\alpha}(1-e^{-\alpha}(\dot \gamma_t+\dot \beta_t))\|\dot X_t\|^2 \nonumber\\
    & \qquad +\sqrt{s}e^{-\alpha}(1-e^{-\alpha}(\dot \gamma_t+\dot \beta_t))\langle \dot X_t,\nabla f(X_t) \rangle - (\frac{1}{C}+\sqrt{s}\dot \beta_te^{-\alpha})\langle \nabla f(X_t) ,X_t-x^*\rangle\nonumber\\
    & \qquad -(\frac{e^{-\alpha}}{C}+\sqrt{s}\dot\beta_te^{-2\alpha})\langle \nabla f(X_t),\dot X_t \rangle-(\frac{\sqrt{s}e^{-\alpha}}{C}+s\dot\beta_te^{-2\alpha})\|\nabla f(X_t)\|^2\Bigg] \nonumber\\
    & \qquad + e^{\beta_t}\langle \nabla f(X_t),\dot X_t \rangle+\dot \beta_t e^{\beta_t} (f(X_t)-f(x^*) )
\end{align*}
Now, using strong convexity of $f$ and applying  $0\leq \dot \beta\leq \dot \gamma_t=e^{\alpha}\leq \mu/C$, gives
\begin{align*}
    \frac{d\varepsilon(t)}{dt}
    & \leq-\frac{C\dot \beta_t e^{ \beta_t-\alpha }}{2}\|\dot X_t\|^2-Ce^{\beta_t-\alpha  }\dot \beta_t\langle \dot X_t,\sqrt{s}\nabla f(X_t) \rangle - e^{\beta_t}\Big(\frac{1}{\sqrt{s}}+\frac{Ce^{-\alpha}\dot \beta_t}{2}\Big)\|\sqrt{s}\nabla f(X_t)\|^2 \nonumber\\
    & \qquad + e^{\beta_t}(\dot\beta_t-e^{\alpha})(f(X_t)-f(x^*)) + \left(\frac{C\dot\beta_te^{\alpha+\beta_t}}{2}-\frac{\mu e^{\alpha+\beta_t}}{2}\right)\|X_t-x^*\|^2\nonumber\\
    & \leq -\frac{C\dot \beta_t e^{ \beta_t-\alpha}}{2} \|   \dot X_t+  \sqrt{s}\nabla f(X_t)\|^2 
    \leq 0.
\end{align*}
Therefore, $e^{\beta_t}(f(X_t)-f(x^*))\leq \varepsilon(t)\leq \varepsilon(0)$ and the proof is complete.

\vspace{1em}

\noindent (b) Consider the Lyapunov function
\begin{align*}
    \begin{aligned}
        \varepsilon(t) = e^{\gamma_t}\Bigg(\frac{Ce^{\alpha_t}}{2}\|X_t-x^*&+e^{-\alpha_t}\dot X_t+\sqrt{s}e^{-\alpha}\nabla f(X_t)\|^2 - C\Big(\frac{e^{\alpha_t}-\dot\beta_t}{2}\Big)\|X_t-x^*\|^2+f(X_t)-f(x^*)\Bigg).
    \end{aligned}
\end{align*}
    We start by taking derivative with respect to time (note that $\alpha_t = \alpha$ and $\Ddot \beta=0$ )
\begin{align} \label{thm32_eqn_5}
    \frac{d\varepsilon(t)}{dt} 
    & = \dot \gamma e^{\gamma_t}\Bigg(\frac{Ce^{\alpha_t}}{2} \|X_t-x^*+e^{-\alpha_t}\dot X_t+\sqrt{s}e^{-\alpha}\nabla f(X_t)\|^2 - C\Big(\frac{e^{\alpha_t}-\dot\beta_t}{2}\Big)\|X_t-x^*\|^2+f(X_t)-f(x^*)\Bigg)\nonumber\\
    & \qquad +e^{\gamma_t}\Bigg(\langle \dot X_t,\nabla f(X_t) \rangle - C\left(e^{\alpha}-\dot \beta_t\right) \langle \dot X_t,X_t-x^*\rangle \\
    & \qquad + Ce^{\alpha}\langle e^{-\alpha}\Ddot{X}_t+\dot X_t+\sqrt{s}e^{-\alpha}  \nabla^2 f(X_t)\dot X_t,X_t-x^*+e^{-\alpha}\dot X_t+\sqrt{s}e^{-\alpha}\nabla f(X_t)\rangle \Bigg).\nonumber
\end{align}
Using the HR-ODE (\ref{sc_eqn1}) we know that
\begin{align*}
    e^{-\alpha}\Ddot{X}_t+\dot X_t+\sqrt{s}e^{-\alpha} \nabla^2 f(X_t)\dot X_t = -\left(e^{-\alpha}\dot \beta_t\dot X_t +\left(\frac{1}{C}+\sqrt{s}e^{-\alpha}\dot\beta_t\right)\nabla f(X_t)\right).
\end{align*}
Replacing in \eqref{thm32_eqn_5} gives
\begin{align*} 
    \frac{d\varepsilon(t)}{dt}
    & = \dot \gamma e^{\gamma_t}\Bigg(\frac{Ce^{\alpha}}{2}\|X_t-x^*+e^{-\alpha}\dot X_t+\sqrt{s}e^{-\alpha}\nabla f(X_t)\|^2 - C\Big(\frac{e^{\alpha}-\dot\beta_t}{2}\Big)\|X_t-x^*\|^2+f(X_t)-f(x^*)\Bigg) \nonumber \\
    & \qquad +e^{\gamma_t}\Bigg(\langle \dot X_t,\nabla f(X_t) \rangle - C\left(e^{\alpha}-\dot \beta_t\right)\langle \dot X_t,X_t-x^*\rangle- C\dot\beta_t\langle \dot X_t,X_t-x^* \rangle  \\
    & \qquad - C\dot \beta_t e^{-\alpha}\|\dot X_t\|^2- C\sqrt{s}\dot \beta_te^{-\alpha} \langle \dot X_t,\nabla f(X_t)\rangle - (1+C\sqrt{s}\dot\beta_te^{-\alpha})\langle \nabla f(X_t) ,\dot X_t \rangle \nonumber\\
    & \qquad -(e^{\alpha}+\sqrt{s}\dot\beta_tC)\langle \nabla f(X_t),X_t-x^* \rangle - \sqrt{s}(1+\sqrt{s}C\dot\beta_te^{-\alpha})\|\nabla f(X_t)\|^2 \Bigg).\nonumber
\end{align*}
After simplification and using the assumption $\dot\gamma_t = e^{\alpha}$ we have
\begin{align*} 
    \begin{aligned}
        \frac{d\varepsilon(t)}{dt} =  &e^{\gamma_t}\Bigg(e^{\alpha}(f(X_t)-f(x^*))+\frac{Ce^{\alpha}\dot \beta_t}{2}\|X_t-x^*\|^2\\
        & +Ce^{-\alpha}(\frac{e^{\alpha}}{2}-\dot \beta_t)\|\dot X_t\|^2 + \left(\frac{Cs}{2}-\sqrt{s}-sC\dot\beta_te^{-\alpha} \right)\|\nabla f(X_t)\|^2   \\
        & + \left(C\sqrt{s}(e^{\alpha}-\dot\beta_t)-e^{\alpha}\right)\langle\nabla f(X_t),X_t-x^*\rangle + C\sqrt{s}e^{-\alpha}2\left(\frac{e^{\alpha}}{2}-\dot\beta_t\right)\langle\dot X_t,\nabla f(X_t)\rangle\Bigg).
    \end{aligned}
\end{align*} 
Using assumptions $-\dot\beta_t\leq -e^{\alpha}/2$ and $\sqrt{s}C\leq 2$, we have
\begin{align*}
    \left(C\sqrt{s}(e^{\alpha}-\dot\beta_t)-e^{\alpha}\right) \leq 0.
\end{align*}
This together with strong convexity of $f$ gives: 
\begin{align} \label{thm32_eqn_8}
    \begin{aligned}
        \frac{d\varepsilon(t)}{dt} 
        & \leq e^{\gamma_t}\Bigg((e^{\alpha}+C\sqrt{s}e^{\alpha}-e^{\alpha}-C\sqrt{s}\dot\beta )(f(X_t)-f(x^*))  +\left[\frac{Ce^{\alpha}\dot\beta_t}{2}-\frac{\mu e^{\alpha}}{2}\right]\|X_t-x^*\|^2   \\
        & \qquad +Ce^{-\alpha}\left(·\frac{e^{\alpha}}{2}-\dot\beta_t\right)\|\dot X_t\|^2 + C\sqrt{s}e^{-\alpha}2\left(\frac{e^{\alpha}}{2}-\dot\beta_t\right)\langle\dot X_t,\nabla f(X_t)\rangle\\
        & \qquad +\frac{\mu C\sqrt{s}}{2}(e^{\alpha}-\dot \beta_t)\|X_t-x^*\|^2+ \left(\frac{Cs}{2}-\sqrt{s}-sC\dot\beta_te^{-\alpha} \right)\|\nabla f(X_t)\|^2 \Bigg)\\
        & =e^{\gamma_t}\Bigg((C\sqrt{s}e^{\alpha}-C\sqrt{s}\dot\beta )(f(X_t)-f(x^*))\\
        &\qquad  +\left[\frac{Ce^{\alpha}\dot\beta_t}{2}-\frac{\mu e^{\alpha}}{2}\right]\|X_t-x^*\|^2 +\frac{\mu C\sqrt{s}}{2}(e^{\alpha}-\dot \beta_t)\|X_t-x^*\|^2  \\
        &\qquad +Ce^{-\alpha}\left(\frac{e^{\alpha}}{2}-\dot\beta_t\right)\left[\|\dot X_t\|^2 + 2\langle \sqrt{s}\nabla f(X_t),\dot X_t\rangle + \|\sqrt{s}\nabla f(X_t)\|^2\right] - \sqrt{s}\|\nabla f(X_t)\|^2 \Bigg).
    \end{aligned}
\end{align}    
Now, note that 
\begin{align*}
    -C\sqrt{s}\dot\beta (f(X_t)-f(x^*))\leq -\frac{\mu C\sqrt{s}\dot\beta }{2}\|X_t-x^*\|^2,
\end{align*}
Simplifying (\ref{thm32_eqn_8}) gives
\begin{align} \label{thm32_eqn_9}
    \begin{aligned}
        \frac{d\varepsilon(t)}{dt} \leq   &e^{\gamma_t}\Bigg(C\sqrt{s}e^{\alpha}(f(X_t)-f(x^*))+\frac{\mu C\sqrt{s}}{2}(e^{\alpha}-2\dot \beta_t)\|X_t-x^*\|^2 + \left[\frac{Ce^{\alpha}\dot\beta_t}{2}-\frac{\mu e^{\alpha}}{2}\right]\|X_t-x^*\|^2 \\
        &+Ce^{-\alpha}\left(\frac{e^{\alpha}}{2}-\dot\beta_t\right)\left[\|\dot X_t + \sqrt{s}\nabla f(X_t)\|^2 \right]- \sqrt{s}\|\nabla f(X_t)\|^2
        \Bigg).
    \end{aligned}
\end{align}     
Further,  if $\dot \beta_te^{-\alpha}\geq 1/2$ we get
\begin{align*}
    Ce^{-\alpha}\left(\frac{e^{\alpha}}{2}-\dot\beta_t\right)\left[\|\dot X_t + \sqrt{s}\nabla f(X_t)\|^2 \right]\leq 0, \quad \frac{\mu C\sqrt{s}}{2}(e^{\alpha}-2\dot \beta_t)\|X_t-x^*\|^2\leq 0
\end{align*}
Also, note that due to strong convexity we can write
\begin{align*}
    - \sqrt{s} \|\nabla f(X_t)\|^2\leq -2\mu\sqrt{s}\left(f(X_t)-f(x^*)\right).
\end{align*}
Utilizing these in (\ref{thm32_eqn_9}) we get
\begin{align} \label{thm32_eqn_10}
    \begin{aligned}
        \frac{d\varepsilon(t)}{dt} \leq   &e^{\gamma_t}\Bigg(C\sqrt{s}e^{\alpha}(f(X_t)-f(x^*))-2\mu\sqrt{s}\left(f(X_t)-f(x^*)\right)
        \Bigg)\leq 0,
    \end{aligned}
\end{align} 
where the last inequality holds due to $\dot \beta_te^{-\alpha}\geq 1/2$ and $C\dot \beta_t\leq\mu$. This gives
\begin{align*}
    e^{\gamma_t}\Bigg(\frac{Ce^{\alpha}}{2}\|X_t-x^*+e^{-\alpha_t}\dot X_t+e^{-\alpha}\sqrt{s}\nabla f(X_t)\|^2\Bigg)\leq \varepsilon(t) \leq \varepsilon(0),
\end{align*}
and the proof is complete. 

\subsubsection{Proof of \Cref{Theorem_ODE_laborde}}\label{thm1_proof}
Consider the Lyapunov function 
\begin{align*} 
    \varepsilon(t)=\frac{1}{2}\|X_t+e^{-\alpha_t}\dot X_t-x^*+\sqrt{s}e^{-\alpha_t}\nabla f(X_t)\|^2+e^{\beta_t}(f(X_t)-f(x^*)).
\end{align*}
Taking derivative with respect to $t$ gives
\begin{align}\label{prf_thm1_eqn1}
    \begin{aligned}
         \frac{d \varepsilon}{dt}=&\langle \frac{d}{dt}(X_t+e^{-\alpha_t}\dot X_t-x^*+\sqrt{s}e^{-\alpha_t}\nabla f(X_t)),X_t+e^{-\alpha_t}\dot X_t-x^*+\sqrt{s}e^{-\alpha_t}\nabla f(X_t)\rangle\\
        & +\dot \beta_t e^{\beta_t}(f(X_t)-f(x^*))+e^{\beta_t}\langle \nabla f(X_t), \dot X_t\rangle.
    \end{aligned}
\end{align}
Note that (\ref{HR_general_ODE_F_laborde}) can be represented as
\begin{align} \label{HR_general_der_format_laborde}
    \frac{d}{dt}\left[X_t+e^{-\alpha_t}\dot X_t+\sqrt{s}e^{-\alpha_t}\nabla f(X_t)\right]=-e^{\alpha_t+\beta_t}\nabla f(X_t).
\end{align}
Using (\ref{HR_general_der_format_laborde}) in (\ref{prf_thm1_eqn1}) we have
\begin{align}
     \frac{d \varepsilon}{dt} 
      =& \langle -e^{\alpha_t+\beta_t}\nabla f(X_t),X_t+e^{-\alpha_t}\dot X_t-x^*+\sqrt{s}e^{-\alpha_t}\nabla f(X_t) \rangle +\dot \beta_t e^{\beta_t}(f(X_t)-f(x^*))+e^{\beta_t}\dot X_t\nabla f(X_t)\nonumber\\
     =& -e^{\alpha_t+\beta_t}\langle\nabla f(X_t),X_t-x^*\rangle -e^{\beta_t}\langle \nabla f(X_t),\dot X_t\rangle -\sqrt{s}e^{\beta_t}\|\nabla f(X_t)\|^2\nonumber\\
     &  +\dot \beta_te^{\beta_t}(f(X_t)-f(x^*))+e^{\beta_t}\langle \nabla f(X_t),\dot X_t \rangle\nonumber\\
      \leq & -e^{\alpha_t+\beta_t}(f(X_t)-f(x^*))+\dot \beta_te^{\beta_t}(f(X_t)-f(x^*))\nonumber\\
      =& -e^{\beta_t}\left[(e^{\alpha_t}-\dot \beta_t)(f(X_t)-f(x^*))\right],\nonumber
\end{align}
where the inequality is due to convexity. Utilizing the ideal scaling condition $\dot \beta_t\leq e^{\alpha_t}$ we have
\begin{align}
     \frac{d \varepsilon}{dt}\leq 0.\nonumber
\end{align}
and the proof is complete.

\subsubsection{Proof of \Cref{Theorem_ODE_Shi}}\label{thm3_proof}
Consider the Lyapunov function 
\begin{align*} 
    \varepsilon(t)=&\frac{1}{2}\|X_t+e^{-\alpha_t}\dot X_t-x^*+\sqrt{s}e^{-\alpha_t}\nabla f(X_t)\|^2+(e^{\beta_t}+\sqrt{s}e^{-2\alpha_t}\dot \beta_t)(f(X_t)-f(x^*)).
\end{align*}
Taking derivative with respect to $t$ gives
\begin{align} \label{prf_thm3_eqn1}
\begin{aligned}
    \frac{d \varepsilon}{dt}=&\langle \frac{d}{dt}(X_t+e^{-\alpha_t}\dot X_t-x^*+\sqrt{s}e^{-\alpha_t}\nabla f(X_t)),X_t+e^{-\alpha_t}\dot X_t-x^*+\sqrt{s}e^{-\alpha_t}\nabla f(X_t)\rangle\\
    & +(\dot \beta_te^{\beta_t}-\sqrt{s}(2\dot \alpha_t)e^{-2\alpha_t}\dot \beta_t +\sqrt{s}e^{-2\alpha_t}\Ddot{\beta_t})(f(X_t)-f(x^*))\\
    &+(e^{\beta_t}+\sqrt{s}e^{-2\alpha_t}\dot \beta_t)\dot X_t\nabla f(X_t).
\end{aligned}
\end{align}
Note that (\ref{HR_general_ODE_F_laborde}) can be represented as
\begin{align} \label{HR_general_der_format_Shi}
    \frac{d}{dt}\left[X_t+e^{-\alpha_t}\dot X_t+\sqrt{s}e^{-\alpha_t}\nabla f(X_t)\right]=-\left(e^{\alpha_t+\beta_t}+\sqrt{s}e^{-\alpha_t}\dot \beta_t \right)\nabla f(X_t).
\end{align}
Using (\ref{HR_general_der_format_Shi}) in (\ref{prf_thm3_eqn1}) we have
\begin{align}
     \frac{d \varepsilon}{dt}=& \langle -\left(e^{\alpha_t+\beta_t}+\sqrt{s}e^{-\alpha_t}\dot \beta_t \right)\nabla f(X_t),X_t+e^{-\alpha_t}\dot X_t-x^*+\sqrt{s}e^{-\alpha_t}\nabla f(X_t) \rangle \nonumber \\
     &+(\dot \beta_te^{\beta_t}-\sqrt{s}(2\dot \alpha_t)e^{-2\alpha_t}\dot \beta_t +\sqrt{s}e^{-2\alpha_t}\Ddot{\beta_t})(f(X_t)-f(x^*))\nonumber\\
    &+(e^{\beta_t}+\sqrt{s}e^{-2\alpha_t}\dot \beta_t)\langle\nabla f(X_t),\dot X_t\rangle\nonumber\\
     =& -\left(e^{\alpha_t+\beta_t}+\sqrt{s}e^{-\alpha_t}\dot \beta_t \right)\langle\nabla f(X_t),X_t-x^*\rangle -\left(e^{\beta_t}+\sqrt{s}e^{-2\alpha_t}\dot \beta_t \right)\langle \nabla f(X_t),\dot X_t\rangle \nonumber\\
     &-\sqrt{s}\left(e^{\beta_t}+\sqrt{s}e^{-2\alpha_t}\dot \beta_t \right)\|\nabla f(X_t)\|^2\nonumber\\
     &+(\dot \beta_te^{\beta_t}-\sqrt{s}(2\dot \alpha_t)e^{-2\alpha_t}\dot \beta_t +\sqrt{s}e^{-2\alpha_t}\Ddot{\beta_t})(f(X_t)-f(x^*))\nonumber\\
    &+(e^{\beta_t}+\sqrt{s}e^{-2\alpha_t}\dot \beta_t)\langle\nabla f(X_t),\dot X_t\rangle\nonumber\\
     \leq &  -\left(e^{\alpha_t+\beta_t}+\sqrt{s}e^{-\alpha_t}\dot \beta_t \right)(f(X_t)-f(x^*))\nonumber\\
     &+(\dot \beta_te^{\beta_t}-\sqrt{s}(2\dot \alpha_t)e^{-2\alpha_t}\dot \beta_t +\sqrt{s}e^{-2\alpha_t}\Ddot{\beta_t})(f(X_t)-f(x^*))\nonumber\\
      =& -\left[e^{\beta_t}(e^{\alpha_t}-\dot \beta_t)+\sqrt{s}e^{-\alpha_t}(\dot \beta_t+2\dot \alpha_t e^{-\alpha_t}\dot \beta_t-e^{-\alpha_t}\Ddot{\beta}_t)\right](f(X_t)-f(x^*)).\nonumber
\end{align}
Utilizing the modified ideal scaling conditions $\dot \beta_t\leq e^{\alpha_t}$ and $\Ddot{\beta}_t\leq e^{\alpha_t}\dot \beta_t + 2\dot \alpha_t \dot \beta_t$ we have
\begin{align}
     \frac{d \varepsilon}{dt}\leq 0.\nonumber
\end{align}
and the proof is complete.

\subsubsection{Proof of \Cref{theorem2}}\label{prf:conti_conv}
    (a) Here, Theorem 6.4 in \citep{doi:10.1137/17M1136845} is used to find both a Lyapunov function and a convergence rate for (\ref{proposed_ODE_conti}). In order to apply the theorem we need to present the state space representation of (\ref{proposed_ODE_conti}). Taking $X,V\in\mathbb{R}^d,$
\begin{align*} 
    \quad\xi = [X,V]^T,\quad A=\left[\begin{array}{cc}
    -nI_d & nI_d \\
    qI_d & -qI_d
\end{array}\right],\quad B=\left[\begin{array}{c}
    -mI_d  \\
   - pI_d
\end{array}\right],\quad C=\left[\begin{array}{c}
    I_d  \\
   0_d
\end{array}\right]^T,
\end{align*} 
allows us to present the state space of (\ref{proposed_ODE_conti}) as
\begin{align} \label{statespace}
    \dot\xi(t)=A\xi(t)+Bu(t),\quad y(t)=C\xi(t),\quad u(t)=\nabla f(y(t))\quad \forall t\geq 0 ,
\end{align}
where $\xi\in\mathbb{R}^{d'}$ is the state, $y(t)\in \mathbb{R}^d$ $(d\leq d')$ is the output, and $u(t)$ is the continuous feedback input. Here, we would have $u^*=0$ and the fixed point of (\ref{statespace}) is
$A\xi^*=0,\quad y^*=C\xi^*.$ For simplicity, we set $\sigma=0$, i.e.  we remove any sign of $L$ from our formulation. Thus, the result holds for $\mu$-strongly convex functions. We need to find $P\succeq 0,\lambda> 0 $ such that $T\preceq 0$.
Consider
\begin{align*} 
 P=\hat{P}\otimes I_d,\quad \hat{P}= \left[\begin{array}{cc}
  p_{11}   &p_{12}  \\
    p_{12} & p_{22}
\end{array}\right],  \qquad  T=\hat{T}\otimes I_d,\quad \hat{T}= \left[\begin{array}{ccc}
  t_{11}   &t_{12}&t_{13}  \\
   t_{12}   &t_{22}&t_{23}  \\
     t_{13}   &t_{23}&t_{33} 
\end{array}\right].  
\end{align*}

\noindent Using the structure in Theorem 6.4 of \citep{doi:10.1137/17M1136845}, one can find the elements $t_{ij}$ as
\begin{align*}
\begin{array}{lrl}
     t_{11}= -(2n-\lambda)p_{11}+2qp_{12}-\frac{\lambda\mu}{2},& t_{12}&= np_{11}+qp_{22}-(n+q-\lambda)p_{12},  \\
     t_{13}= -mp_{11}-pp_{12}+\frac{\lambda-n}{2},& t_{22}&= 2(np_{12}-qp_{22})+\lambda p_{22},  \\
     t_{23}= -mp_{12}-pp_{22}+\frac{n}{2}, & t_{33}&= -m.
\end{array}
\end{align*}
Next step is to ensure that $\hat{P}\succeq 0$. If we take $\det(\hat{P})=0$ and find one of the diagonal elements such that it would be positive, then one of the eigenvalues of $\hat{P}$ is zero and the other one is positive which results in our favor. Before doing so, we will find $p_{11}$ and $p_{22}$ as a function of $p_{12}$. The latter is done by setting $t_{13}=0$ and $t_{23}=0$. Then
\begin{align*}
    \begin{aligned}
        & t_{13} = 0    \qquad  \Longrightarrow \quad p_{11}=\frac{-pp_{12}+\frac{\lambda-n}{2}}{m},  \\
        & t_{23} = 0    \qquad  \Longrightarrow \quad p_{22}=\frac{-mp_{12}+\frac{n}{2}}{p}.
    \end{aligned}
\end{align*}
These choices will lead to a block diagonal $\hat{T}$ which is easier to handle later. Now, we will find $p_{12}$ such that $\det(\hat{P})=0$,
\begin{align*}
    \det(\hat{P})=0\rightarrow p_{11}p_{22}-p_{12}^2=0\rightarrow p_{12}=\frac{n}{4}\left(\frac{n-\lambda}{\frac{m(n-\lambda)-np}{2}}\right).
\end{align*}
From the quadratics analysis, we expect the fastest convergence rate to relate to $q$ with condition $n=q$. Therefore, we set $\lambda=q$ and $n=q$. These two conditions lead to
$ p_{12}=0,\quad p_{11}=0,\quad p_{22}=\frac{n}{2p},$
and since $\tfrac{n}{2p}\geq 0$ we get $\hat{P}\succeq 0$. Also, we have
$
    t_{11}=-\frac{q\mu}{2},\quad t_{12}=\frac{q^2}{2p},\quad t_{13}=0, \quad t_{22}=\frac{-q^2}{2p},\quad t_{23}=0,\quad t_{33}=-m.
$
Now, to establish $\hat{T}\preceq 0$, consider 
$$\hat{T}_1=\left[\begin{array}{cc}
   t_{11}  & t_{12} \\
  t_{21}   & t_{22}
\end{array}\right],\quad \hat{T}_2=[t_{33}],$$ 
as the blocks in the block diagonal matrix $\hat{T}$. We know that $t_{33}\leq 0$. Also if $\mathrm{Tr}(.)$ denotes the trace operator, $\mathrm{Tr}(A)$ is equal to the sum of the eigenvalues of the matrix $A$ and the determinant of $A$ is equal to the multiplication of its eigenvalues. Therefore, if the determinant of the first $2\times 2$ block matrix, $\hat{T}_1$, is positive and $t_{11},t_{22}$ are negative, we ensure that $\mathrm{Tr}(\hat{T}_1)=\Gamma_1+\Gamma_2\leq 0$, $\Gamma_1\Gamma_2\geq 0$, and $\Gamma_3\leq 0$ where $\Gamma_i,i=1,2,3$ are the eigenvalues of $\hat{T}$. Hence, $\Gamma_i\leq 0 \text{ for }i=1,2,3$ and this means $\hat{T}\preceq 0$. One can formulate the above arguments as
$t_{11}t_{22}-t_{12}^2\geq 0\rightarrow (\frac{q\mu}{2})(\frac{q^2}{2p})-(\frac{q^2}{2p})^2\geq 0\rightarrow q\leq \mu p.$
Now, using Theorem 6.4 in \citep{doi:10.1137/17M1136845} we get
\begin{align*}
    f(x(t))-f(x^*)\leq e^{-q t}(f(x(0))-f(x^*)+\frac{q}{2p}\|V-x^*\|^2).
\end{align*}
Note that due to the form of (\ref{proposed_ODE_conti}), $x^*=v^*$ (see Figure \ref{fig1}). Therefore, the claim of the theorem is proved.

\noindent(b) The result follows directly by setting $\alpha = \log n,$ $ \beta = qt,$ $\gamma = nt,$ $m=\sqrt{s}e^{\alpha_t}/C,$ and $C=1/p$ in \Cref{Theorem3_1} (b).

\subsubsection{Proof of \Cref{theorem4}}\label{prf:disc_conv}
 To show the claim of the theorem, we will bound the difference $\varepsilon(k+1)-\varepsilon(k)$ such that
$\varepsilon(k+1) \leq (1-q\sqrt{s})^{(k+1)}\varepsilon(0),$
 with $\varepsilon(k)$ as the Lyapunov function
\begin{align} \label{lyap_discrete2}
    \varepsilon(k)= f(x_k)-f(x^*)+\frac{B}{2}\|v_k-x^*\|_2^2-\frac{Bp^2s}{2}\|\nabla f(x_{k})\|^2,
\end{align}
with $B$ as a positive constant to be found. Using (\ref{lyap_discrete2}) we have
\begin{align} \label{calc7}
        \begin{aligned}
            \varepsilon(k+1)-\varepsilon(k)=& f(x_{k+1})-f(x_k)+\underbrace{\frac{B}{2}(\|v_{k+1}-v_k\|^2+2\langle v_{k+1}-v_k,v_{k}-x^* \rangle)}_{\textbf{(I)}} \\
            &-\frac{Bp^2s}{2}\|\nabla f(x_{k+1})\|^2+\frac{Bp^2s}{2}\|\nabla f(x_{k})\|^2.
        \end{aligned}
    \end{align}
    Note that for \textbf{(I)} we have
\begin{align}
    \textbf{(I)}\overset{\ref{SIE_newalg}}{=}& \frac{B}{2}(p^2s\|\nabla f(x_{k+1})\|^2+q^2 s\|v_k-x_{k+1}\|^2+2pqs\langle \nabla f(x_{k+1}),v_k-x_{k+1} \rangle \nonumber\\
    &-2p\sqrt{s}\langle \nabla f(x_{k+1}),v_{k}-x_{k+1} \rangle-2p\sqrt{s}\langle \nabla f(x_{k+1}),x_{k+1}-x^* \rangle\nonumber\\
    &-2q\sqrt{s}\langle v_k-x_{k+1},v_k-x^*\rangle),\nonumber
\end{align}
    where we have added and subtracted $x_{k+1}$ in $\langle\nabla f(x_{k+1}),v_k-x^*\rangle$. Next, using
    $$\langle a-b,a-c\rangle = \frac{1}{2}\|a-b\|^2+\frac{1}{2}\|a-c\|^2-\frac{1}{2}\|b-c\|^2,$$
    we get
\begin{align*} 
    \textbf{(I)}\overset{ }{=}& \frac{B}{2}(p^2s\|\nabla f(x_{k+1})\|^2+q^2 s\|v_k-x_{k+1}\|^2+2pqs\langle \nabla f(x_{k+1}),v_k-x_{k+1} \rangle \\
    &-2p\sqrt{s}\langle \nabla f(x_{k+1}),v_{k}-x_{k+1} \rangle-2p\sqrt{s}\langle \nabla f(x_{k+1}),x_{k+1}-x^* \rangle\\
    & -q\sqrt{s}\|v_k-x_{k+1}\|^2-q\sqrt{s}\|v_k-x^*\|^2+q\sqrt{s}\|x_{k+1}-x^*\|^2).
\end{align*}
Utilizing strong convexity of $f(x)$ we have
$$\langle \nabla f(x_{k+1}),x_{k+1}-x^* \rangle\geq f(x_{k+1})-f(x^*)+\frac{\mu}{2}\|x_{k+1}-x^*\|^2,\quad \quad \quad \textbf{(S.C)}$$
thus, we can upper bound \textbf{(I)} as

\begin{align} \label{calc9}
    \begin{aligned}
        \textbf{(I)}\overset{\textbf{(S.C)}}{\leq}& \frac{B}{2}(p^2s\|\nabla f(x_{k+1})\|^2+q^2 s\|v_k-x_{k+1}\|^2+2pqs\langle \nabla f(x_{k+1}),v_k-x_{k+1} \rangle \\
        &-2p\sqrt{s}\langle \nabla f(x_{k+1}),v_{k}-x_{k+1} \rangle-2p\sqrt{s}(f(x_{k+1})-f(x^*))-\mu p\sqrt{s}\|x_{k+1}-x^*\|^2\\
        & -q\sqrt{s}\|v_k-x_{k+1}\|^2-q\sqrt{s}\|v_k-x^*\|^2+q\sqrt{s}\|x_{k+1}-x^*\|^2)\\
        \overset{\pm f(x_k)}{=}&\frac{B}{2}(p^2s\|\nabla f(x_{k+1})\|^2+q^2 s\|v_k-x_{k+1}\|^2+2pqs\langle \nabla f(x_{k+1}),v_k-x_{k+1} \rangle \\
        &-2p\sqrt{s}\langle \nabla f(x_{k+1}),v_{k}-x_{k+1} \rangle-2p\sqrt{s}(f(x_{k+1})-f(x_k))-2p\sqrt{s}(f(x_k)-f(x^*))\\
        & -\mu p\sqrt{s}\|x_{k+1}-x^*\|^2-q\sqrt{s}\|v_k-x_{k+1}\|^2-q\sqrt{s}\|v_k-x^*\|^2+q\sqrt{s}\|x_{k+1}-x^*\|^2).
    \end{aligned}
\end{align}
Replacing (\ref{calc9}) in (\ref{calc7}) we have
\begin{align} \label{calc10}
    \begin{aligned}
         \varepsilon(k+1)-\varepsilon(k)\overset{\ref{calc9}}{\leq}& f(x_{k+1})-f(x_k)+\frac{Bp^2s}{2}\|\nabla f(x_{k+1})\|^2+\frac{Bq^2 s}{2}\|v_k-x_{k+1}\|^2 \\
        &+Bpqs\langle \nabla f(x_{k+1}),v_k-x_{k+1} \rangle-Bp\sqrt{s}\langle \nabla f(x_{k+1}),v_{k}-x_{k+1} \rangle\\
        &-Bp\sqrt{s}(f(x_{k+1})-f(x_k))-Bp\sqrt{s}(f(x_k)-f(x^*))\\
        & -\frac{B\mu p\sqrt{s}}{2}\|x_{k+1}-x^*\|^2-\frac{Bq\sqrt{s}}{2}\|v_k-x_{k+1}\|^2-\frac{Bq\sqrt{s}}{2}\|v_k-x^*\|^2 \\
        &+\frac{Bq\sqrt{s}}{2}\|x_{k+1}-x^*\|^2-\frac{Bp^2s}{2}\|\nabla f(x_{k+1})\|^2+\frac{Bp^2s}{2}\|\nabla f(x_{k})\|^2.
    \end{aligned}
\end{align}
By setting $q/p\leq \mu$ we get
$\frac{Bq\sqrt{s}}{2}\|x_{k+1}-x^*\|^2\leq \frac{B\mu p\sqrt{s}}{2}\|x_{k+1}-x^*\|^2,$
and thus simplifying (\ref{calc10}) results in
\begin{align} \label{calc11}
    \begin{aligned}
        \varepsilon(k+1)-\varepsilon(k)\overset{\ref{calc10}}{\leq}& Bp\sqrt{s}(q\sqrt{s}-1)\langle \nabla f(x_{k+1}),v_k-x_{k+1} \rangle\\
        &-Bp\sqrt{s}(f(x_k)-f(x^*))-\frac{Bq \sqrt{s}}{2}(1-q\sqrt{s})\|v_k-x_{k+1}\|^2\\
        &-\frac{Bq\sqrt{s}}{2}\|v_k-x^*\|^2+\frac{Bp^2s}{2}\|\nabla f(x_{k})\|^2 +(1-Bp\sqrt{s})(f(x_{k+1})-f(x_k)).
    \end{aligned}
\end{align}
Next, using smoothness of $f(x)$ we get
\begin{align} \label{calc12}
    \begin{aligned}
        f(x_{k+1})-f(x_k)\leq& \langle \nabla f(x_{k+1}), x_{k+1}-x_k\rangle - \frac{1}{2L}\|\nabla f(x_{k+1})-\nabla f(x_k)\|^2 \\
        \overset{\ref{SIE_newalg}}{=}& -m\sqrt{s}\langle \nabla f(x_{k+1}),\nabla f(x_k) \rangle-n\sqrt{s}\langle \nabla f(x_{k+1}),x_{k+1}-v_k\rangle\\
        &-\frac{1}{2L}\|\nabla f(x_{k+1})-\nabla f(x_k)\|^2.
    \end{aligned}\tag{\textbf{S.L}}
\end{align}
Upper-bounding (\ref{calc11}) using (\ref{calc12}) and considering $q\sqrt{s}\leq 1$ leads to
\begin{align} \label{calc13}
\begin{aligned}
    \varepsilon(k+1)-\varepsilon(k)\overset{\ref{calc12}}{\leq}& -m\sqrt{s}(1-Bp\sqrt{s})\langle \nabla f(x_{k+1}),\nabla f(x_k) \rangle+\frac{Bp^2s}{2}\|\nabla f(x_{k})\|^2\\
    & +[Bp\sqrt{s}(q\sqrt{s}-1)+n\sqrt{s}(1-Bp\sqrt{s})]\langle \nabla f(x_{k+1}),v_k-x_{k+1} \rangle\\
    &-Bp\sqrt{s}(f(x_k)-f(x^*))-\frac{Bq\sqrt{s}}{2}\|v_k-x^*\|^2-\frac{(1-Bp\sqrt{s})}{2L}\|\nabla f(x_{k+1})-\nabla f(x_k)\|^2.
\end{aligned}
\end{align}
Next, by setting $B=\frac{n}{p}$, $n=q$, and $q\sqrt{s}<1$ we have

\begin{align} \label{calc14}
    \begin{aligned}
        \varepsilon(k+1)-\varepsilon(k)\overset{\ref{calc13}}{\leq}& -m\sqrt{s}(1-Bp\sqrt{s})\langle \nabla f(x_{k+1}),\nabla f(x_k) \rangle\\
         &-n\sqrt{s}(f(x_k)-f(x^*))-\frac{Bq\sqrt{s}}{2}\|v_k-x^*\|^2+\frac{Bp^2s}{2}\|\nabla f(x_{k})\|^2\\
         =& -m\sqrt{s}(1-Bp\sqrt{s})\langle \nabla f(x_{k+1}),\nabla f(x_k) \rangle-\frac{(1-Bp\sqrt{s})}{2L}\|\nabla f(x_{k+1})-\nabla f(x_k)\|^2\\
         &-n\sqrt{s}(f(x_k)-f(x^*))-\frac{Bq\sqrt{s}}{2}\|v_k-x^*\|^2+\frac{Bp^2s}{2}(q\sqrt{s})\|\nabla f(x_{k})\|^2\\
         &+\frac{Bp^2s}{2}(1-q\sqrt{s})\|\nabla f(x_{k})\|^2-\frac{(1-Bp\sqrt{s})}{2L}\|\nabla f(x_{k+1})-\nabla f(x_k)\|^2.
    \end{aligned}
\end{align}
Now, adding $\frac{Bp^2s(1-q\sqrt{s})}{2}\|\nabla f(x_{k+1})\|^2$ to (\ref{calc14}) results in
\begin{align} \label{calc15}
\begin{aligned}
    \varepsilon(k+1)-\varepsilon(k)\overset{\ref{calc14}}{\leq}&-m\sqrt{s}(1-Bp\sqrt{s})\langle \nabla f(x_{k+1}),\nabla f(x_k) \rangle\\
     &-\frac{(1-Bp\sqrt{s})}{2L}\|\nabla f(x_{k+1})-\nabla f(x_k)\|^2-n\sqrt{s}(f(x_k)-f(x^*))\\
     &-\frac{Bq\sqrt{s}}{2}\|v_k-x^*\|^2+\frac{Bp^2s}{2}(1-q\sqrt{s})\|\nabla f(x_{k+1})\|^2\\
     &+\frac{Bp^2s}{2}(1-q\sqrt{s})\|\nabla f(x_{k})\|^2+\frac{Bp^2s}{2}(q\sqrt{s})\|\nabla f(x_{k})\|^2.
\end{aligned}
\end{align}
Next, setting $nps\leq m\sqrt{s}$ and noting that $n=q$, we get 
\begin{align} \label{calc16}
    \begin{aligned}
        &-m\sqrt{s}(1-Bp\sqrt{s})\langle \nabla f(x_{k+1}),\nabla f(x_k) \rangle+\tfrac{Bp^2s(1-q\sqrt{s})}{2}\|\nabla f(x_{k+1})\|^2
        +\tfrac{Bp^2s(1-q\sqrt{s})}{2}\|\nabla f(x_{k})\|^2\\
        = & (1-q\sqrt{s})\left[-m\sqrt{s}\langle \nabla f(x_{k+1}),\nabla f(x_k) \rangle +\frac{nps}{2}\|\nabla f(x_{k+1})\|^2+\frac{nps}{2}\|\nabla f(x_{k})\|^2\right]\\
        \leq & \frac{(1-q\sqrt{s})m\sqrt{s}}{2}\|\nabla f(x_{k+1})-\nabla f(x_k)\|^2.
    \end{aligned}
\end{align}
Using (\ref{calc16}) in (\ref{calc15}) we get
\begin{align*}
\varepsilon(k+1)-\varepsilon(k)\overset{\ref{calc16}}{\leq}& (1-q\sqrt{s})(\frac{m\sqrt{s}}{2}-\frac{1}{2L})\|\nabla f(x_{k+1})-\nabla f(x_k)\|^2\\
& -n\sqrt{s}(f(x_k)-f(x^*))-\frac{Bq\sqrt{s}}{2}\|v_k-x^*\|^2+\frac{Bp^2s}{2}(q\sqrt{s})\|\nabla f(x_{k})\|^2.
\end{align*}
Setting $m\sqrt{s}\leq \frac{1}{L}$ leads to
\begin{align*}
\varepsilon(k+1)-\varepsilon(k)&\leq -n\sqrt{s}(f(x_k)-f(x^*))-\frac{Bq\sqrt{s}}{2}\|v_k-x^*\|^2 +\frac{Bp^2s}{2}(q\sqrt{s})\|\nabla f(x_{k})\|^2
\end{align*}
which can be expressed in a more favourable way 
\begin{align*} 
\varepsilon(k+1)-\varepsilon(k)\overset{ }{\leq}& -(q\sqrt{s})\left[f(x_k)-f(x^*)+ \frac{B}{2}\|v_k-x^*\|^2-\frac{Bp^2s}{2}\|\nabla f(x_{k})\|^2\right]= -(q\sqrt{s}) \varepsilon(k), 
\end{align*}
which gives $\varepsilon(k+1)\leq (1-q\sqrt{s})\varepsilon(k),$ and thus
\begin{align} \label{calc18}
    \varepsilon(k+1)\leq (1-q\sqrt{s})^{k+1}\varepsilon(0). 
\end{align}
Using the form of $\varepsilon(k)$ in (\ref{lyap_discrete2}) and the inequality
$f(x_k)-f(x^*)\geq \frac{1}{2L}\|\nabla f(x_k)\|^2,$
which is true for any $L$-smooth function with $x^*$ such that $\nabla f(x^*)=0$, we get
\begin{align} \label{proof_lem2_eqn1}
    \begin{aligned}
        \varepsilon(k)=&(f(x_k)-f(x^*))+\frac{B}{2}\|v_k-x^*\|_2^2-\frac{Bp^2s}{2}\|\nabla f(x_{k})\|^2\\
        &\geq (1-Bp^2sL)(f(x_k)-f(x^*))+\frac{B}{2}\|v_k-x^*\|^2.
    \end{aligned}
\end{align}
Note that
$1-Bp^2sL=1-npsL,$
and under the conditions in Theorem \ref{theorem4} we have $npsL\leq 1$ and thus,
$1-Bp^2sL=1-npsL\geq 0.$
Hence, $(1-Bp^2sL)(f(x_k)-f(x^*))\geq 0$ and (\ref{proof_lem2_eqn1}) leads to
\begin{align} \label{proof_lem2_eqn2}
     \varepsilon(k)&\geq \frac{B}{2}\|v_k-x^*\|_2^2.
\end{align}
With $v_0=x_0-\tfrac{m}{n}\nabla f(x_0)$ in (\ref{calc18}) we have
\begin{align} \label{proof_lem2_eqn3}
    \begin{aligned}
        \varepsilon(k)\leq & (1-q\sqrt{s})^k\varepsilon(0)\\
        =&(1-q\sqrt{s})^k\left[f(x_0)-f(x^*)+\frac{B}{2}\|v_0-x^*\|_2^2-\frac{Bp^2s}{2}\|\nabla f(x_{0})\|^2\right]\\
        =&(1-q\sqrt{s})^k\left[f(x_0)-f(x^*)+\frac{B}{2}\|x_0-x^*\|^2+(\frac{Bm^2}{2n^2}-\frac{Bp^2s}{2})\|\nabla f(x_0)\|^2\right.\\
        &\left.-\frac{Bm}{n}\langle \nabla f(x_0),x_0-x^*\rangle \right]\\
        \overset{\textbf{S.C}}{\leq}& (1-q\sqrt{s})^k\left[f(x_0)-f(x^*)+\frac{B}{2}\|x_0-x^*\|^2+(\frac{Bm^2}{2n^2}-\frac{Bp^2s}{2})\|\nabla f(x_0)\|^2\right.\\
        &\left.-\frac{Bm}{n}(f(x_0)-f(x^*))-\frac{B\mu m}{2n}\|x_0-x^*\|^2\right]\\
        =&(1-q\sqrt{s})^k\left[(1-B(\frac{m}{n}))(f(x_0)-f(x^*))+\frac{B}{2}(1-\frac{\mu m}{2n})\|x_0-x^*\|^2\right.\\
        &\left. +(\frac{Bm^2}{2n^2}-\frac{Bp^2s}{2})\|\nabla f(x_0)\|^2\right]\\
        \overset{\textbf{S.L}}{\leq}& (1-q\sqrt{s})^k\left[\frac{L}{2}(1-B(\frac{m}{n}))+\frac{B}{2}(1-\frac{\mu m}{2n})+(\frac{Bm^2}{2n^2}-\frac{Bp^2s}{2})L^2)\right]\|x_0-x^*\|^2.
    \end{aligned}
\end{align}
Therefore, from (\ref{proof_lem2_eqn2}) and (\ref{proof_lem2_eqn3}) we get
\begin{align*}
    \frac{B}{2}\|v_k-x^*\|_2^2&\leq (1-q\sqrt{s})^k\left[\frac{L}{2}(1-B(\frac{m}{n})) + \frac{B}{2}(1-\frac{\mu m}{2n})+(\frac{Bm^2}{2n^2}-\frac{Bp^2s}{2})L^2)\right]\|x_0-x^*\|^2.
\end{align*}
The positivity of the coefficient 
$$\left[\frac{L}{2}(1-B(\frac{m}{n})) +\frac{B}{2}(1-\frac{\mu m}{2n})+(\frac{Bm^2}{2n^2}-\frac{Bp^2s}{2})L^2)\right]\|x_0-x^*\|^2$$
is guaranteed under the conditions of Theorem \ref{theorem4}. Multiplying both sides by $2/B$ gives 
$\|v_k-x^*\|_2^2 \leq (1-q\sqrt{s})^kM'\|x_0-x^*\|^2$
for 
\begin{align*}
    M'=\left[\frac{L}{B}(1-B(\frac{m}{n})) +(1-\frac{\mu m}{2n})+(\frac{m^2}{n^2}-p^2s)L^2)\right],
\end{align*}
which is nonnegative due to the conditions of the theorem. Now, we know that $\lim_{k\rightarrow \infty}v_k=x^*$. By representing algorithm (\ref{SIE_newalg}) in one-line format of sequence $v_k$ we get
\begin{align} \label{eqn2_thm6}
 v_{k+1}-v_k=-p\sqrt{s}\nabla f(x_{k+1})+\frac{1-q\sqrt{s}}{1+n\sqrt{s}}(v_k-v_{k-1})+\frac{q\sqrt{s}}{1+n\sqrt{s}}(\frac{p}{q}-m\sqrt{s})\nabla f(x_k).
\end{align}
Analyzing (\ref{eqn2_thm6}) in limit and Taking $\iota=\tfrac{1}{1+n\sqrt{s}}(1-\tfrac{qm\sqrt{s}}{p})$ we have
\begin{align*}
    \lim_{k\rightarrow \infty}  (\nabla f(x_{k+1})-\iota\nabla f(x_{k}))  =0 .
\end{align*}
for $\iota< 1$. Thus, 
\begin{align}
   \forall{\epsilon}\quad \exists k_0\quad \text{s.t. }& \|\nabla f(x_{k+1})-\iota \nabla f(x_k)\|\leq \epsilon\qquad k\geq k_0, \nonumber
\end{align}
From the above argument we get
$$\|\nabla f(x_{k+1})\|-\iota\|\nabla f(x_k)\|\leq\epsilon\rightarrow \|\nabla f(x_{k+1})\|\leq\iota\|\nabla f(x_k)\|+\epsilon.$$
Unrolling the above inequality results in
$$\|\nabla f(x_{k+1})\|\leq \iota^{k-k_0+1}\|\nabla f(x_{k_0})\|+\epsilon(1+\iota+\iota^2+\ldots+\iota^{k-k_0}),$$
$$\rightarrow \|\nabla f(x_{k+1})\|\leq \iota^{k-k_0+1}\|\nabla f(x_{k_0})\|+\frac{\epsilon}{1-\iota}.$$
Taking $\epsilon = \tfrac{(1-\iota)\delta}{2}$ and noting that
$$\forall \delta \quad\exists k_0\leq k:\iota^{k-k_0+1}\|\nabla f(x_{k_0})\|\leq \delta/2,$$ 
we get
$\|\nabla f(x_{k+1})\|\leq \delta,$
and therefore, 
$\lim_{k\rightarrow  \infty}\|\nabla f(x_{k+1})\|=0.$
Thus, the limit of $\nabla f(x_k)$ exists and since $f$ is strongly convex $x_k\rightarrow x^*$ as $k\rightarrow \infty$. 
We therefore proved that sequence $x_k$ converges to $x^*$.

\subsubsection{Proof of \Cref{theorem_GDnorm}}\label{thm5_proof}
    Take the Lyapunov function 
\begin{align} \label{Lyapunov_function1}
    \varepsilon(k) = \frac{s(k+2)k}{4}(f(x_k)-f(x^*))+\frac{1}{2}\|x_{k+1}-x^*+\frac{k}{2}(x_{k+1}-x_k)+\frac{ks}{2}\nabla f(x_k)\|^2.
\end{align}
The choice of Lyapunov function is the same as \citep{Shi2021UnderstandingTA}. Note that the second term is equivalent to $\frac{1}{2}\|v_k-x^*\|^2$ through the first line of the update rule (\ref{new_algorithm}). Next, we will show that 
\begin{align} \label{Lyap_1}
    \varepsilon(k+1)-\varepsilon(k)\leq -\frac{s^2k(k+2)}{8}\|\nabla f(x_k)\|^2.
\end{align}
Using (\ref{Lyapunov_function1}) we have
\begin{align} \label{Lyap_2}
    \begin{aligned}
        \varepsilon(k+1)-\varepsilon(k)=&\frac{s(k+3)(k+1)}{4}(f(x_{k+1})-f(x^*))+\frac{1}{2}\|v_{k+1}-x^*\|^2\\
        & - \frac{s(k+2)(k)}{4}(f(x_{k})-f(x^*))-\frac{1}{2}\|v_{k}-x^*\|^2\\
        =&\frac{s(k+2)k}{4}(f(x_{k+1})-f(x_k))+\frac{s(2k+3)}{4}(f(x_{k+1})-f(x^*))\\
        &+\frac{1}{2}(2\langle v_{k+1}-v_k,v_k-x^* \rangle+\|v_{k+1}-v_k\|^2)\\
        =& \frac{s(k+2)k}{4}(f(x_{k+1})-f(x_k))+\frac{s(2k+3)}{4}(f(x_{k+1})-f(x^*))\\
        &+\frac{1}{2}(2\langle -s(\frac{k+2}{2})\nabla f(x_{k+1}),x_{k+1}-x^*+\frac{k}{2}(x_{k+1}-x_k)+\frac{ks}{2}\nabla f(x_k) \rangle \\
        &+\|s(\frac{k+2}{2})\nabla f(x_{k+1})\|^2)\\
        =& \frac{s(k+2)k}{4}(f(x_{k+1})-f(x_k))+\frac{s(2k+3)}{4}(f(x_{k+1})-f(x^*))\\
        &-s(\frac{k+2}{2})\langle \nabla f(x_{k+1}),x_{k+1}-x^*\rangle-s\frac{k(k+2)}{4}\langle \nabla f(x_{k+1}),x_{k+1}-x_k\rangle\\
        & -s^2(\frac{k(k+2)}{4})\langle \nabla f(x_{k+1}),\nabla f(x_k)\rangle+ \frac{(s(k+2))^2}{8}\|\nabla f(x_{k+1})\|^2
    \end{aligned}
\end{align}
Now, from convexity and smoothness of the function $f$ we have
\begin{align} \label{smooth_convex}
    f(x_{k+1})-f(x_k)\leq \langle \nabla f(x_{k+1}),x_{k+1}-x_k\rangle -\frac{1}{2L}\|\nabla f(x_{k+1})-\nabla f(x_k)\|^2.
\end{align}
Applying (\ref{smooth_convex}) in (\ref{Lyap_2}) we get
\begin{align*} 
     \varepsilon(k+1)-\varepsilon(k)\leq& \frac{s(k+2)k}{4}\left[ \langle \nabla f(x_{k+1}),x_{k+1}-x_k\rangle-\frac{1}{2L}\|\nabla f(x_{k+1})-\nabla f(x_k)\|^2 \right]\\
     & +\frac{s(2k+3)}{4} \left[ \langle \nabla f(x_{k+1}),x_{k+1}-x^*\rangle-\frac{1}{2L}\|\nabla f(x_{k+1})\|^2 \right]\\
     & -s(\frac{k(k+2)}{4})\langle \nabla f(x_{k+1}),x_{k+1}-x_k \rangle-s(\frac{k+2}{2})\langle \nabla f(x_{k+1}),x_{k+1}-x^* \rangle\\
     & -s^2(\frac{k(k+2)}{4})\langle \nabla f(x_{k+1}),\nabla f(x_k)\rangle+ \frac{(s(k+2))^2}{8}\|\nabla f(x_{k+1})\|^2.
\end{align*}
Now, using $$-\langle \nabla f(x_{k+1}),x_{k+1}-x^* \rangle\leq -\frac{1}{2L}\|\nabla f(x_{k+1})\|^2,$$ 
\vspace{-2cm}
we have
\begin{equation*}
\begin{aligned}
    \varepsilon(k+1)-\varepsilon(k)\leq& -\frac{s(k+2)k}{8L}\|\nabla f(x_{k+1})-\nabla f(x_k)\|^2-\frac{s(2k+4)}{8L}\|\nabla f(x_{k+1})\|^2 \\
     & -2s^2(\frac{k(k+2)}{8})\langle \nabla f(x_{k+1}),\nabla f(x_k)\rangle+\frac{s^2k(k+2)}{8}\|\nabla f(x_{k+1})\|^2\\
     &+\frac{s^2(k+2)}{4}\|\nabla f(x_{k+1})\|^2\\
     =&  -\frac{s(k+2)k}{8}(\frac{1}{L}-s)\|\nabla f(x_{k+1})\|^2-\frac{s(k+2)}{4}(\frac{1}{L}-s)\|\nabla f(x_{k+1})\|^2\\
     &+2s(\frac{k(k+2)}{8})(\frac{1}{L}-s)\langle \nabla f(x_{k+1}),\nabla f(x_k)\rangle\\
     &-\frac{s(k+2)k}{8}(\frac{1}{L}-s)\|\nabla f(x_{k})\|^2-s^2(\frac{k(k+2)}{8})\|\nabla f(x_{k})\|^2\\
     =& -\frac{s(k+2)k}{8}(\frac{1}{L}-s)\|\nabla f(x_{k+1})-\nabla f(x_k)\|^2\\
     &-\frac{s(k+2)}{4}(\frac{1}{L}-s)\|\nabla f(x_{k+1})\|^2\\
     &-s^2(\frac{k(k+2)}{8})\|\nabla f(x_{k})\|^2\leq -s^2(\frac{k(k+2)}{8})\|\nabla f(x_{k})\|^2,
\end{aligned}     
\end{equation*}
where the last inequality holds as long as $s\leq 1/L$.\par
\noindent With (\ref{Lyap_1}) at hand, we can make sum both sides from $i=0$ till $i=k-1$ and get
\begin{align*} 
        \varepsilon(k)-\varepsilon(0)&\leq -\frac{s^2}{8}\sum_{i=0}^{k} i(i+2)\|\nabla f(x_i)\|^2\\
        & \leq -\frac{s^2}{8}\min_{0\leq i\leq k}\|\nabla f(x_i)\|^2\sum_{i=0}^{k} i(i+2)\\
        &=-\frac{s^2}{8}\min_{0\leq i\leq k}\|\nabla f(x_i)\|^2\sum_{i=1}^{k} i(i+2)\\
        &=-\frac{s^2}{8}\min_{0\leq i\leq k}\|\nabla f(x_i)\|^2\left[ \frac{k(k+1)(2k+1)}{6}+k(k+1) \right]\\
        &\leq -\frac{s^2}{8}\min_{0\leq i\leq k}\|\nabla f(x_i)\|^2\left[ \frac{k(k+1)(2k+1)}{6}\right]\\
        &\leq -\frac{k^3s^2}{24}\min_{0\leq i\leq k}\|\nabla f(x_i)\|^2.
\end{align*}
This completes the proof.

\subsubsection{Proof of \Cref{prop1}}\label{proof_prop1}
    From the update rule (\ref{rate_match_2}) we have
    \begin{align} \label{prop1_eqn1}
        z_k = x_{k+1}+\frac{k}{2}(x_{k+1}-y_k)=x_{k+1}+\frac{k}{2}(x_{k+1}-x_k+s \nabla f(x_k)).
    \end{align}
    Replacing (\ref{prop1_eqn1}) in the update rule of $z_k$ in (\ref{rate_match_2}), we get
    \begin{align*} 
    \begin{aligned}
        &x_{k+1}-x_k+\frac{k}{2}(x_{k+1}-x_k+s\nabla f(x_k))-\frac{k-1}{2}(x_k-x_{k-1}+s \nabla f(x_{k-1}))\\
        =&-\frac{sk}{2}\nabla f(y_k).
    \end{aligned}
    \end{align*}
    By rearranging we have
      \begin{align*} 
      \begin{aligned}
          & x_{k+1}-x_k+\frac{1}{2}(x_{k}-x_{k-1})+\frac{s}{2}\nabla f(x_{k-1})+\frac{k}{2}(x_{k+1}+x_{k-1}-2x_k)\\
        + & \frac{ks}{2}(\nabla f(x_{k})- \nabla f(x_{k-1}))=-\frac{s k}{2}\nabla f(y_k),\\
        \rightarrow & \frac{2}{k\sqrt{s}}(\frac{x_{k+1}-x_k}{\sqrt{s}})+\frac{1}{k\sqrt{s}}(\frac{x_k-x_{k-1}}{\sqrt{s}}) +\frac{1}{k}\nabla f(x_{k-1}) \\
         +&\frac{x_{k+1}-2x_k+x_{k-1}}{s}+\nabla f(x_k)-\nabla f(x_{k-1})=-\nabla f(y_k).
      \end{aligned}
    \end{align*}
    Then, we use the following approximations with $t_k = (k+1/2)\sqrt{s}$
    \begin{align*}
        \begin{aligned} 
           &\frac{2}{k\sqrt{s}}(\frac{x_{k+1}-x_k}{\sqrt{s}})+\frac{1}{k\sqrt{s}}(\frac{x_k-x_{k-1}}{\sqrt{s}}) \approx \frac{1}{t_k - \sqrt{s}/2}(3 \dot X(t_k) + \frac{\sqrt{s}\Ddot{X}(t_k)}{2})\\
            &\frac{x_{k+1}-2x_k+x_{k-1}}{s}\approx \Ddot{X}(t_k)\\
            &\nabla f(x_k)-\nabla f(x_{k-1})\approx \sqrt{s}\nabla^2 f(X(t_k))\dot X(t_k)\\
            &\nabla f(y_k)=\nabla f(x_k-s\nabla f(x_k))\approx \nabla f(x_k)-s\nabla^2 f(x_k)\nabla f(x_k)\approx \nabla f(X(t_k))\\
            &X(t)\approx X(t_k)\qquad \dot X(t)\approx\dot X(t_k)\qquad \Ddot X(t)\approx\Ddot X(t_k) \qquad Y(t_k)\approx Y(t)
    \end{aligned}
    \end{align*}
and we get
\begin{align*}
    \frac{1}{t_k - \sqrt{s}/2}(3 \dot X(t_k) + \frac{\sqrt{s}\Ddot{X}(t_k)}{2}) + \frac{\sqrt{s}}{t_k-\sqrt{s}/2}\nabla f(X(t_k)) + \Ddot{X}(t_k) + \sqrt{s}\nabla^2 f(X(t_k))\dot X(t_k) + \nabla f(X(t_k)) = 0.
\end{align*}
By rearranging we have
\begin{align*}
    \frac{t_k}{t_k - \sqrt{s}/2}\Ddot{X}(t_k) + \frac{3}{t_k - \sqrt{s}/2}( \dot X(t_k) ) + \frac{t_k + \sqrt{s}/2}{t_k-\sqrt{s}/2}\nabla f(X(t_k)) + \sqrt{s}\nabla^2 f(X(t_k))\dot X(t_k) = 0,
\end{align*}
which gives
\begin{align} \label{prop1_eqn5}
    \Ddot{X}(t)+(\frac{3}{t}+\sqrt{s}\nabla^2 f(X(t)))\dot X(t)+(1+\frac{\sqrt{s}}{2t})\nabla f(X(t))=0.
\end{align}
The HR-ODE \Cref{prop1_eqn5} can be modified using $(k+1)/(k+1/2)\approx (k+1)/(k)$ as
\begin{align*} 
    \Ddot{X}(t)+(\frac{3}{t}+\sqrt{s}\nabla^2 f(X(t)))\dot X(t)+(1+\frac{\sqrt{s}}{t})\nabla f(X(t))=0.
\end{align*}
Similar modifications were previously made as in \citep{shi2019acceleration}[Section 4] to demonstrate acceleration of the discretization of the HR-ODE. 
This concludes the proof.

\subsubsection{Proof of Corollary \ref{corl_TM}}\label{proof_TM}
By comparing the one line representation (\ref{SC_ODE_represented}) with (\ref{HR_TM}) we can easily check the validity of parameters stated in the corollary. All that remains is to show that these parameters satisfy the conditions of Theorem \ref{theorem2}. Due to positiveness of the parameters $(\alpha,\beta,\gamma,\delta)$ in (\ref{TM_method}) and $n\leq 2q$ for $\xi\geq 2/3$ the first condition in Theorem \ref{theorem2} is satisfied. Also, for $m\leq 2p$ to hold we need 
\begin{align*}
   \gamma \sqrt{\alpha}(1+\sqrt{M\alpha}) & \leq 2\frac{(1-\xi\gamma\sqrt{M\alpha})(1+\sqrt{M\alpha})}{(2-\xi)\sqrt{M}}\nonumber\\
   \Rightarrow  \gamma\sqrt{\alpha}& \leq 2\frac{(1-\xi\gamma\sqrt{M\alpha})}{(2-\xi)\sqrt{M}}.
\end{align*}
Now, note that $0 \leq \gamma\leq 1/2$. For $\gamma=0$ the result trivially holds. For $\gamma=1/2$ we get 
\begin{align*}
    \sqrt{M\alpha}\leq \frac{2}{1+\frac{\xi}{2}}
\end{align*}
which holds for $\xi\leq 2$ since $\sqrt{M\alpha}\leq 1$.
The last condition is $q/p\leq \mu$. Using the values for $q$ and $p$ we have
\begin{align} \label{eqn1_prf_cor1}
    \begin{aligned}
        \frac{q}{p}=\frac{\xi(2-\xi)M}{(1-\xi\gamma\sqrt{M\alpha})(1+\sqrt{M \alpha})}&\leq \mu\\
        \Rightarrow \frac{(1-\xi\gamma\sqrt{M\alpha})(1+\sqrt{M \alpha})}{\xi(2-\xi)} - \frac{M}{\mu}&\geq 0
    \end{aligned}
\end{align}
Replacing TM method parameters and defining function $G(\xi)$ we get
\begin{align}
    G(\xi,\kappa)=\frac{\left(1-\xi(\frac{\rho^2}{(1+\rho)(2-\rho)})\sqrt{(1+\rho)\tfrac{M}{L}}\right)\left(1+\sqrt{(1+\rho)\tfrac{M}{L}}\right)}{\xi(2-\xi)} - \frac{M}{\mu}.\nonumber
\end{align}
Now, using $M=\left(\tfrac{1-\beta}{\sqrt{\alpha}(1+\beta)}\right)^2$ and equation (8) in \citep{DBLP:conf/cdc/SunGK20} we have 
\begin{align} \label{eqn2_prf_cor1}
    \begin{aligned}
        \frac{M}{\mu}&=\frac{9\kappa^3\sqrt{\kappa}-6\kappa^3+\kappa^2\sqrt{\kappa}}{8\kappa^3\sqrt{\kappa}-12\kappa^3+14\kappa^2\sqrt{\kappa}-9\kappa^2+4\kappa\sqrt{\kappa}-\kappa},\\
        \frac{M}{L}&=\frac{9\kappa^2\sqrt{\kappa}-6\kappa^2+\kappa\sqrt{\kappa}}{8\kappa^3\sqrt{\kappa}-12\kappa^3+14\kappa^2\sqrt{\kappa}-9\kappa^2+4\kappa\sqrt{\kappa}-\kappa}.
    \end{aligned}
\end{align}
With $\rho=1-\tfrac{1}{\sqrt{\kappa}}$, all the terms in (\ref{eqn2_prf_cor1}) depend on $\kappa$ and we need $G(\xi,\kappa)\geq 0$ for (\ref{eqn1_prf_cor1}) to hold. By analyzing $\lim_{\kappa\rightarrow \infty}G(\xi,\kappa)$ we see that for $\xi=\{\tfrac{2}{3},\tfrac{4}{3}\}$ we get $\lim_{\kappa\rightarrow \infty}G(\xi,\kappa)=0$. However, any solution between $\xi=\tfrac{2}{3}$ and $\xi=\tfrac{4}{3}$ is not acceptable since it leads to negativity of $G$. Also, for $\kappa =1$ we have $$G(\xi,1)=\frac{3}{\xi(2-\xi)}-3$$ 
which remains non-negative for $0<\xi< 2$. 
Therefore, invoking Theorem \ref{theorem2} concludes the proof and the best possible result is achieved for $\xi =2/3$.
\begin{figure}
    \centering
    \includegraphics[width = \textwidth]{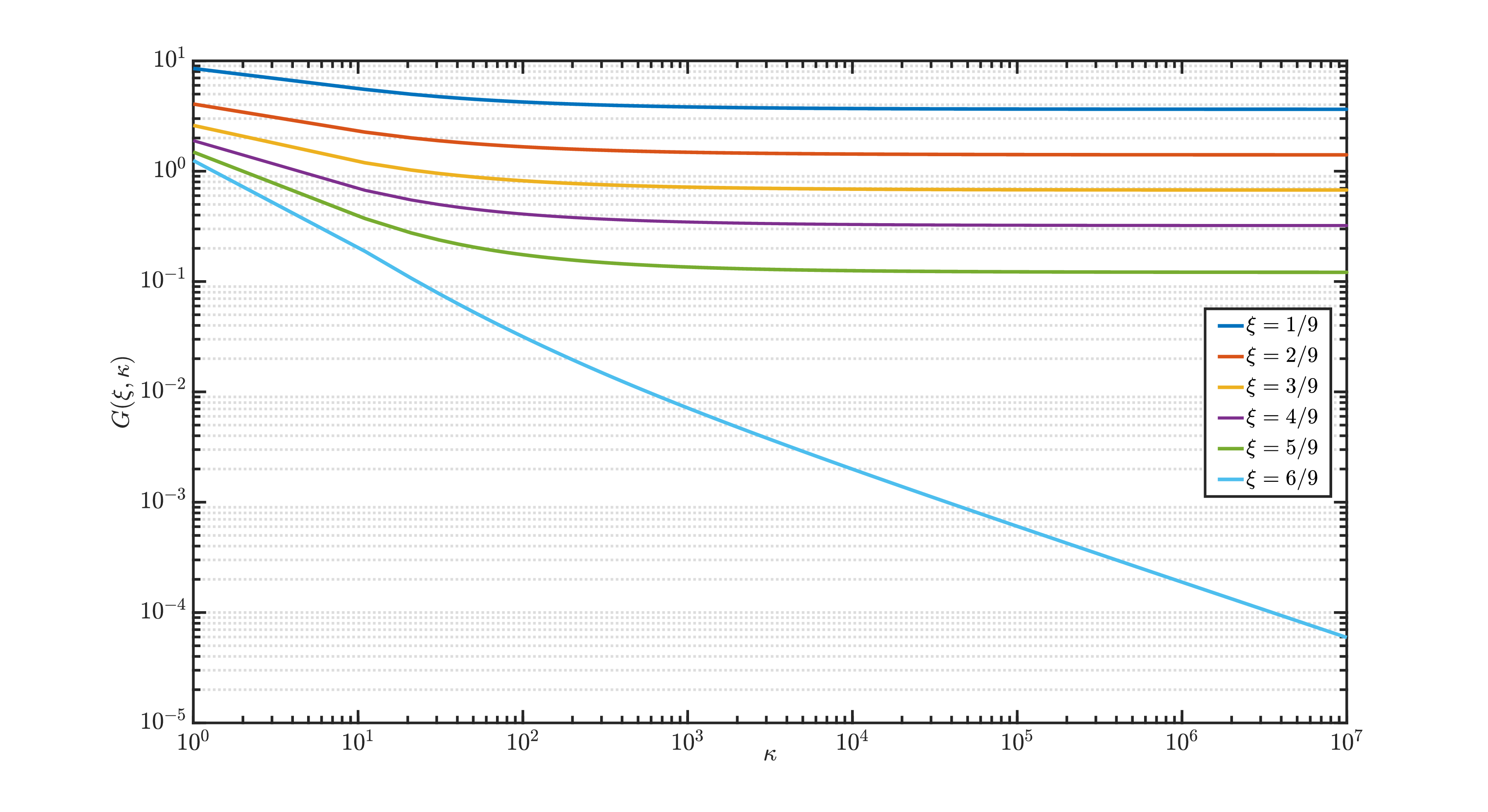}
    \caption{Positivity of $G(\xi,\kappa)$ for various values of $\kappa$}
    \label{fig:enter-label}
\end{figure}

\subsubsection{Proof of Corollary \ref{coro_QHM}}\label{proof_QHM}
The proof is based on the QHM representation in \citep{zhang2021revisiting}. It is possible to rewrite (\ref{QHM}) in one-line format
$$x_{k+1}-x_k=b(x_k-x_{k-1})-s\nabla f(x_k)+sb(1-a)\nabla f(x_{k-1}).$$
On the other hand, (\ref{SIE_newalg}) has the one-line format
$$x_{k+1}-x_k=\frac{1-q\sqrt{s}}{1+n\sqrt{s}}(x_k-x_{k-1})-(\frac{m\sqrt{s}+nps}{1+n\sqrt{s}})\nabla f(x_k)+\frac{m\sqrt{s}(1-q\sqrt{s})}{1+n\sqrt{s}}\nabla f(x_{k-1}).$$
Taking $b=\frac{1-q\sqrt{s}}{1+n\sqrt{s}}$, $m=(1-a)\sqrt{s}$, $n=q$, $p=\frac{a}{n}+\sqrt{s}$ result in
$$p>m,nps=s(a+n\sqrt{s})\leq (1-a)s=m\sqrt{s}\rightarrow a\leq \frac{1-q\sqrt{s}}{2},$$
Also, since $n\sqrt{s}\leq 1/2$ we have $(1-n\sqrt{s})/2\geq 1/4$ and therefore $a\leq 1/4$. For $q/p\leq \mu$ to hold one needs 
$$\frac{q}{\frac{a}{q}+\sqrt{s}}\leq \frac{q^2}{a}\leq \mu\rightarrow q\leq\sqrt{a\mu }\rightarrow q=\sqrt{a\mu},$$
and for $m\sqrt{s}=(1-a)s\leq 1/L$ we need
$$s\leq \frac{1}{L(1-a)}\overset{\frac{1}{1-a}\leq \frac{4}{3}}{\Longrightarrow} s\leq \frac{4}{3L}.$$
Since all the conditions of Theorem \ref{theorem4} are satisfied, for (\ref{QHM}) we have
$$f(x_{k})-f(x^*)\leq C(1-\sqrt{a \mu s})^k.$$

\subsection{Comparison Between GM-ODE and GM2-ODE}\label{sec:gmvsgm2}

A recent analysis \citep{zhang2021revisiting} showed that a suitable choice of the HR-ODE and discretization method can yield an accelerated algorithm, challenging the conventional view that acceleration is inherently tied to the integrator’s properties. 
This connection is investigated through a general ODE, defined as follows: 
\begin{align} 
    \left\{
    \begin{array}{rcl}
    \dot{U}_t &= & -m'\nabla f(U_t)-n'W_t, \\
    \dot{W}_t  &=& \nabla f(U_t)-q'W_t,
    \end{array}
    \right.\tag{GM-ODE}
\end{align}
where $m',n',q'\geq 0$, $U_t=U(t)$, and $W_t=W(t)$; which can be equivalently written as a single line equation: 
\begin{align} \label{zhang_general_ode} 
    \ddot U_t +(q'+m'\nabla^2 f(U_t))\dot{U}_t+(n'+m'q')\nabla f(U_t)=0. 
\end{align}
Of particular interest, \eqref{GM-ODE} was utilized to approximate \eqref{NAG_ODE} in continuous-time, and the NAG algorithm in discrete-time using the SIE integrator.
However, their approximation is inexact due to a coefficient mismatch between the recovered ODE parameters and those of the original \eqref{NAG_ODE}. 
Specifically, the following parameters are selected in their continuous-time analysis: 
\begin{align} \label{eqn:parameters-zhang-cont}
    m'=\sqrt{s}, ~~~ q'=2\sqrt{\mu}, ~~~ \text{and} ~~~n'=1. 
\end{align} 
Substituting these values into \eqref{zhang_general_ode} gives
\begin{align*} 
    \ddot U_t +(2\sqrt{\mu}+\sqrt{s}\nabla^2 f(U_t))\dot{U}_t+(1+2\sqrt{\mu s})\nabla f(U_t)=0,
\end{align*}
which differs from \eqref{NAG_ODE} in the coefficient of $\nabla f(U_t)$. 
Furthermore, a different set of parameters, given by
\begin{align} \label{eqn:parameters-zhang-disc}
    m'=\sqrt{s}, ~~~ q'=2\sqrt{\mu}, ~~~ \text{and} ~~~ n'=1-2\sqrt{\mu s}
\end{align}
is used with the SIE discretization to recover the NAG algorithm, creating a discrepancy between the continuous- and discrete-time models. 

\noindent\eqref{GM2-ODE} not only precisely recovers \eqref{NAG_ODE} and the NAG algorithm with a consistent set of parameters (see \Cref{table_compare1}), but also offers improved convergence guarantees in both continuous and discrete-time. 
To provide context, we first review the convergence results presented in \citep{zhang2021revisiting}. 

\begin{theorem}[{[\citenum{zhang2021revisiting}, Theorem 1]}] \label{theorem1_zhang}

    Let $f$ be a $\mu$-strongly convex and $L$-smooth function. 
    Assume that $m',n',q'\geq 0$. 
    Then, the point ${(x^*,0)\in\mathbb{R}^{2d}}$ is globally asymptotically stable for (\ref{GM-ODE}) as 
    \begin{align*}
        \varepsilon'(U_t,W_t)\leq e^{-\gamma_1t}\varepsilon'(U_0,W_0),
    \end{align*}
    where the Lyapunov function is defined as:
    \begin{align*}
        \varepsilon'(U_t,W_t)=(q'm'+n')(f(U_t)-f(x^*))+\tfrac{1}{4}\|q'(U_t-x^*)-n'W_t\|^2+\tfrac{n'(q'm'+n')}{4}\|W_t\|^2, 
    \end{align*}
    and $\gamma_1:=\min\left(\frac{\mu(n'+q'm')}{2q'},\frac{q'}{2}\right)$.
\end{theorem}

\noindent Substituting the parameters from \eqref{eqn:parameters-zhang-cont} into \Cref{theorem1_zhang} results in the convergence rate 
\begin{align*}
    f(X_t)-f(x^*) \leq \mathcal{O}(e^{-\gamma_1 t}), ~~~\text{with}~~~ \gamma_1 = \min \left( \tfrac{1}{4}\sqrt{\mu}(1 + 2\sqrt{\mu s}), \mu \right).
\end{align*}
Note that $\gamma_1 \leq \tfrac{3}{4}\sqrt{\mu}$ due to the step-size condition $s \leq 1/L$ and the fact that $\mu/L \leq 1$. 
However, this rate is suboptimal: 
while the HR-ODE reduces to Polyak's ODE as $s \to 0$; it does not recover the $\mathcal{O}(e^{-\sqrt{\mu}t})$ rate of Polyak’s ODE \citep{wilson2021lyapunov} in this limit. 

\noindent In contrast, \eqref{proposed_ODE_conti} recovers (\ref{NAG_ODE}) by using the parameter setting
\begin{align} \label{eqn:param-gm2ode-nag}
    n=q=\sqrt{\mu}, ~~~ p=1/\sqrt{\mu} ~~~ \text{and} ~~~ m=\sqrt{s}. 
\end{align}
Using these parameters in \Cref{theorem2} results in a convergence rate of $\mathcal{O}(e^{-\sqrt{\mu} t})$, matching the rate of Polyak’s ODE. 

\begin{remark}
    It is worth noting that the parameter set ($n'=1-\sqrt{\mu s}$, ~ $m'=\sqrt{s}$, and $q'=2\sqrt{\mu}$) in \eqref{GM-ODE} would recover \eqref{NAG_ODE}. 
    However, in this setting,  we obtain  $\gamma_1 = \min \left( \frac{1}{4}\sqrt{\mu}(1+\sqrt{\mu s}),\sqrt{\mu} \right) \leq \frac{1}{2}\sqrt{\mu}$, resulting in an even slower rate than that of \eqref{eqn:parameters-zhang-cont}. 
\end{remark}

\noindent We also demonstrate significant improvements in the discrete-time analysis. 
First, we review the existing result. 

\begin{theorem}[{[\citenum{zhang2021revisiting}, Theorem 3]}] \label{theorem3_zhang}
    Let $f\in \mathcal{F}_{\mu,L}$, and suppose $(u_k)_{k=1}^{\infty}$ is the sequence generated by the SIE discretization of (\ref{GM-ODE}) with a step-size of $\sqrt{s}$. 
    If $0 \leq m' \leq \frac{1}{2L}$, $0< n's\leq m'\sqrt{s}$, and $q'\sqrt{s}\leq 1/2$, then there exists a constant $C>0$ such that, for any $k\geq 1$, the following holds: 
    \begin{align*}
        f(u_k)-f(x^*) \leq C(1+\gamma_2\sqrt{s})^{-k},
        \quad \text{with} \quad
        \gamma_2 = \tfrac{1}{5}\min\left(  \tfrac{n'\mu}{q'},\tfrac{q'}{(1+q'^2/(n'L))} \right).
    \end{align*}
\end{theorem}

\noindent By using the parameters in \eqref{eqn:parameters-zhang-disc}, this implies a convergence rate of 
\begin{align*}
    f(u_{k})-f(x^*) \leq \mathcal{O} \big(1+\sqrt{ 1/\kappa}/30)^{-k}\big), \quad \text{\textit{assuming} $\kappa \geq 9$.}
\end{align*}

\noindent Our analysis with \eqref{GM2-ODE} not only removes the restriction on $\kappa$, but also improves the convergence rate to $\mathcal{O}((1-\sqrt{ 1/\kappa})^k)$, as shown in \Cref{theorem4}. 
Notably, we used the same parameters, given in \eqref{eqn:param-gm2ode-nag}, in both the continuous and discrete-time settings. 

\subsection{SIE exactly recovers NAG}\label{sec:SIE_rec_NAG}
In a recent analysis, \eqref{NAG_ODE} was studied through a phase space representation as:
\begin{gather}
    \begin{aligned}\label{eqn:Shi_PH_repres}
            \left\{\begin{array}{ll}
                \dot Q_t = J_t,   \\ 
                \dot J_t = -(2\sqrt{\mu}+\sqrt{s}\nabla^2 f(Q_t))J_t-(1+\sqrt{\mu s})\nabla f(Q_t),
            \end{array}\right.
    \end{aligned}
\end{gather}
where $Q_t:=Q(t) $ and $ J_t:=J(t)$ are known as the momentum term in the literature \citep{WibisonoE7351,Shi2021UnderstandingTA}. Then, they studied the convergence behaviour of various integrators (SIE, explicit Euler, and implicit Euler). As expressed by the authors, none of the integrators successfully recovered the NAG algorithm. In addition, convergence rate of the resulting methods were worse than the NAG method's. Here, we will justify this mismatch, provide a solution to fix it, and show that SIE discretization of the \eqref{GM2-ODE} recovers the NAG method.   

\noindent To begin, reformulate \eqref{proposed_ODE_conti} as
\begin{gather}
    \begin{aligned}\label{proposed_ODE_conti3}
            \left\{\begin{array}{ll}
                \dot Q_t = J_t,   \\ 
                \dot J_t = -((n+q)+m\nabla^2 f(Q_t))J_t-(np+mq)\nabla f(Q_t) .
            \end{array}\right.
    \end{aligned}
\end{gather}
This is the same phase space representation used by \citet{shi2019acceleration} to study various discretizations. 

\noindent Due to the presence of a Hessian term, naively discretizing (\ref{eqn:Shi_PH_repres}) will lead to a second-order algorithm. This issue was resolved by \citet{Shi2021UnderstandingTA} through approximating the Hessian as
$$\sqrt{s}\nabla^2f(X_t)\dot{X}_t\approx   \nabla f(x_{k+1})-\nabla f(x_k).$$
Utilizing this approach and discretizing \eqref{proposed_ODE_conti3} gives
\begin{gather}
    \begin{aligned}\label{proposed_alg_1}
        \left\{
        \begin{array}{l l}
            q_{k+1}-q_k =  j_k\sqrt{s},  &\\
            j_{k+1}-j_k  = -(n+q)\sqrt{s}j_{k+1}
            &-\sqrt{s}(np+mq)\nabla f(q_{k+1})-m(\nabla f(q_{k+1})-\nabla f(q_k)).
        \end{array}
        \right.
    \end{aligned}
\end{gather}
Now, when ${n=\sqrt{\mu},q=\sqrt{\mu},p=1/\sqrt{\mu},m=\sqrt{ s},}$ \eqref{proposed_ODE_conti3} becomes \eqref{NAG_ODE}, while \eqref{proposed_alg_1} does NOT result in the NAG algorithm. We believe this is due to the aforementioned Hessian approximation. One may choose $${n=\frac{\sqrt{\mu}}{1-\sqrt{\mu s}},q=\frac{\sqrt{\mu}}{1-\sqrt{\mu s}},p=\frac{1}{\sqrt{\mu}},m=\sqrt{ s}},$$ so that (\ref{proposed_alg_1}) coincides with the NAG method. However, in this case, \eqref{proposed_ODE_conti3} no longer relates to \eqref{NAG_ODE}. This is what we refer to as \emph{coefficient inconsistency}. A more elegant way to solve the issue is through a Hessian-free reformulation of \eqref{NAG_ODE}. This is what we already have in \eqref{GM2-ODE}. Applying the SIE discretization to \eqref{proposed_ODE_conti} results in \eqref{SIE_newalg}
where both reduce to \eqref{NAG_ODE} and the NAG algorithm by setting 
$n=\sqrt{\mu},q=\sqrt{\mu},p=1/\sqrt{\mu},m=\sqrt{ s}.$ Thus, SIE discretization exactly recovers the NAG method.

\noindent This paves the way for a more consistent Lyapunov analysis and essentially, better convergence rates than \citep{shi2019acceleration}. Specifically, in discrete-time, \Cref{theorem4} proves the convergence rate
$${f(x_{k})-f(x^*)\leq C''_{GM}(1-\sqrt{ 1/\kappa})^k},$$ for the NAG's method
which is faster than the rate $f(x_{k})-f(x^*)\leq C'(1+\sqrt{ 1/\kappa}/9)^{-k}$ proven for the numerical discretization in \citep{Shi2021UnderstandingTA}. 

\subsection{Continuous-Time Analysis of Quadratics}\label{app2}
Consider the $\mu$-strongly convex function $f : \mathbb{R}^d\rightarrow \mathbb{R}$ to be a quadratic of the form $f(X)=\tfrac{1}{2}X^T\boldsymbol A X$. Then the second-order ODE (\ref{proposed_ODE_conti3}) can be rewritten as
\begin{align} \label{quadratic_continuous}
\left[
\begin{array}{c}
     \dot{X}  \\
      \dot{V}
\end{array}\right]=\left[\begin{array}{cc}
    \boldsymbol 0_d & \boldsymbol I_d \\
    -(np+mq)\boldsymbol A & -(n+q)\boldsymbol I_d-m \boldsymbol A
\end{array} \right]   \left[\begin{array}{c}
     X  \\
     V
\end{array}\right],
\end{align}
where $\boldsymbol 0_d$ is $d \times d$ zero matrix, and $\boldsymbol I_d$ is $d \times d$ identity matrix. Without loss of generality, we assume that $\boldsymbol A$ is a diagonal matrix and its entries are sorted decreasingly. The equation (\ref{quadratic_continuous}) is solved by taking $Z = [X^T \quad V^T]^T$ and 
$$\boldsymbol M=-\left[\begin{array}{cc}
    \boldsymbol 0_d & -\boldsymbol I_d \\
    (np+mq)\boldsymbol A & (n+q)\boldsymbol I_d+m\boldsymbol A
\end{array} \right],  $$
so that $\dot Z =-\boldsymbol MZ$ which has the well-known solution $Z=e^{-\boldsymbol Mt}Z_0$ for an initialization $Z_0$. Note that we have
\begin{align*}
    \|Z\|_2& \leq \|e^{-\boldsymbol Mt}\|_{2\rightarrow 2}\|Z_0\|_2\leq\sum_{k=0}^{\infty}\frac{\|(-\boldsymbol M)^k\|t^k}{k!}\|Z_0\|_2\leq\sum_{k=0}^{\infty} \frac{(\rho(-\boldsymbol M)+o(1))^kt^k}{k!}\|Z_0\|_2\\
    &\leq e^{(\rho(-\boldsymbol M+o(1)))t}\|Z_0\|_2
\end{align*}
where $\rho(\boldsymbol M)$ is the spectral radius of $\boldsymbol M$, $\|.\|_{2\rightarrow 2}$ denotes the spectral norm, and the last inequality is true asymptotically as $k\rightarrow \infty$ since ($\|A^k\|\leq (\rho(A)+o(1))^k$). To find the convergence rate we need maximum eigenvalue of $ -\boldsymbol M$ (minimum eigenvalue of $\boldsymbol M$) which corresponds to the largest spectral radius of ($ -\boldsymbol M$). 
Matrix $\boldsymbol M$ is
\begin{align}
\left[
\begin{array}{cc}
 \begin{array}{rrr}
     0 &  &  \\
     & \ddots & \\
     &  &  0
\end{array}     &  \begin{array}{rrr}
     -1 &  &  \\
     & \ddots & \\
     &  &  -1
\end{array} \\
    \begin{array}{rrr}
     (np+mq)a_{11} &  &  \\
     & \ddots & \\
     &  &    (np+mq)a_{dd}
\end{array}  &  \begin{array}{rrr}
     (n+q)+ma_{11} &  &  \\
     & \ddots & \\
     &  &   (n+q)+ma_{dd}
\end{array}
\end{array}
\right],\nonumber
\end{align}
which after permuting its rows and columns becomes
\begin{align}
\left[
\begin{array}{ccc}
 \begin{array}{cc}
     0 &  -1 \\
   (np+mq)a_{11} &  (n+q)+ma_{11}
\end{array}     &  \hdots &  \begin{array}{cc}
     0 &  0  \\
    0 & 0
\end{array}\\ 
\vdots & \ddots & \vdots \\
      \begin{array}{rr}
     0 &  0  \\
    0 & 0
\end{array} & \hdots & \begin{array}{cc}
     0 &  -1  \\
    (np+mq)a_{dd} &  (n+q)+ma_{dd}
\end{array}
\end{array}
\right],\nonumber
\end{align}
such that $a_{11}\geq a_{22}\geq \ldots\geq a_{dd}$. Due to the $\mu$-strong convexity of $f$, we have $a_{dd}=\mu$. Next, the eigenvalues of each $2\times 2 $ matrix in the block matrix $\boldsymbol M$ will lead to the eigenvalues of the whole matrix. The matrix $\boldsymbol M$ has $d$ blocks and each block has 2 eigenvalues. The eigenvalues of $i$-th block are noted with $ \lambda^i_{1,2}$ for $i\in\{1,..,d\}$. This will lead to
\begin{align*}
    \lambda^i_{1,2} = \tfrac{1}{2}\left(ma_{ii} + (n+q) \pm \sqrt{((ma_{ii})+(n+q))^2-4a_{ii}(np+mq)}\right)\quad \forall i\in\{1,..,d\}.
\end{align*}
Now, taking $n=q,p=\frac{q}{a_{ii}}+\frac{(ma_{ii})^2}{4qa_{ii}}$ results in the critical damping setting i.e.
\begin{align*}
    \sqrt{((ma_{ii})+(n+q))^2-4a_{ii}(np+mq)}=0\qquad \forall i\in\{1,..,d\}.
\end{align*}
Note that under this setting, all the eigenvalues are real and nonnegative. Since $m\geq 0$, choosing $a_{dd}=\mu$ will lead to the smallest eigenvalue (slowest one in convergence) which is
\begin{align*} 
    \lambda^d_{1,2} = \tfrac{1}{2}\left(m\mu + (2q) \right) = \tfrac{m\mu}{2}  +q .
\end{align*}
The analysis above gives us
\begin{align*} 
    \|Z_t\|_2\leq e^{-(\tfrac{m\mu}{2}  +q )t}\|Z_0\|_2,
\end{align*}
and 
\begin{align*}
    f(X_t)\leq \tfrac{\|A\|_{2\rightarrow 2}}{2}\|X_t\|^2_2\leq\tfrac{\|A\|_{2\rightarrow 2}}{2}\|Z_t\|^2_2\leq\tfrac{\|A\|_{2\rightarrow 2}}{2}e^{-2(\tfrac{m\mu}{2}  +q )t}\|Z_0\|^2_2.
\end{align*}
In particular, note that increasing $q$ and $m$ can lead to arbitrary fast convergence with the rate of $e^{-(m\mu +2q )t}$ under the conditions mentioned.
Note that this rate is slightly faster than the rate of the TM method and the best possible rate through Theorem \ref{theorem2}.

\subsection{Discrete-Time Analysis of the Quadratics}\label{app3}
We consider discretizing (\ref{proposed_ODE_conti}) and investigate the convergence behaviour of it for $\mu$-strongly convex $L$-smooth quadratic function of the form $f(X)=\frac{1}{2}X^T\boldsymbol AX$. Applying the SIE discretization on (\ref{proposed_ODE_conti}) we get
\begin{align} \label{proposed_alg_2_quad}
\left\{
\begin{array}{ll}
    x_{k+1}-x_k = & -m\sqrt{s}\boldsymbol Ax_k -n\sqrt{s}(x_{k+1}-v_k), \\
   v_{k+1}-v_k  =& -p\sqrt{s}\boldsymbol Ax_{k+1}-q\sqrt{s}(v_k-x_{k+1}).
\end{array}
    \right.
\end{align}
Without loss of generality we can assume that $\boldsymbol A$ is a diagonal matrix in which case the diagonal elements of $\boldsymbol A$ are its eigenvalues. The one line representation of (\ref{proposed_alg_2_quad}) is
\begin{align} \label{one_line_rep_alg2}
x_{k+1}^i&=\left(1+\frac{(1-q\sqrt{s}-npsa_{ii}-m\sqrt{s}a_{ii})}{1+n\sqrt{s}}\right)x_k^i +\left(\frac{m\sqrt{s}a_{ii}-1+q\sqrt{s}(1-m\sqrt{s}a_{ii})}{1+n\sqrt{s}}\right)x_{k-1}^i,
\end{align}
where upper index $i$ denotes the $i$'th element and $a_{ii}$ is the $i$'th element of $\boldsymbol A$'s diagonal elements for $i=1,\ldots,d$. For comparison, the one line representation of the NAG algorithm for quadratic function $f(X)$ is 
$$x_{k+1}^i=(\frac{2}{1+\sqrt{\mu s}}(1-sa_{ii}))x_k^i+\frac{1-\sqrt{\mu s}}{1+\sqrt{\mu s}}(sa_{ii} -1)x_{k-1}^i\qquad \forall i\in \{1,\ldots d\},$$
which can be derived from (\ref{one_line_rep_alg2}) by setting
$$n=\sqrt{\mu},q=\sqrt{\mu},p=\frac{1}{\sqrt{\mu}},m=\sqrt{ s}.$$
To study the convergence rate of (\ref{proposed_alg_2_quad}), we reformulate (\ref{one_line_rep_alg2}) as
\begin{align*}
    \begin{aligned}
        y_k&=\left[\begin{array}{c}
            x_{k+1}   \\
            x_k  
        \end{array}\right]\\
        &=\left[
        \begin{array}{cc}
           \left((1+\frac{(1-q\sqrt{s})}{1+n\sqrt{s}})\boldsymbol I_d-\frac{(nps+m\sqrt{s})}{1+n\sqrt{s}}\boldsymbol A\right)  & \left(\frac{(1-q\sqrt{s})(\boldsymbol Am\sqrt{s}-\boldsymbol I_d)}{1+n\sqrt{s}}\right) \\
            \boldsymbol I_d &\boldsymbol 0_d
        \end{array}
        \right]\left[\begin{array}{c}
            x_{k}   \\
            x_{k-1} 
        \end{array}\right]\\
        &=\boldsymbol Ty_{k-1}
    \end{aligned}
\end{align*}
with $\boldsymbol 0_d$ and $\boldsymbol I_d$ as a $d\times d$ zero matrix and d-dimension identity matrix. Next, we have
$$\|y_k\|_2=\|\boldsymbol Ty_{k-1}\|_2=\|\boldsymbol T^ky_0\|_2\leq \|\boldsymbol T^k\|_{2}\|y_0\|_2\leq (\rho(\boldsymbol T)+o(1))^k\|y_0\|_2,$$
where $\rho(\boldsymbol T)$ is the spectral radius of $\boldsymbol T$ and the last inequality is true asymptotically as $k\rightarrow \infty$ through Gelfand's formula \citep{horn2012matrix}. To define convergence rate we need the largest eigenvalue of $\boldsymbol T$ which corresponds to the largest spectral radius. Note that $\boldsymbol T$ is the block diagonal matrix 
\begin{align*}
    \boldsymbol T&=\left[ 
    \begin{array}{cccc}
      \boldsymbol T_1  & 0 & \ldots & 0 \\
        0 &\boldsymbol T_2  & \ldots & 0 \\
        \vdots & \vdots & \ddots & \vdots\\
        0 & 0 &\ldots & \boldsymbol T_d
    \end{array}
    \right],\\
    &\boldsymbol T_i= \left[\begin{array}{cc}
       \left((1+\frac{(1-q\sqrt{s})}{1+n\sqrt{s}})-\frac{(nps+m\sqrt{s})}{1+n\sqrt{s}}a_{ii} \right)  & \left(\frac{(1-q\sqrt{s})(a_{ii} m\sqrt{s}-1)}{1+n\sqrt{s}}\right) \\
        1 & 0
    \end{array}
    \right],
\end{align*}
for $i\in\{1,\ldots ,d\}$. Hence, the eigenvalues of $\boldsymbol T$ are the union of the eigenvalues of $\boldsymbol T_i$'s. For each $\boldsymbol T_i$ there exist two eigenvalues as the solutions of
$$r^2-\left(1+\frac{1-q\sqrt{s}}{1+n\sqrt{s}}-\frac{(nps+m\sqrt{s}a_{ii})}{1+n\sqrt{s}}\right)r-\frac{(1-q\sqrt{s})(a_{ii}m\sqrt{s}-1)}{1+n\sqrt{s}}=0,$$
with
$$\Delta = \left(1+\frac{1-q\sqrt{s}}{1+n\sqrt{s}}-\frac{(nps+m\sqrt{s}a_{ii})}{1+n\sqrt{s}}\right)^2-4\frac{(1-q\sqrt{s})(a_{ii}m\sqrt{s}-1)}{1+n\sqrt{s}}.$$

\noindent Taking $n=q=\sqrt{a_{ii}},np=1,m=\sqrt{s}$ leads to $\Delta=0$ and $r_{1,2}=1-\sqrt{sa_{ii}}$. With this choice of parameters, the convergence rate of (\ref{proposed_alg_2_quad}) for $\mu$-strongly convex and $L$-Lipschitz quadratics of the form $f(X)=\frac{1}{2} X^T\boldsymbol AX$ will be
$$f(x_k)-f(x^*)=\frac{1}{2}\sum_{i=1}^d a_{ii}(x_k^i)^2\leq \max_i a_{ii}\|x_k\|^2\leq C\max_i(1-\sqrt{sa_{ii}})^{2k},$$
for $\sqrt{s}\leq \tfrac{1}{\sqrt{L}}$ (due to $1-\sqrt{sa_{ii}}\geq 0$). The worst case scenario happens for $a_{ii}=\mu$ (closest possible rate to 1) which leads to the rate $O((1-\sqrt{\tfrac{1}{\kappa}})^{2k})$. 

\subsection{Convergence of Optimization Algorithms through Dynamical Systems}\label{app_subsec_QCI}
In the state space, dynamical systems are usually presented in the form of 
\begin{align}\label{statespace2}
    \dot\xi(t)=A\xi(t)+Bu(t),\quad y(t)=C\xi(t),\quad u(t)=\nabla f(y(t))\quad \forall t\geq 0 ,
\end{align}
where $\xi\in\mathbb{R}^{r}$ is the state, $y(t)\in \mathbb{R}^d\:(d\leq r)$ is the output, and $u(t)$ is the continuous feedback input. Here, we would have $u^*=0$ and the fixed point of (\ref{statespace2}) is
$$A\xi^*=0,\quad y^*=C\xi^*.$$
Consider the nonnegative function 
$$\varepsilon(t)=e^{\lambda t}\left(f(y(t))-f(y^*)+(\xi(t)-\xi^*)^TP(\xi(t)-\xi^*)\right),$$
with $\lambda>0$, $y^*=x^*$ and $P\succeq 0$ where $A\succeq B$ denotes that $A-B$ is positive semi-definite. If when $\xi\rightarrow \xi^*$ we have $\tfrac{d}{dt}\varepsilon(t)\leq 0$, then $\varepsilon(t)\leq \varepsilon(0)$. This results in
$$f(y(t))-f(y^*)\leq e^{-\lambda t}\varepsilon(0).$$
The following result from \citep{doi:10.1137/17M1136845} proposes a Linear Matrix Inequality (LMI) that guarantees the existence of a Lyapunov function through which we can show that $f(x)$ converges exponentially fast. For simplicity, we adopt the presentation of \citep{doi:10.1137/20M1364138}.
\begin{theorem}[Theorem 6.4 in \citep{doi:10.1137/17M1136845}]\label{theorem1}
Suppose that for (\ref{statespace2}) there exists $\lambda>0,P\succeq 0$ and $\sigma\geq 0$ such that
$\boldsymbol T=M^{(0)}+M^{(1)}+\lambda M^{(2)}+\sigma M^{(3)}\preceq 0$
where 
\begin{align}
M^{(0)}&=\left[
    \begin{array}{cc}
     PA+A^TP+\lambda P    & PB \\
        B^TP & 0 
    \end{array}\right] ,\nonumber\\
M^{(1)}&=\frac{1}{2}\left[
    \begin{array}{cc}
         0&(CA)^T  \\
         CA&CB+(CB)^T 
    \end{array}\right] ,\nonumber\\
    M^{(2)}&=\left[
    \begin{array}{cc}
        C^T & 0 \\
        0 & I_d
    \end{array}\right]\left[
    \begin{array}{cc}
        -\frac{\mu}{2}I_d & \frac{1}{2}I_d \\
        \frac{1}{2}I_d & 0
    \end{array}\right]\left[
    \begin{array}{cc}
        C & 0 \\
        0 & I_d
    \end{array}\right] ,\nonumber\\
    M^{(3)}&=\left[
    \begin{array}{cc}
        C^T & 0 \\
        0 & I_d
    \end{array}\right]\left[
    \begin{array}{cc}
        -\frac{\mu L}{\mu+L}I_d & \frac{1}{2}I_d \\
        \frac{1}{2}I_d & \frac{-1}{\mu+L}
    \end{array}\right]\left[
    \begin{array}{cc}
        C & 0 \\
        0 & I_d
    \end{array}\right] ,\nonumber
\end{align}
and $(.)^T$ denotes the transpose operator and $I_d$ is the identity matrix of size $d$. Then for $\mu$-strongly convex $L$-smooth $f$ we have
\begin{align}
  f(y(t))-f(y^*)&\leq e^{-\lambda t}\varepsilon(0).  \nonumber
\end{align}
\end{theorem}
\subsection{Acceleration of the Explicit Euler discretization}\label{appendix_EE}
We would like to show the correspondence between EE and SIE discretizations of (\ref{proposed_ODE_conti}). The following lemma shows how to update the coefficients of EE method such that SIE is derived.
\begin{lemma}\label{lem3}
Consider the parameters of EE discretization of (\ref{proposed_ODE_conti}) as $n_{EE},m_{EE},q,p$ and the parameters of SIE discretization of (\ref{proposed_ODE_conti}) as $n_{SIE},$ $m_{SIE} ,q,p.$ Then by taking
\begin{align} \label{param_update}
    \left\{
\begin{array}{c}
     n_{EE}=\frac{n_{SIE}-qn_{SIE}\sqrt{s}}{1+n_{SIE}\sqrt{s}}, \\
      m_{EE}=\frac{m_{SIE}+n_{SIE}p\sqrt{s}}{1+n_{SIE}\sqrt{s}},
\end{array}
    \right.
\end{align}
SIE discretization of (\ref{proposed_ODE_conti}) will be the same as its EE discretization with step-size $\sqrt{s}$.
\end{lemma}
\noindent For proving the result, note that the EE discretization of (\ref{proposed_ODE_conti}) is \begin{align} \label{EE_newalg}
\left\{
\begin{array}{ll}
    x_{k+1}-x_k = &  -m\sqrt{s}\nabla f(x_k)-n\sqrt{s}(x_k-v_k), \\
   v_{k+1}-v_k  =& -p\sqrt{s}\nabla f(x_{k_{k}})-\sqrt{s}q(v_k-x_{k}) ,
\end{array}\right.
\end{align}
which can be written in one line format
\begin{align}
   &x_{k+1}=x_k-m\sqrt{s}\nabla f(x_k)+(1-q\sqrt{s}-n\sqrt{s})(x_k-x_{k-1})+(m\sqrt{s}(1-q\sqrt{s})-nps)\nabla f(x_{k-1}),\nonumber
\end{align}
replacing the coefficient updates from (\ref{param_update}) in above gives the SIE one line update of (\ref{SIE_newalg}) which is
\begin{align*}
    x_{k+1}&=x_k-\tfrac{m\sqrt{s}+nps}{1+n\sqrt{s}}\nabla f(x_k)+\tfrac{1-q\sqrt{s}}{1+n\sqrt{s}}(x_k-x_{k-1})+\tfrac{m\sqrt{s}(1-q\sqrt{s})}{1+n\sqrt{s}}\nabla f(x_{k-1}),
\end{align*}

\noindent With Lemma \ref{lem3} and Theorem \ref{theorem4}, we establish the convergence result for (\ref{EE_newalg}) as follows.
\begin{corollary}[Convergence of (\ref{EE_newalg})]\label{coro_EE}
    For $\mu$-strongly convex $L$-smooth function $f$ with $0\leq \mu< L$ and parameters 
    $m,n,p,q$ such that
    $$q/p\leq \mu,0\leq qps\leq m\sqrt{s}(1+q\sqrt{s})-qps\leq \frac{1}{L},n=q\frac{1-q\sqrt{s}}{1+q\sqrt{s}}, 0\leq q\sqrt{s}<1,p>0,$$
    the sequence $x_k$ in (\ref{EE_newalg}) will satisfy 
    $$f(x_{k})-f(x^*)\leq LC'_{GM}(1-q\sqrt{s})^k,$$
for constant $C'_{GM}>0$ and any $x_0,v_0=x_0-(\frac{m}{n}-\frac{2qps}{n\sqrt{s}(1+q\sqrt{s})})\nabla f(x_0)$.
\end{corollary}
\noindent The proof is simply done by using (\ref{param_update}) in Theorem \ref{theorem4}. Note that for the initial condition we need to get the same result as in Theorem \ref{theorem4} with the new coefficients and the new update rule (\ref{EE_newalg}).

\noindent Corollary \Cref{coro_EE} suggests that choosing 
$$q=\sqrt{\mu},p=\frac{1}{\sqrt{\mu}},m=\frac{2s}{1+\sqrt{\mu s}}, n=\sqrt{\mu}\frac{1-\sqrt{\mu s}}{1+\sqrt{\mu s}},$$
in (\ref{EE_newalg}) will recover the NAG algorithm.
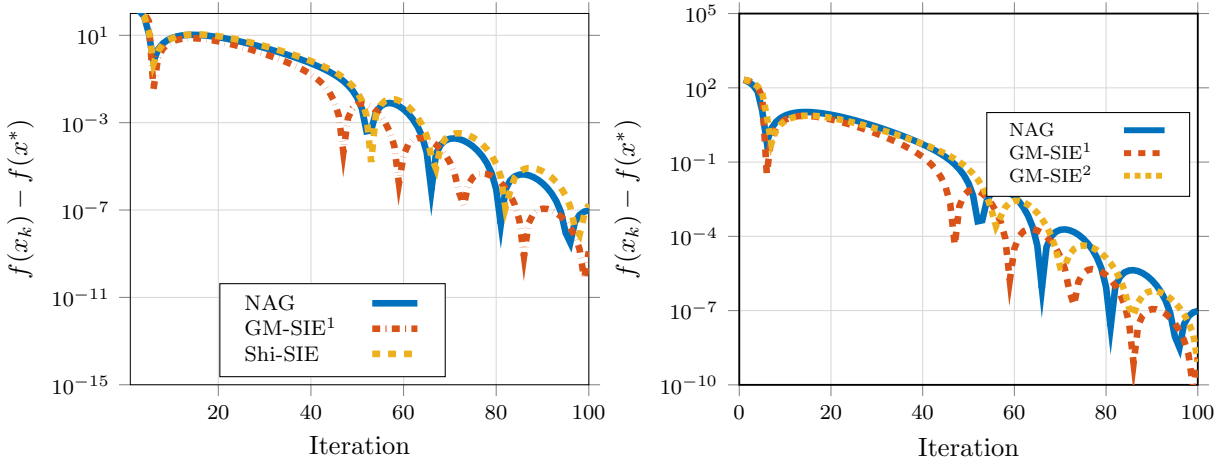
\begin{figure}
    \centering
\definecolor{mycolor1}{rgb}{0.00000,0.44700,0.74100}
\definecolor{mycolor2}{rgb}{0.85000,0.32500,0.09800}
\definecolor{mycolor3}{rgb}{0.92900,0.69400,0.12500}
\begin{tikzpicture}

\begin{groupplot}[
        group style={group name=fig3,group size=2 by 1, horizontal sep=2cm
        },   
        width=0.465\textwidth, 
        height=0.395\textwidth,        
        grid=both, grid style={gray!30},        
        tick label style={font=\footnotesize}, 
        ]
    
\nextgroupplot[
    ymin=1e-15, ymax=1e2,
    xmin=1e0, xmax=1e2,  
    ylabel={$f(x_k)-f(x^*)$}, 
    xlabel={Iteration},
    line width=0.1pt,
    xtick align=outside,
    ytick align=outside,
    ymode=log,
    ]
\addplot [color=mycolor1, line width=2.5pt]
  table[row sep=crcr]{
1	210.095562436466\\
2	174.267152871458\\
3	126.941603897643\\
4	71.7324572481976\\
5	11.9346632646616\\
6	0.643617944910188\\
7	2.26593197483689\\
8	4.16391748760623\\
9	5.9993287826018\\
10	7.58574489921796\\
11	8.83676333006101\\
12	9.72960808869808\\
13	10.2799278421884\\
14	10.5246609906501\\
15	10.510662481992\\
16	10.2874003299502\\
17	9.90248779423209\\
18	9.39915784486114\\
19	8.81503875060043\\
20	8.18177538799921\\
21	7.52517685848794\\
22	6.86566987281051\\
23	6.21890866856726\\
24	5.59644316304305\\
25	5.00638301095583\\
26	4.45402027476312\\
27	3.94239053244771\\
28	3.47276367498876\\
29	3.04506302830117\\
30	2.658215978608\\
31	2.31044186443264\\
32	1.99948415734652\\
33	1.72279434185219\\
34	1.47767474501624\\
35	1.26138708629761\\
36	1.07123287816467\\
37	0.904611123522105\\
38	0.759058109831239\\
39	0.632273548677836\\
40	0.522136878513104\\
41	0.42671721733901\\
42	0.344280139797628\\
43	0.273294005761856\\
44	0.212437754916451\\
45	0.160610585467698\\
46	0.116941312205603\\
47	0.0807908716467565\\
48	0.0517349786586425\\
49	0.029506287951665\\
50	0.0138718807635229\\
51	0.00443581272405988\\
52	0.000405747876149803\\
53	0.000443905761938007\\
54	0.00276726365101393\\
55	0.00555588963231965\\
56	0.00747780504603743\\
57	0.00799103331403131\\
58	0.00725272359162116\\
59	0.00577645790523598\\
60	0.00409747662463517\\
61	0.00259112757962357\\
62	0.00143878391923521\\
63	0.000672168914798099\\
64	0.000236521269773193\\
65	4.27581116657916e-05\\
66	8.05116577823339e-07\\
67	3.6506961688304e-05\\
68	9.72987949741377e-05\\
69	0.000151359939729617\\
70	0.000183643291051389\\
71	0.000190925578930329\\
72	0.000177101601337187\\
73	0.000149332143164621\\
74	0.000115268107690608\\
75	8.13367394512654e-05\\
76	5.1937813102948e-05\\
77	2.9327401049134e-05\\
78	1.39465069800282e-05\\
79	4.96931433602299e-06\\
80	8.89700076167621e-07\\
81	2.27064870794713e-08\\
82	8.57408629567624e-07\\
83	2.24858649885973e-06\\
84	3.47024628627457e-06\\
85	4.17246849349073e-06\\
86	4.28649264580527e-06\\
87	3.91603779278937e-06\\
88	3.24087839306819e-06\\
89	2.4460937029791e-06\\
90	1.68006966480561e-06\\
91	1.03750882229914e-06\\
92	5.60383039249768e-07\\
93	2.4919075217178e-07\\
94	7.80374178055432e-08\\
95	8.99711341539877e-09\\
96	3.22838580868456e-09\\
97	2.7985179507084e-08\\
98	5.98074949564875e-08\\
99	8.4795231919399e-08\\
100	9.63718636837196e-08\\
};\label{figNag3}

\addplot [color=mycolor2, dashdotted, line width=2.5pt]
  table[row sep=crcr]{
1	210.095562436466\\
2	192.310448333424\\
3	155.211091896182\\
4	104.276550830535\\
5	44.2253111777479\\
6	0.0336653136927416\\
7	1.00105716395343\\
8	2.42635112215863\\
9	3.882010749747\\
10	5.1545785197657\\
11	6.14707568252227\\
12	6.83369658335104\\
13	7.2295977202392\\
14	7.37118767254353\\
15	7.30364153914042\\
16	7.07331816083508\\
17	6.72344776634928\\
18	6.29195296585163\\
19	5.81061980792334\\
20	5.30508665891702\\
21	4.79529542149194\\
22	4.29617292742777\\
23	3.81839543247452\\
24	3.3691470687553\\
25	2.95282187982538\\
26	2.57164444266347\\
27	2.22620020156702\\
28	1.91587644259716\\
29	1.63922040301953\\
30	1.39422379665988\\
31	1.1785440594262\\
32	0.989672584891191\\
33	0.825059618297787\\
34	0.682204648895392\\
35	0.558720308357574\\
36	0.452377054995153\\
37	0.361135269680054\\
38	0.283170614664236\\
39	0.216897233958728\\
40	0.160991026831913\\
41	0.11441096537433\\
42	0.0764091029924071\\
43	0.0465085702070367\\
44	0.0244152900393521\\
45	0.00982466940047605\\
46	0.00211771322490945\\
47	4.46516842905198e-05\\
48	0.00163448391721482\\
49	0.00456445500960301\\
50	0.00689661155254484\\
51	0.00767895504621979\\
52	0.00697271818035508\\
53	0.00539825978246053\\
54	0.00362797620762034\\
55	0.00211080066519534\\
56	0.00103194129288026\\
57	0.000388880948158615\\
58	8.55171654220166e-05\\
59	1.70082730363563e-06\\
60	3.17894311252231e-05\\
61	9.97576802620748e-05\\
62	0.000160298912705109\\
63	0.000193158037306879\\
64	0.000195274474078799\\
65	0.000173270969659606\\
66	0.000137486022265543\\
67	9.79323553117473e-05\\
68	6.20880795473289e-05\\
69	3.416805249154e-05\\
70	1.54110138688157e-05\\
71	4.91377969649709e-06\\
72	6.17030762595228e-07\\
73	1.66639491872012e-07\\
74	1.50662734027529e-06\\
75	3.17504522337475e-06\\
76	4.353092092374e-06\\
77	4.7550876305924e-06\\
78	4.44887926051574e-06\\
79	3.67596705450057e-06\\
80	2.7119724739233e-06\\
81	1.78161990366243e-06\\
82	1.02377429681e-06\\
83	4.92502746057832e-07\\
84	1.77992959621287e-07\\
85	3.37122864624284e-08\\
86	7.55407764119198e-10\\
87	2.4952412081003e-08\\
88	6.59030695293605e-08\\
89	9.9306634082641e-08\\
90	1.14886703983519e-07\\
91	1.12218752651438e-07\\
92	9.62550244787419e-08\\
93	7.36478176777133e-08\\
94	5.0329638817681e-08\\
95	3.03377770960456e-08\\
96	1.56025845854391e-08\\
97	6.32636423802957e-09\\
98	1.60914667701739e-09\\
99	7.06291414243054e-11\\
100	1.2720110242892e-09\\
};\label{figGMSIE12}

\addplot [color=mycolor3, dashed, line width=2.5pt]
  table[row sep=crcr]{
1	210.095562436466\\
2	174.267152871458\\
3	127.980266203922\\
4	74.3545003883637\\
5	16.1535214325929\\
6	0.490294572156493\\
7	1.97157294565933\\
8	3.79345197592255\\
9	5.63156041503103\\
10	7.2917008807982\\
11	8.66997933130001\\
12	9.72352197571708\\
13	10.4491921897046\\
14	10.8683429999437\\
15	11.0160876184419\\
16	10.933921665481\\
17	10.664805628456\\
18	10.250030469102\\
19	9.72735573219248\\
20	9.13003813616659\\
21	8.48646752131536\\
22	7.82020263456053\\
23	7.1502566594202\\
24	6.49152571000438\\
25	5.85528590043916\\
26	5.2497085953561\\
27	4.68036101815712\\
28	4.15067208398013\\
29	3.66235233274037\\
30	3.21576309215645\\
31	2.81023421551131\\
32	2.44433246067944\\
33	2.11608422140244\\
34	1.82315720303525\\
35	1.56300598904181\\
36	1.3329864485938\\
37	1.13044372429899\\
38	0.952778218007156\\
39	0.797493648774435\\
40	0.662230962482316\\
41	0.544791677353198\\
42	0.443154160715038\\
43	0.355486282037072\\
44	0.280157700148391\\
45	0.215754410385839\\
46	0.161096629813816\\
47	0.115257939690241\\
48	0.0775778454671363\\
49	0.0476505473511532\\
50	0.0252603944573988\\
51	0.01022637793234\\
52	0.0021351022441711\\
53	1.59277894344889e-05\\
54	0.00214247143292806\\
55	0.00621504527362621\\
56	0.00999048374764144\\
57	0.0120156071220268\\
58	0.011945628657265\\
59	0.0102845764447449\\
60	0.00785469186159757\\
61	0.00537403842419354\\
62	0.00328237254977562\\
63	0.00175188051348202\\
64	0.000772873713916128\\
65	0.000242467256778151\\
66	2.8731768621254e-05\\
67	7.99936665890622e-06\\
68	8.19543216838692e-05\\
69	0.000182048957193759\\
70	0.000267049584158463\\
71	0.000317514579977196\\
72	0.000329462344018094\\
73	0.000308462357188954\\
74	0.000264750865485469\\
75	0.000209595266945695\\
76	0.000152906130338137\\
77	0.000101960067661833\\
78	6.1017039011263e-05\\
79	3.15770829886763e-05\\
80	1.30177536180609e-05\\
81	3.3802176547526e-06\\
82	1.22101013133591e-07\\
83	7.1844438170876e-07\\
84	3.05644094461544e-06\\
85	5.62403982584669e-06\\
86	7.52986120586274e-06\\
87	8.41006355375984e-06\\
88	8.27925120916451e-06\\
89	7.37246363996524e-06\\
90	6.00976062450909e-06\\
91	4.49902365992627e-06\\
92	3.07966765072942e-06\\
93	1.90144623984922e-06\\
94	1.02840295862494e-06\\
95	4.57349076360414e-07\\
96	1.41793459557116e-07\\
97	1.48363061835965e-08\\
98	7.27312540749203e-09\\
99	5.948165679448e-08\\
100	1.83301257006274e-07\\
};\label{figshiSIE2}

\nextgroupplot[
        ymin=1e-10, ymax=1e5,   
        xmin=0, xmax=1e2,  
        ylabel={$f(x_k)-f(x^*)$},    
        xlabel={Iteration},
        line width=0.7pt,
        xtick align=outside,
        ytick align=outside,
        ymode = log,
        ]

\addplot [color=mycolor1, line width=2.5pt]
  table[row sep=crcr]{
1	210.095562436466\\
2	174.267152871458\\
3	126.941603897643\\
4	71.7324572481976\\
5	11.9346632646616\\
6	0.643617944910188\\
7	2.26593197483689\\
8	4.16391748760623\\
9	5.9993287826018\\
10	7.58574489921796\\
11	8.83676333006101\\
12	9.72960808869808\\
13	10.2799278421884\\
14	10.5246609906501\\
15	10.510662481992\\
16	10.2874003299502\\
17	9.90248779423209\\
18	9.39915784486114\\
19	8.81503875060043\\
20	8.18177538799921\\
21	7.52517685848794\\
22	6.86566987281051\\
23	6.21890866856726\\
24	5.59644316304305\\
25	5.00638301095583\\
26	4.45402027476312\\
27	3.94239053244771\\
28	3.47276367498876\\
29	3.04506302830117\\
30	2.658215978608\\
31	2.31044186443264\\
32	1.99948415734652\\
33	1.72279434185219\\
34	1.47767474501624\\
35	1.26138708629761\\
36	1.07123287816467\\
37	0.904611123522105\\
38	0.759058109831239\\
39	0.632273548677836\\
40	0.522136878513104\\
41	0.42671721733901\\
42	0.344280139797628\\
43	0.273294005761856\\
44	0.212437754916451\\
45	0.160610585467698\\
46	0.116941312205603\\
47	0.0807908716467565\\
48	0.0517349786586425\\
49	0.029506287951665\\
50	0.0138718807635229\\
51	0.00443581272405988\\
52	0.000405747876149803\\
53	0.000443905761938007\\
54	0.00276726365101393\\
55	0.00555588963231965\\
56	0.00747780504603743\\
57	0.00799103331403131\\
58	0.00725272359162116\\
59	0.00577645790523598\\
60	0.00409747662463517\\
61	0.00259112757962357\\
62	0.00143878391923521\\
63	0.000672168914798099\\
64	0.000236521269773193\\
65	4.27581116657916e-05\\
66	8.05116577823339e-07\\
67	3.6506961688304e-05\\
68	9.72987949741377e-05\\
69	0.000151359939729617\\
70	0.000183643291051389\\
71	0.000190925578930329\\
72	0.000177101601337187\\
73	0.000149332143164621\\
74	0.000115268107690608\\
75	8.13367394512654e-05\\
76	5.1937813102948e-05\\
77	2.9327401049134e-05\\
78	1.39465069800282e-05\\
79	4.96931433602299e-06\\
80	8.89700076167621e-07\\
81	2.27064870794713e-08\\
82	8.57408629567624e-07\\
83	2.24858649885973e-06\\
84	3.47024628627457e-06\\
85	4.17246849349073e-06\\
86	4.28649264580527e-06\\
87	3.91603779278937e-06\\
88	3.24087839306819e-06\\
89	2.4460937029791e-06\\
90	1.68006966480561e-06\\
91	1.03750882229914e-06\\
92	5.60383039249768e-07\\
93	2.4919075217178e-07\\
94	7.80374178055432e-08\\
95	8.99711341539877e-09\\
96	3.22838580868456e-09\\
97	2.7985179507084e-08\\
98	5.98074949564875e-08\\
99	8.4795231919399e-08\\
100	9.63718636837196e-08\\
};\label{figNag4}

\addplot [color=mycolor2, dashed, line width=2.5pt]
  table[row sep=crcr]{
1	210.095562436466\\
2	192.310448333424\\
3	155.211091896182\\
4	104.276550830535\\
5	44.2253111777479\\
6	0.0336653136927416\\
7	1.00105716395343\\
8	2.42635112215863\\
9	3.882010749747\\
10	5.1545785197657\\
11	6.14707568252227\\
12	6.83369658335104\\
13	7.2295977202392\\
14	7.37118767254353\\
15	7.30364153914042\\
16	7.07331816083508\\
17	6.72344776634928\\
18	6.29195296585163\\
19	5.81061980792334\\
20	5.30508665891702\\
21	4.79529542149194\\
22	4.29617292742777\\
23	3.81839543247452\\
24	3.3691470687553\\
25	2.95282187982538\\
26	2.57164444266347\\
27	2.22620020156702\\
28	1.91587644259716\\
29	1.63922040301953\\
30	1.39422379665988\\
31	1.1785440594262\\
32	0.989672584891191\\
33	0.825059618297787\\
34	0.682204648895392\\
35	0.558720308357574\\
36	0.452377054995153\\
37	0.361135269680054\\
38	0.283170614664236\\
39	0.216897233958728\\
40	0.160991026831913\\
41	0.11441096537433\\
42	0.0764091029924071\\
43	0.0465085702070367\\
44	0.0244152900393521\\
45	0.00982466940047605\\
46	0.00211771322490945\\
47	4.46516842905198e-05\\
48	0.00163448391721482\\
49	0.00456445500960301\\
50	0.00689661155254484\\
51	0.00767895504621979\\
52	0.00697271818035508\\
53	0.00539825978246053\\
54	0.00362797620762034\\
55	0.00211080066519534\\
56	0.00103194129288026\\
57	0.000388880948158615\\
58	8.55171654220166e-05\\
59	1.70082730363563e-06\\
60	3.17894311252231e-05\\
61	9.97576802620748e-05\\
62	0.000160298912705109\\
63	0.000193158037306879\\
64	0.000195274474078799\\
65	0.000173270969659606\\
66	0.000137486022265543\\
67	9.79323553117473e-05\\
68	6.20880795473289e-05\\
69	3.416805249154e-05\\
70	1.54110138688157e-05\\
71	4.91377969649709e-06\\
72	6.17030762595228e-07\\
73	1.66639491872012e-07\\
74	1.50662734027529e-06\\
75	3.17504522337475e-06\\
76	4.353092092374e-06\\
77	4.7550876305924e-06\\
78	4.44887926051574e-06\\
79	3.67596705450057e-06\\
80	2.7119724739233e-06\\
81	1.78161990366243e-06\\
82	1.02377429681e-06\\
83	4.92502746057832e-07\\
84	1.77992959621287e-07\\
85	3.37122864624284e-08\\
86	7.55407764119198e-10\\
87	2.4952412081003e-08\\
88	6.59030695293605e-08\\
89	9.9306634082641e-08\\
90	1.14886703983519e-07\\
91	1.12218752651438e-07\\
92	9.62550244787419e-08\\
93	7.36478176777133e-08\\
94	5.0329638817681e-08\\
95	3.03377770960456e-08\\
96	1.56025845854391e-08\\
97	6.32636423802957e-09\\
98	1.60914667701739e-09\\
99	7.06291414243054e-11\\
100	1.2720110242892e-09\\
};\label{figGMSIE1}

\addplot [color=mycolor3, dotted, line width=2.5pt]
  table[row sep=crcr]{
1	210.432534638163\\
2	192.551006016613\\
3	159.207462826889\\
4	114.698015876021\\
5	62.6695695791468\\
6	7.2397078152237\\
7	0.563326298300102\\
8	1.8534392147372\\
9	3.28768932347961\\
10	4.61622093229191\\
11	5.716232576669\\
12	6.54272539270462\\
13	7.09614808552883\\
14	7.40150989109017\\
15	7.49521875538816\\
16	7.41706559473942\\
17	7.20558829196232\\
18	6.89561468919847\\
19	6.51717570911527\\
20	6.0952498363338\\
21	5.64998521514121\\
22	5.19717142628723\\
23	4.74881777325988\\
24	4.31375138098826\\
25	3.8981854934072\\
26	3.50623225154009\\
27	3.14034923937474\\
28	2.80171816566059\\
29	2.49055923831591\\
30	2.20638746754354\\
31	1.94821823691032\\
32	1.71472963089259\\
33	1.50438861226963\\
34	1.31554746817818\\
35	1.14651615725251\\
36	0.995615394962181\\
37	0.86121457016122\\
38	0.741757924617854\\
39	0.635781860705137\\
40	0.54192576725371\\
41	0.458938354037583\\
42	0.385681134873781\\
43	0.321130361961393\\
44	0.264378345780283\\
45	0.214634642870394\\
46	0.171226994538578\\
47	0.133601080218179\\
48	0.101317051566404\\
49	0.0740394759994315\\
50	0.0515160626165021\\
51	0.0335402589327239\\
52	0.0198951423465242\\
53	0.0102829394683012\\
54	0.00425638541055817\\
55	0.00118065608456311\\
56	0.000257122905792218\\
57	0.000621235533577874\\
58	0.00148812534624609\\
59	0.00228524329576898\\
60	0.00271276516268443\\
61	0.00271394665440766\\
62	0.00238630217363459\\
63	0.0018840965738221\\
64	0.00134840980810499\\
65	0.000874593253390149\\
66	0.000508548145103344\\
67	0.000257971351253078\\
68	0.000107561775685078\\
69	3.21614893047939e-05\\
70	5.74420644477991e-06\\
71	6.36650771751457e-06\\
72	1.80693871067172e-05\\
73	3.08445937799162e-05\\
74	3.95892827150246e-05\\
75	4.26999855943078e-05\\
76	4.0707489801034e-05\\
77	3.51614921683674e-05\\
78	2.78398594763252e-05\\
79	2.02728513814532e-05\\
80	1.3526428497751e-05\\
81	8.17034144606765e-06\\
82	4.35715395638581e-06\\
83	1.95022491827901e-06\\
84	6.55869239502005e-07\\
85	1.32775205385904e-07\\
86	6.71292349485597e-08\\
87	2.13127866244278e-07\\
88	4.05268921171986e-07\\
89	5.51569778239891e-07\\
90	6.1673393583761e-07\\
91	6.02476767419535e-07\\
92	5.29772368196291e-07\\
93	4.25439230428148e-07\\
94	3.13657453426819e-07\\
95	2.11841357206088e-07\\
96	1.2972741927908e-07\\
97	7.04294077780165e-08\\
98	3.23859806194982e-08\\
99	1.14267944273383e-08\\
100	8.63525639882567e-10\\
};\label{figGMSIE2}
\end{groupplot}
\node[anchor=south east, draw = black, line width=0.4pt, fill=white, font=\footnotesize]  (legend) at ([shift={(-1.9cm,0.1cm)}]fig3 c1r1.south east) {\begin{tabular}{l l}
NAG & \ref*{figNag3} \\
GM-SIE$^1$ & \ref*{figGMSIE12}\\
Shi-SIE & \ref*{figshiSIE2}
\end{tabular}};

\node[anchor=south east, draw = black, line width=0.4pt, fill=white, font=\footnotesize,scale=0.9]  (legend) at ([shift={(-0.1cm,2.5cm)}]fig3 c2r1.south east) {\begin{tabular}{l l}
NAG & \ref*{figNag4} \\
GM-SIE$^1$ & \ref*{figGMSIE1}\\
GM-SIE$^2$ & \ref*{figGMSIE2}
\end{tabular}};

\end{tikzpicture}     \caption{Comparison between the NAG, (\ref{NAG_SIE_shi}) as Shi-SIE, SIE discretization of (\ref{GM-ODE}) as GM-SIE with superscript 1 and 2 when $n'=1,m'=\sqrt{s},q'=2\sqrt{\mu}$ and $n'=1-2\sqrt{\mu s},m'=\sqrt{s},q'=2\sqrt{\mu}$ respectively. The simulation function was $f(x)=4(L-\mu)\log(1+e^{-x})+\frac{\mu}{2}x^2$ with $L=1,\mu=0.01$. (a) The effect of the approximations ($1/(1-\sqrt{\mu s})\approx 1$ in Shi-SIE and coefficient deviation in GM-SIE$^1$) in the ODE trajectories. (b) Different coefficients used for discretizing (\ref{GM-ODE}). GM-SIE$^1$ is the SIE discretization of GM-ODE$^1$ (the recovered high-resolution NAG ODE from (\ref{GM-ODE})) and GM-SIE$^2$ is the SIE discretization of GM-ODE$^2$ (the ODE used to recover the NAG algorithm).}  
    \label{fig:three graphs}
\end{figure}

\subsection{Numerical Results}\label{app6}
In this section, numerical experiments are designed for further illustration of the previous findings. An important note is that (\ref{GM-ODE}) in \citep{zhang2021revisiting} uses different parameters to recover (\ref{NAG_ODE}) ($m'=\sqrt{s},q'=2\sqrt{\mu},n'=1$)  and the NAG algorithm after discretization ($m'=\sqrt{s},q'=2\sqrt{\mu},n'=1-2\sqrt{\mu s}$). In Figure \ref{fig:three graphs} we have considered two SIE discretizations of (\ref{GM-ODE}) and 
   \begin{align} \label{NAG_SIE_shi}
\left\{
    \begin{array}{ll}
    q_{k+1}-q_k = &  j_k\sqrt{s}, \\
   j_{k+1}-j_k  =& -2\sqrt{\mu s}j_{k+1}-\sqrt{s}(1+\sqrt{\mu s})\nabla f(q_{k+1})-\sqrt{s}(\nabla f(q_{k+1})-\nabla f(q_k)),
    \end{array}\right.\tag{Shi-SIE}
\end{align} 
which is the SIE discretization of (\ref{proposed_ODE_conti3}) used in \citep{Shi2021UnderstandingTA}. The discretizations of (\ref{GM-ODE}) are shown with GM-SIE$^1$ and GM-SIE$^2$. The aim of Figure \ref{fig:three graphs} is to highlight two things; First, the effect of \emph{coefficient inconsistency} (see Section \ref{sec:SIE_rec_NAG}) before and after discretization of (\ref{GM-ODE}) and second, to depict the approximation $1/(1-\sqrt{\mu s})\approx 1$ made in \citep{Shi2021UnderstandingTA} (see Section \ref{section4} for more detail). The step-size was $s=1/L$ in all simulations. All algorithms are simulated with the parameters they use to recover (\ref{NAG_ODE}) except for $\text{GM-SIE}^2$ which uses $n'=1-2\sqrt{\mu s}$ for the sake of comparison with $\text{GM-SIE}^1$. GM-SIE$^2$ does not exactly follow the NAG algorithm due to different initializations. We did not simulate  (\ref{SIE_newalg}) due to its exact match with the NAG method.\par

\end{document}